\documentclass[11pt]{article}
\usepackage[english]{babel}
\usepackage{indentfirst}
\usepackage{latexsym}
\usepackage[usenames,dvipsnames]{color}
\usepackage[colorlinks=true, bookmarksnumbered=true, bookmarksopen=true,
bookmarksopenlevel=3, pdfstartview=FitH, linkcolor=cyan, pdfmenubar=true,
pdftoolbar=true, bookmarks=true,citecolor=cyan, urlcolor=magenta,
filecolor=cyan,menucolor=black,plainpages=false,pdfpagelabels]{hyperref}
\usepackage[lmargin=2cm,tmargin=2cm,bmargin=2cm,rmargin=2cm]{geometry}
\usepackage{amsbsy,amsmath,amsfonts,amssymb,amsthm}
\usepackage{times}
\usepackage{graphicx}
\usepackage[utf8]{inputenc}
\usepackage[T1]{fontenc}
\usepackage{amsmath}
\usepackage{amsfonts}
\usepackage{amssymb}

\numberwithin{equation}{section}
\newtheorem{prop}{Proposition}[section]
\newtheorem{definition}[prop]{Definition}
\newtheorem{theorem}[prop]{Theorem}

\newtheorem{remark}[prop]{Remark}
\newtheorem{lemma}[prop]{Lemma}

\newtheorem{proposition}[prop]{Proposition}

\def\fin { \vskip 0pt \hfill $\diamond$ \vskip 12pt}

\begin{document}

\title{Local and global analysis in Besov-Morrey spaces for \\
inhomogeneous Navier-Stokes equations}
\author{{Lucas C. F. Ferreira$^{1}$}{\thanks{%
LCFF was partially supported by CNPq 308799/2019-4, BR. Email:
lcff@ime.unicamp.br (Corresponding author).}}, \ Daniel F. Machado$^{2}${%
\thanks{%
DFM was supported by CAPES (Finance Code 001), BR. Email:
daniellmath@gmail.com.}} \\
{\small $^{1,2}$ Universidade Estadual de Campinas, IMECC-Departamento de
Matem\'{a}tica,}\\
{\small Rua S\'{e}rgio Buarque de Holanda, CEP 13083-859, Campinas, SP,
Brazil.}}
\date{}
\maketitle

\begin{abstract}
In this paper we consider the incompressible inhomogeneous Navier-Stokes
equations in the whole space with dimension $n\geq 3$. We present local and
global well-posedness results in a new framework for inhomogeneous fluids,
namely Besov-Morrey spaces $\mathcal{N}_{p,q,r}^{s}$ that are Besov spaces
based on Morrey ones. In comparison with the previous works in Sobolev and
Besov spaces, our results provide a larger initial-data class for both the
velocity and density, constructing a unique global-in-time flow under
smallness conditions on weaker initial-data norms. In particular, we can
consider some kind of initial discontinuous densities, since our density
class $\mathcal{N}_{p,q,\infty }^{n/p}\cap L^{\infty }$ is not contained in
any space of continuous functions. From a technical viewpoint, the Morrey underlying norms prevent the common use of energy-type and integration by parts arguments, and then we need to obtain some estimates
for the localizations of the heat semigroup, the commutator, and the
volume-preserving map in our setting, as well as estimates for transport
equations and the linearized inhomogeneous Navier-Stokes system.

\medskip


\end{abstract}


\renewcommand{\abstractname}{Abstract}


\section{Introduction}

\bigskip We are concerned with the incompressible inhomogeneous
Navier-Stokes equations in the whole space
\begin{equation}
\left\{
\begin{array}{l}
\partial_{t}\rho+\mathrm{div}(\rho u)=0, \\
\partial_{t}(\rho u)+\mathrm{div}(\rho u\otimes u)-\mathrm{div}(\mu (\rho)%
\mathcal{M})+\nabla\pi=0, \\
\mathrm{div}\;u=0, \\
(\rho,u)_{|_{t=0}}=(\rho_{0},u_{0}),%
\end{array}
\right. \hspace{0.5cm}(x,t)\in\mathbb{R}^{n}\times\mathbb{R}^{+},
\label{sist:navier_stokes_0}
\end{equation}
where $n\geq3$ and the unknowns $\rho,$ $u=(u_{1},\dots,u_{n})$ and $\pi$
denote the density, velocity and scalar pressure of the fluid, respectively.
Moreover, the viscosity coefficient $\mu(\rho)$ is a positive smooth
function on $[0,\infty)$ and $\mathcal{M}=(\partial_{i}u_{j}+%
\partial_{j}u_{i})/2$ is the deformation tensor. This system describes the
evolution of incompressible fluids with non-constant density such as a fluid
containing a melted substance or mixing two miscible fluids with different
densities, as well as a multi-phase flow consisting of various
incompressible immiscible fluids with constant (different) densities and
viscosities. For more details on the physical model, see \cite{Lions-1}.

For simplicity, we use the same notation for scalar and vector function
spaces, that is, we write $u_{0}\in X$ instead of $u_{0}\in(X)^{n}$ for a
Banach space $X$. Also, we focus on the case $\mu(\rho)=1$. Although the
generalization to the variable viscosity case $\mu(\rho)\geq\mu>0$ with $\mu$
being a constant and $\mu(\cdot)$ a smooth function on $[0,\infty)$ is not
so immediate, we believe that it would be possible to treat it. In this way,
considering positive density and making the change $a=1/\rho-1$, system (\ref%
{sist:navier_stokes_0}) can be equivalently rewritten as
\begin{equation}
\left\{
\begin{array}{l}
\partial_{t}a+u\cdot\nabla a=0, \\
\partial_{t}u+u\cdot\nabla u+(1+a)(\nabla\pi-\Delta u)=0, \\
\mathrm{div}\;u=0, \\
(a,u)_{|_{t=0}}=(a_{0},u_{0}),%
\end{array}
\right. \hspace{0.5cm}(x,t)\in\mathbb{R}^{n}\times\mathbb{R}^{+}.
\label{sist:navier_stokes_1}
\end{equation}

Systems (\ref{sist:navier_stokes_0}) and (\ref{sist:navier_stokes_1}) have
been extensively studied in the framework of Sobolev and Besov spaces. In
what follows, without making a complete list, we recall some works of the
literature. In \cite{DiPerna-Lions, Lions-1}, DiPerna and Lions obtained
global existence of weak solutions in the Leray spirit for (\ref%
{sist:navier_stokes_0}) with $n\geq2$ (see also \cite{Simon1}). Considering $%
n=2,$ periodic conditions and $u_{0}\in H^{1}(\mathbb{T}^{2}),\rho_{0}\in
L^{\infty}(\mathbb{T}^{2})$ with $\mu(\rho_{0})\simeq c>0,$ Desjardins \cite%
{Desjardins-ARMA1997} showed some regularity properties of weak solutions,
say $u\in L^{\infty}([0,T];H^{1}(\mathbb{T}^{2}))$ and $\rho\in
L^{\infty}([0,T]\times\mathbb{T}^{2}).$ After, considering a modified model
of (\ref{sist:navier_stokes_0}) with $\mathrm{div}(\mu(\rho )\mathcal{M})$
replaced by $\nabla^{\perp}(\mu(\rho)\omega),$ the global smoothness was
proved by Zhang \cite{Zhang-2D-2008}, where $\omega$ is the vorticity, and
both terms coincide when the viscosity $\mu(\rho)$ is a constant (see also
\cite{Gui-Zhang}). In fact, the issues about existence, uniqueness and
regularity of solutions for (\ref{sist:navier_stokes_0}) and (\ref%
{sist:navier_stokes_1}) are far from completely settled for the full
dimension range $n\geq2.$

This scenario of development of the theory of weak and smooth solutions has
motivated to study the problem in distinct functional spaces, exploring
borderline cases of regularity to obtain the well-posedness property (i.e.,
existence, uniqueness, and persistence). In this direction, a local
well-posedness result for (\ref{sist:navier_stokes_1}) was obtained by
Danchin \cite{Danchin-localandglobal} for $(a_{0},u_{0})\in H^{n/2+\gamma
_{1}}(\mathbb{R}^{n})\times H^{n/2-1+\gamma _{2}}(\mathbb{R}^{n})$ with
sufficiently small $\left\Vert a_{0}\right\Vert _{H^{n/2+\gamma _{1}}}$,
where $\gamma _{1},\gamma _{2}>0$ with $\gamma _{2}\in (\gamma _{1}-1,\gamma
_{1}+1]$ and $\inf_{x}(1+a_{0}(x))=\underline{a}>0.$ Assuming in addition a
small condition on $\left\Vert u_{0}\right\Vert _{H^{n/2+\gamma _{2}}}$, he
obtained the global well-posedness. See also \cite{Danchin-inviscid} for
further results in Sobolev spaces and a study of inviscid limit properties
as $\mu (\rho )=\mu \rightarrow 0$. The results in \cite%
{Danchin-localandglobal} was extended to the context of homogeneous Besov
spaces in \cite{Danchin-density} by considering $(a_{0},u_{0})\in \dot{B}%
_{2,\infty }^{n/2}\cap L^{\infty }(\mathbb{R}^{n})\times \dot{B}%
_{2,1}^{-1+n/2}(\mathbb{R}^{n})$ and covering the critical regularities $%
(s_{1},s_{2})=(n/2,-1+n/2)$. After, considering small $a_{0},$ Abidi \cite%
{Abidi} showed local existence for (\ref{sist:navier_stokes_1}) in the
homogeneous Besov space $\dot{B}_{q,1}^{n/q}(\mathbb{R}^{n})\times \dot{B}%
_{q,1}^{n/q-1}(\mathbb{R}^{n})$ with $q\in (1,2n),$ and assuming $q\in (1,n]$
for the uniqueness property (and then for the well-posedness). The global
counterpart of the results is obtained by considering additionally a
smallness condition on $\left\Vert u_{0}\right\Vert _{\dot{B}%
_{q,1}^{n/q-1}}. $ Moreover, he proved a well-posedness result with an
improved regularity in Sobolev spaces $H^{n/2+\alpha }(\mathbb{R}^{n})\times
H^{n/2-1+\alpha }(\mathbb{R}^{n})$ for $\alpha >0$. The restriction $q\in
(1,n]$ for the velocity $u_{0}$ was improved in \cite{Abidi-existglobale} by
obtaining existence in the mixed space $\dot{B}_{l,1}^{n/l}(\mathbb{R}%
^{n})\times \dot{B}_{q,1}^{n/q-1}(\mathbb{R}^{n})$ with $l\in \lbrack
1,\infty )$ and $q\in (1,\infty )$ such that $\sup \{l^{-1},q^{-1}\}<n^{-1}
+\inf \{l^{-1},q^{-1}\}$ and assuming $l^{-1}+q^{-1}>n^{-1}$ for the
uniqueness property. By requiring small initial density in the multiplier
space of $\dot{B}_{q,1}^{n/q-1}$ and employing a Lagrangian approach,
Danchin and Mucha \cite{Danchin-Mucha} extended the uniqueness result in
\cite{Abidi} to the range $q\in (n,2n)$ and thus completing the
well-posedness result for $q\in (1,2n).$ Abidi, Gui and Zhang \cite%
{Abidi-onthewell} proved global well-posedness for $(a_{0},u_{0})\in
B_{2,1}^{3/2}(\mathbb{R}^{3})\times \dot{B}_{2,1}^{1/2}(\mathbb{R}^{3})$
with a smallness condition only on the norm of $u_{0}$ where such condition
depends on $\left\Vert a_{0}\right\Vert _{B_{2,1}^{3/2}}$ \ and $B_{q,r}^{s}$
stands for the inhomogeneous Besov spaces. In \cite{Abidi-wellposed}, they
considered (\ref{sist:navier_stokes_1}) with the initial data $(a_{0},u_{0})$
belonging to the mixed Besov space $B_{l,1}^{3/l}(\mathbb{R}^{3})\times \dot{%
B}_{q,1}^{3/q-1}(\mathbb{R}^{3})$ with $l\in \lbrack 1,2]$ and $q\in \lbrack
3,4]$ satisfying $l^{-1}+q^{-1}>5/6$ and $l^{-1}-q^{-1}\leq 1/3$, and proved
global well-posedness with the smallness condition $\left\Vert
u_{0}\right\Vert _{\dot{B}_{q,1}^{3/q-1}}\leq c$ with $c$ depending on the
size of $a_{0}$ in $B_{l,1}^{3/l}$ . First this result was improved by Zhai
and Yin \cite{Zhai-JDE2017} where the authors extended the range of $q$ to $%
(1,\frac{5+\sqrt{17}}{2})$ with $1<l\leq q,$ $l^{-1}+q^{-1}\geq 1/2$ and $%
l^{-1}-q^{-1}\leq 1/3.$ After, analyzing the inhomogeneous MHD system and in
particular system (\ref{sist:navier_stokes_1}), Huang and Qian \cite%
{Huang-Qian} complemented the previous ranges by considering the space $%
B_{l,1}^{3/l}(\mathbb{R}^{3})\times \dot{B}_{q,1}^{3/q-1}(\mathbb{R}^{3})$
with $q,l\in (1,6),$ $l^{-1}+q^{-1}>1/2$, $l^{-1}-q^{-1}<1/3,$ and $%
q^{-1}-l^{-1}<6^{-1}.$ In \cite{Burtea2017}, by using the Lagrangian
formulation as in \cite{Danchin-Mucha}, Burtea proved the local
well-posedness in the homogeneous space $\dot{B}_{q,1}^{3/q}(\mathbb{R}%
^{3})\times \dot{B}_{q,1}^{3/q-1}(\mathbb{R}^{3})$ for $q\in (6/5,4),$
without any additional condition on the index $q$, nor smallness assumptions
on initial data. As well as \cite{Danchin-localandglobal}, it is worth
noting that the above-mentioned works also assumed explicitly (or
implicitly) the basic condition $\inf_{x}(1+a_{0}(x))=\underline{a}>0.$
Also, energy-type or suitable integration by parts arguments for
frequency-localizations are commonly employed in order to obtain the key
estimates for (\ref{sist:navier_stokes_1}), except for the Lagrangian
approach in \cite{Burtea2017} and \cite{Danchin-Mucha}. In \cite%
{Danchin_2000}, Danchin considered the case of a barotropic compressible
fluid and obtained the global well-posedness in a critical framework based
on homogeneous Besov spaces $\dot{B}_{2,1}^{s}(\mathbb{R}^{n})$. Still in
the compressible case, a local theory for a general model of viscous and
heat-conductive gases was developed by him in \cite%
{Danchin-Compressible-CPDE-2001}.

In this paper we show local and global well-posedness for (\ref%
{sist:navier_stokes_1}) in a new framework, namely homogeneous Besov-Morrey
spaces $\mathcal{N}_{p,q,r}^{s}$ which are homogeneous Besov spaces based on
Morrey spaces $\mathcal{M}_{q}^{p}$. This setting allows us to consider a
larger initial-data class for both the velocity and density, providing a
unique global-in-time flow under a smallness condition on the weaker
initial-data norms. As a matter of fact, the homogeneous Besov-Morrey space $%
\mathcal{N}_{p,q,r}^{s}$ is strictly larger than the homogeneous Besov space
$\dot{B}_{p,r}^{s}$ for $q<p$, $1\leq r\leq \infty ,$ and $s\in \mathbb{R}$
(see Remark \ref{rem:1.3} $(i)$-$(ii)$ for further details). It is worth
mentioning that Besov-Morrey spaces were initially introduced in \cite%
{Kozono} to analyze Navier-Stokes equations via the so-called Kato approach
(see also \cite{Mazzucato}). Moreover, other fluid dynamics models have been
studied in those spaces, see, e.g., \cite{Bie-inviscid, DFV2020,
Fer-Perez-1, Jiang-morrey} and their references.

For $0<T\leq \infty $, we denote%
\begin{equation}
\widetilde{C}([0,T),X)=C([0,T),X)\cap \widetilde{L}_{T}^{\infty }(X),
\label{class-1}
\end{equation}%
where $X$ is a Banach space and $\widetilde{L}_{T}^{\infty }(X)$ is the
Chemin-Lerner type space (see, e.g., (\ref{chemin-lerner}) more below). In
fact, due to the space $\widetilde{L}_{T}^{\infty }(X),$ we should consider $%
X$ in (\ref{class-1}) as some type of Besov space such as $\mathcal{N}%
_{p,q,r}^{s}.$ Moreover, the time-continuity at $t=0^{+}$ in (\ref{class-1})
should be meant in the sense of tempered distributions when the summation
index $r$ is equal to $\infty $.

Our main result reads as follows.

\begin{theorem}
\label{the:navier_stokes_2} Let $n\geq 3$, $1<q\leq p<\infty $ and $1\leq
r\leq \infty $. Assume either {$n/p-1<s\leq n/p$ with $n/p\geq 1$ or $%
s=n/p-1 $ with $n/p>1$}. Consider $u_{0}\in \mathcal{N}_{p,q,1}^{s}$ with $%
\mathrm{div}\hspace{0.05cm}u_{0}=0$ and $a_{0}\in \mathcal{N}%
_{p,q,r}^{n/p}\cap L^{\infty }$. There exist $T\in (0,\infty ]$ and a small
constant $c>0$ such that if
\begin{equation*}
\Vert a_{0}\Vert _{\mathcal{N}_{p,q,\infty }^{n/p}\cap L^{\infty }}\leq c,
\end{equation*}%
then system (\ref{sist:navier_stokes_1}) has a unique solution $(a,u,\nabla
\pi )$ satisfying%
\begin{equation}
a\in \widetilde{C}([0,T);\mathcal{N}_{p,q,r}^{n/p})\cap L_{T}^{\infty
}(L^{\infty }),\hspace{0.25cm}u\in \widetilde{C}([0,T);\mathcal{N}%
_{p,q,1}^{s})\cap L_{T}^{1}(\mathcal{N}_{p,q,1}^{s+2})\text{ \ and \ }\nabla
\pi \in L_{T}^{1}(\mathcal{N}_{p,q,1}^{s}),  \label{class-sol-1}
\end{equation}%
Moreover, in the critical case $s=n/p-1$, there exists a small constant $%
c^{\prime }>0$ such that if $\Vert u_{0}\Vert _{\mathcal{N}_{p,q,1}^{s}}\leq
c^{\prime }$, then $T=\infty $.
\end{theorem}

Further comments on the above result are in order.

\begin{remark}
\label{rem:1.3}

\begin{itemize}
\item[(i)] (Velocity class) As already pointed out more above, the strict
inclusion $\dot{B}_{p,r}^{s}\hookrightarrow \mathcal{N}_{p,q,r}^{s}$ holds
for $q<p,$ $r\in \lbrack 1,\infty ],$ and $s\in \mathbb{R}$. Then, in
comparison with previous works, we are providing in Theorem \ref%
{the:navier_stokes_2} a larger initial velocity class when $s>0$ $(\,$i.e., $%
1<p<n)$. On the other hand, considering jointly the results of \cite{Abidi},
\cite{Danchin-density}, and \cite{Danchin-Mucha}, we have the
existence-uniqueness for (\ref{sist:navier_stokes_1}) in $\dot{B}%
_{q,1}^{n/q}(\mathbb{R}^{n})\times \dot{B}_{q,1}^{n/q-1}(\mathbb{R}^{n})$
with $q\in (1,2n)$ which covers some negative regularity indexes, that is, $%
s=n/q-1\in (-1/2,0).$ Although we treat with $s>0,$ our initial-data class $%
\mathcal{N}_{p,q,1}^{s}$ with $1<q<p<\infty $ cannot be embedded in any
homogeneous Besov space $\dot{B}_{l,r}^{\tilde{s}}$ with negative regularity
$\tilde{s}$. Even more, we have that $\mathcal{N}_{p,q,1}^{s}\not\subset
\dot{B}_{l,r}^{\tilde{s}}$ for any $\tilde{s}\in \mathbb{R}$, $l\in \lbrack
1,\infty )$ and $r\in \lbrack 1,\infty ]$.

\item[(ii)] (Density class) Concerning the density, in view of the strict
inclusions $\dot{B}_{p,r}^{n/p}\hookrightarrow \mathcal{N}_{p,q,r
}^{n/p}\hookrightarrow\mathcal{N}_{p,q,\infty}^{n/p}$ for all $r\in \lbrack 1,\infty ]$, our initial data class extends those of \cite{Danchin-density} and \cite{Abidi} (and some other works mentioned above), where the authors
considered $a_{0}\in \dot{B}_{2,\infty }^{n/2}\cap L^{\infty }$ and $%
a_{0}\in \dot{B}_{p,1}^{n/p},$ respectively. Moreover, we have that $%
\mathcal{N}_{p,q,\infty }^{n/p}\subset \dot{B}_{\infty ,\infty }^{0}$
(homogeneous Zygmund spaces) but $\mathcal{N}_{p,q,\infty }^{n/p}$ is not
contained in any space of continuous functions. In this way, we can consider
discontinuous initial densities $a_{0}$. In \cite{Danchin-Mucha}, Danchin
and Mucha considered the initial density belonging to the multiplier space
of $\dot{B}_{q,1}^{n/q-1},$ denoted by $\mathcal{M}(\dot{B}_{q,1}^{n/q-1}),$
where $q\in \lbrack 1,2n)$. For $l^{-1}>n^{-1}-p^{-1}$ and $l^{-1}\geq
p^{-1}-n^{-1},$ we have the embedding $\dot{B}_{l,\infty }^{n/l}\cap
L^{\infty }\hookrightarrow \mathcal{M}(\dot{B}_{q,1}^{n/q-1}).$ However, as
far as we know, there is no inclusion relations between $\mathcal{M}(\dot{B}%
_{q,1}^{n/q-1})$ and $\mathcal{N}_{p,q,\infty }^{n/p}\cap L^{\infty }.$
Finally, we highlight that a smallness condition in the weaker norm of $%
\mathcal{N}_{p,q,\infty }^{s}$ allows us to consider some classes of initial
data that could be large in other setting such as Sobolev $\dot{H}_{p}^{s}$ and
Besov $\dot{B}_{p,\infty }^{s}$ spaces.

\item[(iii)] (Forcing terms) We can also treat (\ref{sist:navier_stokes_0})
with a coupling term $\rho f$ on the right-hand side of the velocity
equation. For that, it is sufficient to assume $f\in L_{loc}^{1}([0,\infty );%
\mathcal{N}_{p,q,1}^{s})$ for the local part in Theorem \ref%
{the:navier_stokes_2}; for the global one, it is necessary to additionally
assume a smallness condition on $\Vert f\Vert _{L^{1}(\mathbb{R}^{+};%
\mathcal{N}_{p,q,1}^{s})}.$
\end{itemize}
\end{remark}

In order to obtain our results we need to perform an adaptation from the
approach in \cite{Abidi,Abidi-wellposed,Danchin-density, Danchin_2000,
Danchin-inviscid, Danchin-Compressible-CPDE-2001} which involves some
different ingredients and certain care. For the reader convenience, we
highlight some points in the proofs and technical difficulties that we
overcome. The $\mathcal{M}_{q}^{p}$-norm involves computations of $L^{q}$%
-norms in closed balls $B(x_{0},R)$ of $\mathbb{R}^{n}$ and a corresponding
weighted supremum over $R>0$ (see Definition \ref{def:morrey_spaces}), which
make it difficult to employ energy-type and integration by parts arguments
such as in previous results in Sobolev and Besov spaces (see, e.g., \cite%
{Abidi, Abidi-wellposed,Abidi-existglobale,Danchin-density, Danchin_2000,
Danchin-localandglobal, Danchin-Compressible-CPDE-2001,
Huang-Qian,Zhai-JDE2017} and their references). In fact, much of the $L^{q}$%
-theory in the whole space $\mathbb{R}^{n}$ does not work in a
straightforward way in the context of Morrey and Besov-Morrey spaces. Thus,
in order to avoid those difficulties and to develop suitable linear
estimates, we handle the equations of the velocity $u$ and $\nabla\pi$ by
performing localized heat-semigroup estimates in our setting, inspired by
some arguments found in \cite{Chemin_2}{, and using the Leray projector }$%
\mathbb{P}$ and some Bernstein-type inequalities (see Lemmas \ref%
{lem:exp_h_2} and \ref{cor:exp_heat_2} and Propositions \ref{prop:stokes}
and \ref{prop:navier_stokes}). In the estimates for $a$ we employ the
volume-preserving map $X$ associated with $u,$ which leads to estimate the
map $X$ in Morrey spaces, and thus preventing increasing the level of
regularity in $a$, in a spirit resembling the Lagrangian approach in \cite%
{Danchin-Mucha} (see Lemma \ref{lem:X_general} and Proposition \ref%
{prop:transport}). Moreover, by adapting the proofs of some previous works
(see, e.g., \cite{Bie-inviscid, Commutator-interms}), we derive some
commutator estimates since we are not able to locate in the literature
estimates in a suitable form for our purposes (see Lemmas \ref%
{lem:commutator_uv} and \ref{lem:commutator_pi}).

\textbf{Organization of the paper:} Section 2 is intended to provide some
useful notations and review definitions and properties about functional
spaces such as Morrey, Besov-Morrey and Chemin-Lerner type spaces. In
Section 3 we present some estimates for the commutator, the
volume-preserving map and the localizations of the heat semigroup in our
framework. Section 4 is devoted to estimates for linear systems associated
with (\ref{sist:navier_stokes_1}). In Section 5 we give the proof of Theorem %
\ref{the:navier_stokes_2} which is based on uniform estimates in our setting
for a hyperbolic-parabolic approximate linear problem.

\section{Preliminaries}

This section is dedicated to collecting basic notations, definitions, tools
and properties about some operators and functional spaces which will be
useful for our ends.

We use the notation $A\lesssim B$ to mean that there exists an absolute
constant $C>0$, which may change from line to line, such that $A\leq CB.$
The commutator between two operators $F_{1}$ and $F_{2}$ is denoted by $%
[F_{1},F_{2}]=F_{1}F_{2}-F_{2}F_{1}$. Also, $C_{0}^{\infty }(\mathbb{R}%
^{n}), $ $\mathcal{S}=\mathcal{S}(\mathbb{R}^{n})$, and $\mathcal{S}^{\prime
}=\mathcal{S}^{\prime }(\mathbb{R}^{n})$ stand for the space of compactly
supported smooth functions, Schwartz space, and the space of tempered
distributions, respectively. The notation $\mathbb{N}_{0}=\{0\}\cup \mathbb{N%
}$ represents the set of non-negative integers.

We start by giving some preliminaries on Morrey spaces. For more details,
see \cite{Kozono, Rosenthal, Sawano-morrey,Zhou-morrey}.

\begin{definition}
\label{def:morrey_spaces} Let $1\leq q\leq p<\infty$. The Morrey space $%
\mathcal{M}_{q}^{p}=\mathcal{M}_{q}^{p}(\mathbb{R}^{n})$ is the set of all $%
u\in L_{\mathrm{loc}}^{q}(\mathbb{R}^{n})$ satisfying
\begin{equation*}
\Vert u\Vert_{\mathcal{M}_{q}^{p}}:=\sup_{x_{0}\in\mathbb{R}^{n}}\sup
_{R>0}R^{n/p-n/q}\left( \int_{B(x_{0},R)}|u(y)|^{q}\;dy\right) ^{1/q}<\infty,
\end{equation*}
where $B(x_{0},R)\subset\mathbb{R}^{n}$ stands for the closed ball with
radius $R>0$ and center $x_{0}$.
\end{definition}

\begin{remark}
It is not difficult to see that $\mathcal{M}_{q}^{p}$ is a Banach space with
the norm $\Vert \cdot \Vert _{\mathcal{M}_{q}^{p}}$. Furthermore, we have
the relation $\mathcal{M}_{r}^{p}\subset \mathcal{M}_{q}^{p}$ for $1\leq
q<r\leq p<\infty $ and that $\mathcal{M}_{p}^{p}=L^{p}(\mathbb{R}^{n})$. The
case $p=\infty $ could be identify with $L^{\infty }(\mathbb{R}^{n})$, i.e.,
$\mathcal{M}_{q}^{\infty }=L^{\infty }.$
\end{remark}

In what follows, we recall the Littlewood-Paley decompositions, see \cite%
{Cannone, Chemin, Lemarie} for further details. For that, consider $\psi $
and $\varphi \in C_{0}^{\infty }(\mathbb{R}^{n})$ functions supported in the
ball $\mathcal{B}:=\{\xi \in \mathbb{R}^{n}:|\xi |\leq 4/3\}$ and in the
ring $\mathcal{C}:=\{\xi \in \mathbb{R}^{n}:3/4\leq |\xi |\leq 8/3\},$
respectively, satisfying
\begin{equation}
\sum_{j\in \mathbb{Z}}\varphi _{j}(\xi )=1,\hspace{0.25cm}%
\mbox{for
$\xi\in\mathbb{R}^n\setminus\{0\}$}\hspace{0.75cm}\text{and}\hspace{0.75cm}%
\psi (\xi )+\sum_{j\in \mathbb{N}_{0}}\varphi _{j}(\xi )=1,\hspace{0.25cm}%
\mbox{for $\xi\in\mathbb{R}^n$}.  \label{decomp:unitaria}
\end{equation}%
where $\varphi _{j}(\xi ):=\varphi (2^{-j}\xi )$. Now, for $u\in \mathcal{S}%
^{\prime }(\mathbb{R}^{n})$, we define the frequency projection by
\begin{equation*}
\Delta _{j}u:=[\varphi _{j}\widehat{u}]^{\vee },\hspace{0.25cm}%
\mbox{for all
$j\in\mathbb{Z}$}\hspace{0.5cm}\mbox{and}\hspace{0.5cm}\bar{\Delta}%
_{j}u:=\left\{
\begin{array}{cl}
0, & j\leq -2, \\
\left[ \psi \hspace{0.1cm}\widehat{u}\right] ^{\vee }, & j=-1, \\
\Delta _{j}u, & j\geq 0,%
\end{array}%
\right.
\end{equation*}%
and also the low-frequency cut-off operators $S_{j}u:=\sum_{l\leq j-1}\Delta
_{l}u$ and $\tilde{S}_{j}u:=\left[ \psi (2^{-j}\xi )\hspace{0.1cm}\widehat{u}%
\right] ^{\vee },$ for $j\in \mathbb{Z}$, where $\widehat{u}$ denotes the
Fourier transform of $u$ on $\mathbb{R}^{n}$ and $u^{\vee }$ denotes the
inverse Fourier transform. Denote by $\mathcal{L}_{h}^{\prime }=\mathcal{L}%
_{h}^{\prime }(\mathbb{R}^{n})$ and $\mathcal{S}_{h}^{\prime }=\mathcal{S}%
_{h}^{\prime }(\mathbb{R}^{n})$ the spaces of all $u\in \mathcal{S}^{\prime
} $ such that $\lim_{j\rightarrow -\infty }\displaystyle\tilde{S}_{j}u=0$ in
$\mathcal{S}^{\prime }$ and $\lim_{j\rightarrow -\infty }||\tilde{S}%
_{j}u||_{L^{\infty }}=0$, respectively. Then, we have the Littlewood-Paley
decompositions
\begin{equation*}
u=\sum_{j\in \mathbb{Z}}\Delta _{j}u,\hspace{0.25cm}\text{for all }u\in
\mathcal{L}_{h}^{\prime }\hspace{0.75cm}\mbox{and}\hspace{0.75cm}%
u=\sum_{j\in \mathbb{Z}}\bar{\Delta}_{j}u,\text{ \ for all }u\in \mathcal{S}%
^{\prime }.
\end{equation*}%
Moreover, the above decompositions have the properties
\begin{equation*}
\Delta _{j}(\Delta _{k}u)\equiv 0\hspace{0.25cm}\mbox{if $|j-k|\geq2$}%
\hspace{0.5cm}\mbox{and}\hspace{0.5cm}\Delta _{j}(S_{k-1}u\Delta
_{k}u)\equiv 0\hspace{0.25cm}\mbox{if $|j-k|\geq5$}.
\end{equation*}

Now we are in position to recall the definition of homogeneous and
inhomogeneous Besov-Morrey spaces as well as some of their basic properties
(see \cite{Bie-inviscid, Commutator-interms, Kozono, Mazzucato}).

\begin{definition}
\label{def:besov_morrey_homog} Consider $s\in \mathbb{R}$, $1\leq q\leq
p<\infty$ and $1\leq r\leq \infty .$ The homogeneous Besov-Morrey space $%
\mathcal{N}_{p,q,r}^{s}=\mathcal{N}_{p,q,r}^{s}(\mathbb{R}^{n})$ is the set
of all $u\in \mathcal{S}_{h}^{\prime }$ such that $\Delta _{j}u\in \mathcal{M%
}_{q}^{p}$ for every $j\in \mathbb{Z}$, and that
\begin{equation*}
\Vert u\Vert _{\mathcal{N}_{p,q,r}^{s}}:=\left\{
\begin{array}{l}
\left( \displaystyle\sum_{j\in \mathbb{Z}}2^{sjr}\Vert \Delta _{j}u\Vert _{%
\mathcal{M}_{q}^{p}}^{r}\right) ^{1/r}<\infty ,\hspace{0.5cm}%
\mbox{if
$r<\infty$}, \\
\sup_{j\in \mathbb{Z}}2^{sj}\Vert \Delta _{j}u\Vert _{\mathcal{M}%
_{q}^{p}}<\infty ,\hspace{0.5cm}\mbox{if $r=\infty$}.%
\end{array}%
\right.
\end{equation*}%
The inhomogeneous Besov-Morrey space $N_{p,q,r}^{s}=N_{p,q,r}^{s}(\mathbb{R}%
^{n})$ is the set of all $u\in \mathcal{S}^{\prime }$ such that $\bar{\Delta}%
_{j}u\in \mathcal{M}_{q}^{p}$ for every $j\geq -1$ and that
\begin{equation*}
\Vert u\Vert _{N_{p,q,r}^{s}}:=\left\{
\begin{array}{l}
\left( \displaystyle\sum_{j\geq -1}2^{sjr}\Vert \bar{\Delta}_{j}u\Vert _{%
\mathcal{M}_{q}^{p}}^{r}\right) ^{1/r}<\infty ,\hspace{0.5cm}%
\mbox{if
$r<\infty$}, \\
\sup_{j\geq -1}2^{sj}\Vert \bar{\Delta}_{j}u\Vert _{\mathcal{M}%
_{q}^{p}}<\infty ,\hspace{0.5cm}\mbox{if $r=\infty$}.%
\end{array}%
\right.
\end{equation*}
\end{definition}

The spaces $\mathcal{N}_{p,q,r}^{s}$ and $N_{p,q,r}^{s}$ are Banach spaces
and the continuous inclusions $\mathcal{N}_{p,q,r}^{s}\subset \mathcal{S}%
_{h}^{\prime }$ and $N_{p,q,r}^{s}\subset \mathcal{S}^{\prime }$ holds. From
the inclusion relations in $\mathcal{M}_{q}^{p}$ and $\ell ^{r}(\mathbb{Z})$%
, it follows that
\begin{equation}
N_{p,q,r_{1}}^{s}\subset N_{p,q,r_{2}}^{s}\hspace{0.5cm}\text{and}\hspace{%
0.5cm}\mathcal{N}_{p,q,r_{1}}^{s}\subset \mathcal{N}_{p,q,r_{2}}^{s},
\label{inclusion-1}
\end{equation}%
for $1\leq r_{1}\leq r_{2}\leq \infty $. Furthermore, if $s_{1}>s_{2}$ we
also have the following inclusion relation $N_{p,q,r_{1}}^{s_{1}}\subset
N_{p,q,r_{2}}^{s_{2}}$ for all $1\leq r_{1},r_{2}\leq \infty $. Also, for $%
s>0$ we have that $N_{p,q,r}^{s}=\mathcal{N}_{p,q,r}^{s}\cap \mathcal{M}%
_{q}^{p}$ and
\begin{equation}
\Vert u\Vert _{N_{p,q,r}^{s}}\simeq \Vert u\Vert _{\mathcal{N}%
_{p,q,r}^{s}}+\Vert u\Vert _{\mathcal{M}_{q}^{p}}.
\label{def:new_norm_inhomog}
\end{equation}%
Using the decompositions in (\ref{decomp:unitaria}), one can show the
following embeddings (see, e.g., \cite{Mazzucato}).

\begin{remark}
Let $1\leq q\leq p<\infty$. For either $s>n/p$ with $1\leq r\leq\infty$ or $%
s=n/p$ with $r=1$, it follows that
\begin{equation}
N_{p,q,r}^{s}\hookrightarrow L^{\infty}\hspace{0.5cm}\mbox{and}\hspace {0.5cm%
}\mathcal{N}_{p,q,1}^{n/p}\hookrightarrow L^{\infty}.
\label{linfty_imersion}
\end{equation}
\end{remark}

As the study of system (\ref{sist:navier_stokes_1}) involves the
time-variable, we need to work with Chemin-Lerner type spaces.

\begin{definition}
Let $s\in \mathbb{R}$, $1\leq q\leq p<\infty $, $1\leq r,\beta \leq \infty $%
, and $0<T\leq \infty $. We define the space $\widetilde{L}_{T}^{\beta }(%
\mathcal{N}_{p,q,r}^{s})$ as the set of all $u\in \mathcal{S}^{\prime }(%
\mathbb{R}^{n}\times (0,T))$ such that $\lim_{j\rightarrow -\infty }%
\displaystyle\tilde{S}_{j}u=0$ in $L^{\beta }((0,T);L^{\infty }(\mathbb{R}%
^{n}))$ and the quantity
\begin{equation}
\Vert u\Vert _{\widetilde{L}_{T}^{\beta }(\mathcal{N}_{p,q,r}^{s})}:=\left(
\sum_{j\in \mathbb{Z}}2^{sjr}\left( \int_{0}^{T}\Vert \Delta _{j}u(t)\Vert _{%
\mathcal{M}_{q}^{p}}^{\beta }\;dt\right) ^{r/\beta }\right) ^{1/r}<\infty ,
\label{chemin-lerner}
\end{equation}%
with the natural change if $r=\infty $. The space $\widetilde{L}_{T}^{\beta
}(\mathcal{N}_{p,q,r}^{s})$ equipped with \eqref{chemin-lerner} is a Banach
space.
\end{definition}

\begin{remark}
\label{rem:2.6} F{or $\theta\in(0,1)$, it is not difficult to see that
\begin{equation}
\Vert u\Vert_{\widetilde{L}_{T}^{\beta}(\mathcal{N}_{p,q,r}^{s})}\leq\Vert
u\Vert_{\widetilde{L}_{T}^{\beta_{1}}(\mathcal{N}_{p,q,r}^{s_{1}})}^{\theta
}\Vert u\Vert_{\widetilde{L}_{T}^{\beta_{2}}(\mathcal{N}%
_{p,q,r}^{s_{2}})}^{1-\theta},  \label{est:interpolation}
\end{equation}
with $1/\beta=\theta/\beta_{1}+(1-\theta)/\beta_{2}$ and $s=\theta
s_{1}+(1-\theta)s_{2}$.} Moreover, Minkowski inequality yields
\begin{equation}
\Vert u\Vert_{\widetilde{L}_{T}^{\beta}(\mathcal{N}_{p,q,r}^{s})}\leq\Vert
u\Vert_{L_{T}^{\beta}(\mathcal{N}_{p,q,r}^{s})},\;\text{ if }\beta\leq r,%
\hspace{0.5cm} \text{and}\hspace{0.5cm} \Vert u\Vert_{L_{T}^{\beta}(\mathcal{%
N}_{p,q,r}^{s})}\leq\Vert u\Vert_{\widetilde{L}_{T}^{\beta}(\mathcal{N}%
_{p,q,r}^{s})},\;\text{ if }r\leq\beta.  \label{est:l_tilde}
\end{equation}
The same inequalities also hold in inhomogeneous Besov-Morrey spaces.
\end{remark}

We finish the section by presenting the Bony paraproduct decomposition. For
simplicity, we only do so in the homogeneous case.

Let $u,v\in\mathcal{S}^{\prime}$. The product of $u$ with $v$ admits the
following decomposition
\begin{equation}
uv=T_{u}v+T_{v}u+\mathcal{R}(u,v)=T_{u}v+R(u,v),  \label{decomp_para_product}
\end{equation}
where
\begin{equation*}
\begin{array}{rclcrcl}
\displaystyle T_{u}v & := & \displaystyle\sum_{j\in\mathbb{Z}%
}S_{j-1}u\Delta_{j}v, & \hspace{1cm} & R(u,v) & := & \displaystyle\sum _{j\in%
\mathbb{Z}}\Delta_{j}uS_{j+2}v, \\
\mathcal{R}(u,v) & := & \displaystyle\sum_{j\in\mathbb{Z}}\Delta_{j}u%
\widetilde{\Delta }_{j}v, & \hspace{1cm} & \widetilde{\Delta}_{j}v & := & %
\displaystyle\sum _{i=-1}^{1}\Delta_{j-i}v%
\end{array}%
\end{equation*}

\subsection{Some basic estimates in Morrey spaces}

We recall that H\"{o}lder-type inequalities work well in Morrey spaces (see,
e.g., {\cite{Kozono}}).

\begin{lemma}
Let $N\in\mathbb{N}$ and suppose that $p,$ $q,$ $\{p_{i}\}_{i=1}^{N}$ and $%
\{q_{i}\}_{i=1}^{N}$ satisfy $1\leq q_{i}\leq p_{i}<\infty$ for every $i =
1,\cdots,N$, $\sum_{i=1}^{N}1/q_{i}\leq1/q\leq1$ and $1/p=%
\sum_{i=1}^{N}1/p_{i},$ where $q\leq p.$ Then we have following:

\begin{itemize}
\item[(i)] If $u_{i}\in \mathcal{M}_{q_{i}}^{p_{i}}$ for all $i=1,\cdots ,N$
and also $u_{0}\in L^{\infty }(\mathbb{R}^{n})$, we have
\begin{equation}
\Vert u_{0}\cdot u_{1}\cdots u_{N}\Vert _{\mathcal{M}_{q}^{p}}\leq C\Vert
u_{0}\Vert _{L^{\infty }(\mathbb{R}^{n})}\Vert u_{1}\Vert _{\mathcal{M}%
_{q_{1}}^{p_{1}}}\cdots \Vert u_{N}\Vert _{\mathcal{M}_{q_{N}}^{p_{N}}}.
\label{est:holder_geral}
\end{equation}

\item[(ii)] In the case $N=1$ and $q_{1}=q$, we have
\begin{equation}
\Vert u_{0}\cdot u_{1}\Vert _{\mathcal{M}_{q}^{p}}\leq C\Vert u_{0}\Vert
_{L^{\infty }(\mathbb{R}^{n})}\Vert u_{1}\Vert _{\mathcal{M}_{q}^{p}}.
\label{est:holder_linfty}
\end{equation}
\end{itemize}
\end{lemma}

The following lemma collects some Bernstein-type inequalities in the
framework of Morrey spaces found in the literature (see {\cite{Bie-inviscid,
Chemin, Kozono, Jiang-morrey}}). For $0<R_{1}<R_{2},$ consider the notations
$\mathcal{B}(0,R_{1})=\left\{ \xi\in\mathbb{R}^{n};\left\vert \xi\right\vert
\leq R_{2}\right\} $ and $\mathcal{C}(0,R_{1},R_{2})=\left\{ \xi \in\mathbb{R%
}^{n};R_{1}\leq\left\vert \xi\right\vert \leq R_{2}\right\} $.

\begin{lemma}
\label{lem:bernstein} Let $1\leq q\leq p<\infty$, $k\in\mathbb{N}$ and $j\in%
\mathbb{Z}$.

\begin{itemize}
\item[(i)] For $u\in \mathcal{M}_{q}^{p}$ satisfying $\mathrm{supp}(\widehat{%
u})\subset \mathcal{B}(0,\lambda R_{1})$ for some $\lambda >0$ and $R_{1}>0$%
, it follows that
\begin{equation}
\sup_{|\alpha |=k}\Vert \partial ^{\alpha }u\Vert _{\mathcal{M}_{q}^{p}}\leq
C^{k+1}\lambda ^{k}\Vert u\Vert _{\mathcal{M}_{q}^{p}},
\label{bernstein_bola}
\end{equation}%
where the constant $C:=C(n,p,q,R_{1},R_{2})$ is independent of $\lambda $.

\item[(ii)] For $u\in \mathcal{M}_{q}^{p}$ satisfying $\mathrm{supp}(%
\widehat{u})\subset \mathcal{C}(0,\lambda R_{1},\lambda R_{2})$ for some $%
\lambda >0$ and $0<R_{1}<R_{2}$, it follows that
\begin{equation}
C^{-(k+1)}\lambda ^{k}\Vert u\Vert _{\mathcal{M}_{q}^{p}}\leq \sup_{|\alpha
|=k}\Vert \partial ^{\alpha }u\Vert _{\mathcal{M}_{q}^{p}}\leq
C^{k+1}\lambda ^{k}\Vert u\Vert _{\mathcal{M}_{q}^{p}},
\label{bernstein_anel}
\end{equation}%
where the constant $C:=C(n,p,q,R_{1},R_{2})$ is independent of $\lambda $.

\item[(iii)] For $u\in \mathcal{M}_{q}^{p}$ satisfying $\mathrm{supp}(%
\widehat{u})\subset \mathcal{B}(0,\lambda ^{j}R_{1})$ for some $\lambda >0$
and $R_{1}>0$, it follows that
\begin{equation}
\Vert u\Vert _{L^{\infty }}\leq C\lambda ^{jn/p}\Vert u\Vert _{\mathcal{M}%
_{q}^{p}}.  \label{est:linfty_morrey}
\end{equation}
\end{itemize}
\end{lemma}

\begin{remark}
As a consequence of the lemma above, we have the estimates
\begin{align}
\sup_{|\alpha|=k}\Vert\partial^{\alpha}u\Vert_{N_{p,q,r}^{s}} & \leq
C^{k+1}\Vert u\Vert_{N_{p,q,r}^{s+k}},  \label{bernstein_bm_inhomog} \\
C^{-(k+1)}\Vert u\Vert_{\mathcal{N}_{p,q,r}^{s+k}}\leq\sup_{|\alpha|=k}\Vert
& \partial^{\alpha}u\Vert_{\mathcal{N}_{p,q,r}^{s}}\leq C^{k+1}\Vert u\Vert_{%
\mathcal{N}_{p,q,r}^{s+k}},  \label{bernstein_bm_homog}
\end{align}
where $1\leq q\leq p<\infty$, $1\leq r\leq\infty$, $k\in\mathbb{N},$ $s\in%
\mathbb{R}$, and $C>0$ is a universal constant.
\end{remark}

\section{Commutator, heat estimates and volume-preserving maps}

This section is devoted to present some commutator and heat{\ estimates} in
the context of Morrey and Besov-Morrey spaces, as well as estimates for
volume-preserving maps.

We start with \ a commutator estimate which we have not been able to locate
them in the literature with the needed hypotheses and conclusions for our
purposes. The reader is referred to \cite{Bie-inviscid, Commutator-interms}
for similar estimates in Besov-Morrey spaces.

\begin{lemma}
\label{lem:commutator_uv} For $0<s<n/p+1$ with $1\leq r\leq\infty$, or $%
s=n/p+1$ with $r=1$, and $1\leq q\leq p<\infty$, there is a constant $C>0$
such that
\begin{equation}
\Vert2^{sj}\Vert\lbrack\Delta_{j},v\cdot\nabla]u\Vert_{\mathcal{M}%
_{q}^{p}}\Vert_{\ell^{r}}\leq C\Vert u\Vert_{\mathcal{N}_{p,q,r}^{s}}\left(%
\Vert v\Vert_{\mathcal{N}_{p,q,\infty}^{n/p+1}} + \Vert\nabla
v\Vert_{L^{\infty}}\right) ,  \label{est:commutator_uv}
\end{equation}
for all $u\in\mathcal{N}_{p,q,r}^{s}$ and $v\in\mathcal{N}%
_{p,q,\infty}^{n/p+1}$ with $\nabla v\in L^{\infty}$ and $\mathrm{div}\;v=0$%
. As a consequence, we have that
\begin{equation}
\Vert2^{sj}\Vert\lbrack\Delta_{j},v\cdot\nabla]u\Vert_{\mathcal{M}%
_{q}^{p}}\Vert_{\ell^{r}}\leq C\Vert u\Vert_{\mathcal{N}_{p,q,r}^{s}} \Vert
v\Vert_{\mathcal{N}_{p,q,1}^{n/p+1}},  \label{est:commutator_uv_2}
\end{equation}
for $u\in\mathcal{N}_{p,q,r}^{s}$ and $v\in \mathcal{N}_{p,q,1}^{n/p+1}$.
\end{lemma}

\textbf{Proof.} By the Bony decomposition (\ref{decomp_para_product}) and
using $\nabla\cdot v=0$, we can write
\begin{equation}
\left[ \Delta_{j},v\cdot\nabla\right] u=\underbrace{\mathrm{div}(\Delta _{j}(%
\mathcal{R}(u,v)))}_{=:\mathcal{R}_{j}^{1}}+\underbrace{\Delta_{j}\left(
T_{\nabla u}v\right) }_{=:\mathcal{R}_{j}^{2}}\;\underbrace{-\;R\left(
v,\Delta_{j}\nabla u\right) }_{=:\mathcal{R}_{j}^{3}}\;\underbrace{-\;\left[
T_{v},\Delta_{j}\right] \nabla u}_{=:\mathcal{R}_{j}^{4}}.
\label{est:comt_decomp}
\end{equation}
For $\mathcal{R}_{j}^{1}$, it follows that
\begin{equation*}
\mathcal{R}_{j}^{1}=\sum_{j-k\leq3}\mathrm{div}(\Delta_{j}(\Delta_{k}u\tilde{%
\Delta}_{k}v))=\sum_{j-k\leq3}\sum_{i=-1}^{1}\mathrm{div}(\Delta
_{j}(\Delta_{k}u\Delta_{k-i}v)).
\end{equation*}
Then, by Bernstein inequality (\ref{bernstein_anel}) and H\"{o}lder
inequality (\ref{est:holder_linfty}), we have that
\begin{align*}
\Vert\mathcal{R}_{j}^{1}\Vert_{\mathcal{M}_{q}^{p}} & \lesssim\sum
_{j-k\leq3}\sum_{i=-1}^{1}2^{j}\Vert\Delta_{k}u\Delta_{k-i}v\Vert _{\mathcal{%
M}_{q}^{p}}\lesssim\sum_{j-k\leq3}\sum_{i=-1}^{1}2^{j}\Vert
\Delta_{k-i}v\Vert_{L^{\infty}}\Vert\Delta_{k}u\Vert_{\mathcal{M}_{q}^{p}} \\
& \lesssim\sum_{j-k\leq3}\sum_{i=-1}^{1}2^{j}2^{(k-i)n/p}\Vert\Delta
_{k-i}v\Vert_{\mathcal{M}_{q}^{p}}\Vert\Delta_{k}u\Vert_{\mathcal{M}_{q}^{p}}
\\
& \lesssim\sum_{j-k\leq3}\sum_{i=-1}^{1}2^{-sk}2^{j}2^{-(k-i)}\left(
2^{(k-i)(n/p+1)}\Vert\Delta_{k-i}v\Vert_{\mathcal{M}_{q}^{p}}\right) \left(
2^{sk}\Vert\Delta_{k}u\Vert_{\mathcal{M}_{q}^{p}}\right) \\
& \lesssim2^{-sj}\sum_{j-k\leq3}s^{(s+1)(j-k)}\left(
\sum_{i=-1}^{1}2^{(k-i)(n/p+1)}\Vert\Delta_{k-i}v\Vert_{\mathcal{M}%
_{q}^{p}}\right) \left( 2^{sk}\Vert\Delta_{k}u\Vert_{\mathcal{M}%
_{q}^{p}}\right) .
\end{align*}
Multiplying both sides by $2^{sj}$ and taking the $\ell^{r}$-norm yield
\begin{equation}
\Vert2^{sj}\Vert\mathcal{R}_{j}^{1}\Vert_{\mathcal{M}_{q}^{p}}\Vert_{%
\ell^{r}}\lesssim\Vert u\Vert_{\mathcal{N}_{p,q,r}^{s}}\Vert v\Vert_{%
\mathcal{N}_{p,q,\infty}^{n/p+1}},\hspace{0.5cm}\mbox{for $s>-1$}.
\label{est:comt_R_1}
\end{equation}

Similarly, for $\mathcal{R}_{j}^{2}$ note that
\begin{equation*}
\mathcal{R}_{j}^{2}=\Delta_{j}\left( T_{\nabla u}v\right) =\sum_{|j-k|\leq
4}\Delta_{j}\left( S_{k-1}\nabla u\Delta_{k}v\right) .
\end{equation*}
Then, by H\"{o}lder inequality (\ref{est:holder_linfty})
\begin{equation*}
\Vert\mathcal{R}_{j}^{2}\Vert_{\mathcal{M}_{q}^{p}}\leq\sum_{|j-k|\leq4}%
\Vert\Delta_{j}\left( S_{k-1}\nabla u\Delta_{k}v\right) \Vert_{\mathcal{M}%
_{q}^{p}}\lesssim\sum_{|j-k|\leq4}\Vert S_{k-1}\nabla
u\Vert_{L^{\infty}}\Vert\Delta_{k}v\Vert_{\mathcal{M}_{q}^{p}}.
\end{equation*}
So, using Bernstein inequality (\ref{bernstein_anel}) and (\ref%
{est:linfty_morrey}), we obtain%
\begin{align*}
\Vert S_{k-1}\nabla u\Vert_{L^{\infty}} & \leq\sum_{l-k\leq-2}\Vert
\Delta_{l}\nabla
u\Vert_{L^{\infty}}\lesssim\sum_{l-k\leq-2}2^{l(n/p+1)}\Vert\Delta_{l}u%
\Vert_{\mathcal{M}_{q}^{p}} \\
& \lesssim2^{k(n/p+1-s)}\sum_{l-k\leq-2}2^{(l-k)(n/p+1-s)}\left(
2^{sl}\Vert\Delta_{l}u\Vert_{\mathcal{M}_{q}^{p}}\right) \\
& \lesssim2^{k(n/p+1-s)}\Vert u\Vert_{\mathcal{N}_{p,q,r}^{s}},
\end{align*}
for $s<n/p+1$ or $s=n/p+1$ with $r=1$. Consequently,
\begin{align*}
\Vert\mathcal{R}_{j}^{2}\Vert_{\mathcal{M}_{q}^{p}} & \lesssim\Vert u\Vert_{%
\mathcal{N}_{p,q,r}^{s}}\sum_{|j-k|\leq4}2^{k(n/p+1-s)}\Vert\Delta
_{k}v\Vert_{\mathcal{M}_{q}^{p}} \\
& \lesssim2^{-sj}\Vert u\Vert_{\mathcal{N}_{p,q,r}^{s}}\sum_{|j-k|\leq
4}2^{s(j-k)}\left( 2^{k(n/p+1)}\Vert\Delta_{k}v\Vert_{\mathcal{M}%
_{q}^{p}}\right) .
\end{align*}
Again multiplying both sides by $2^{sj}$ and taking the $\ell^{r}$-norm, we
arrive at
\begin{equation}
\Vert2^{sj}\Vert\mathcal{R}_{j}^{2}\Vert_{\mathcal{M}_{q}^{p}}\Vert_{%
\ell^{r}}\lesssim\Vert u\Vert_{\mathcal{N}_{p,q,r}^{s}} \Vert v\Vert_{%
\mathcal{N}_{p,q,\infty}^{n/p+1}},  \label{est:comt_R_2}
\end{equation}
for $s<n/p+1$ or $s=n/p+1$ with $r=1$.

For the parcel $\mathcal{R}_{j}^{3}$, recalling the support of the Fourier
transform of $S_{k+2}\Delta_{j}\nabla u$, we may write
\begin{equation*}
\mathcal{R}_{j}^{3}=-R\left( v,\Delta_{j}\nabla u\right) =-\sum_{j-k\leq
2}\Delta_{k}vS_{k+2}(\Delta_{j}\nabla u),
\end{equation*}
which implies
\begin{equation*}
\Vert\mathcal{R}_{j}^{3}\Vert_{\mathcal{M}_{q}^{p}}\lesssim\sum_{j-k\leq
2}\Vert\Delta_{k}vS_{k+2}(\Delta_{j}\nabla u)\Vert_{\mathcal{M}%
_{q}^{p}}\lesssim\sum_{j-k\leq2}\Vert\Delta_{k}v\Vert_{\mathcal{M}%
_{q}^{p}}\Vert S_{k+2}(\Delta_{j}\nabla u)\Vert_{L^{\infty}},
\end{equation*}
where the last inequality follows from the H\"{o}lder inequality (\ref%
{est:holder_linfty}). So, by Bernstein inequality (\ref{bernstein_anel}),
\begin{align*}
\Vert S_{k+2}(\Delta_{j}\nabla u)\Vert_{L^{\infty}} & \leq\sum_{l-k\leq
1}\Vert\Delta_{l}(\Delta_{j}\nabla u)\Vert_{L^{\infty}}\lesssim\sum_{l-k\leq
1}2^{l(n/p+1)}\Vert\Delta_{l}u\Vert_{\mathcal{M}_{q}^{p}} \\
& \lesssim2^{k(n/p+1-s)}\sum_{l-k\leq1}2^{(l-k)(n/p+1-s)}\left(
2^{sl}\Vert\Delta_{l}u\Vert_{\mathcal{M}_{q}^{p}}\right) \\
& \lesssim2^{k(n/p+1-s)}\Vert u\Vert_{\mathcal{N}_{p,q,r}^{s}},
\end{align*}
for $s<n/p+1$, or $s=n/p+1$ with $r=1$, where in the last inequality we use
the H\"{o}lder inequality. Consequently, we obtain that
\begin{align*}
\Vert\mathcal{R}_{j}^{3}\Vert_{\mathcal{M}_{q}^{p}} & \lesssim\Vert u\Vert_{%
\mathcal{N}_{p,q,r}^{s}}\sum_{j-k\leq2}2^{k(n/p+1-s)}\Vert\Delta _{k}v\Vert_{%
\mathcal{M}_{q}^{p}} \\
& \lesssim2^{-sj}\Vert u\Vert_{\mathcal{N}_{p,q,r}^{s}}\sum_{j-k\leq
2}2^{s(j-k)}\left( 2^{k(n/p+1)}\Vert\Delta_{k}v\Vert_{\mathcal{M}%
_{q}^{p}}\right) .
\end{align*}
Therefore, multiplying both sides of the above inequality by $2^{sj}$,
applying the $\ell^{r}$-norm and using H\"{o}lder inequality in $\ell^{p}$%
-spaces, it follows that
\begin{equation*}
\Vert2^{sj}\Vert\mathcal{R}_{j}^{3}\Vert_{\mathcal{M}_{q}^{p}}\Vert_{%
\ell^{r}}\lesssim\Vert u\Vert_{\mathcal{N}_{p,q,r}^{s}} \Vert v\Vert_{%
\mathcal{N}_{p,q,\infty}^{n/p+1}},
\end{equation*}
for $0<s<n/p+1$ or $s=n/p+1$ with $r=1$. For the last term in (\ref%
{est:comt_decomp}), by simply changing of variables and convolution
properties, we can write
\begin{align*}
-[T_{v},\Delta_{j}]\nabla u & =-\left( T_{v}(\Delta_{j}\nabla u)-\Delta
_{j}(T_{v}\nabla u)\right) \\
& =-\sum_{|j-k|\leq4}\left( S_{k-1}v\Delta_{k}(\Delta_{j}\nabla
u)-\Delta_{j}(S_{k-1}v\Delta_{k}\nabla u)\right) \\
& =-\sum_{|j-k|\leq4}S_{k-1}v(x)\left( \varphi_{j}^{\vee}\ast\Delta
_{k}\nabla u\right) (x)-\varphi_{j}^{\vee}\ast(S_{k-1}v\Delta_{k}\nabla u)(x)
\\
& =-\sum_{|j-k|\leq4}2^{-j}\int_{\mathbb{R}^{n}}\varphi^{\vee}(y)%
\int_{0}^{1}(y\cdot\nabla)S_{k-1}v(x-2^{-j}y\tau)\;d\tau\;\Delta_{k}\nabla
u(x-2^{-j}y)dy.
\end{align*}
Then, applying the Morrey norm and using Young inequality (see \cite[Lemma
1.8]{Kozono}), since $\varphi\in\mathcal{S}$, we can estimate
\begin{align*}
\Vert\mathcal{R}_{j}^{4}\Vert_{\mathcal{M}_{q}^{p}} & \lesssim
\sum_{|j-k|\leq4}2^{-j}\Vert\nabla S_{k-1}v\Vert_{L^{\infty}}\Vert
y\cdot\varphi^{\vee}\Vert_{L^{1}}\Vert\Delta_{k}\nabla u\Vert_{\mathcal{M}%
_{q}^{p}} \\
& \lesssim\Vert\nabla v\Vert_{L^{\infty}}\sum_{|j-k|\leq4}2^{-j}\Vert
\Delta_{k}\nabla u\Vert_{\mathcal{M}_{q}^{p}} \\
& \lesssim\Vert\nabla v\Vert_{L^{\infty}}\sum_{|j-k|\leq4}2^{k-j}\Vert
\Delta_{k}u\Vert_{\mathcal{M}_{q}^{p}},
\end{align*}
where the last inequality follows from the Bernstein inequality (\ref%
{bernstein_anel}). Then
\begin{equation*}
\Vert\mathcal{R}_{j}^{4}\Vert_{\mathcal{M}_{q}^{p}}\lesssim2^{-sj}\Vert%
\nabla v\Vert_{L^{\infty}}\sum_{|j-k|\leq4}2^{(s-1)(j-k)}\left( 2^{sk}\Vert
\Delta_{k}u\Vert_{\mathcal{M}_{q}^{p}}\right) ,
\end{equation*}
and consequently
\begin{equation}
\Vert2^{sj}\Vert\mathcal{R}_{j}^{4}\Vert_{\mathcal{M}_{q}^{p}}\Vert_{%
\ell^{r}}\lesssim\Vert u\Vert_{\mathcal{N}_{p,q,r}^{s}}\Vert\nabla
v\Vert_{L^{\infty}},  \label{est:comt_R_4}
\end{equation}
for $s\in\mathbb{R}.$ From (\ref{est:comt_R_1})-(\ref{est:comt_R_4}) and (%
\ref{est:comt_decomp}), we conclude (\ref{est:commutator_uv}). Estimate (\ref%
{est:commutator_uv_2}) follows from (\ref{est:commutator_uv}) together with %
\eqref{inclusion-1} and \eqref{linfty_imersion}. \fin

\begin{remark}
\label{Rem-fields}The same result holds for $u$ and $v$ vector fields, just
considering partial derivatives in the decomposition given in (\ref%
{est:comt_decomp}).
\end{remark}

The following lemma provides a commutator estimate for the pressure term in
Besov-Morrey spaces.

\begin{lemma}
\label{lem:commutator_pi} Let $0<s<n/p$ with $1\leq r\leq\infty$, or $s=n/p$
with $r=1$, $1\leq q\leq p<\infty$ and $T>0$. Then, we have that
\begin{equation}
\Vert2^{sj}\Vert\lbrack\Delta_{j},a]\nabla\pi\Vert_{L_{T}^{1}(\mathcal{M}%
_{q}^{p})}\Vert_{\ell^{r}}\leq C \Vert a\Vert_{\widetilde{L}_{T}^{\infty }(%
\mathcal{N}_{p,q,\infty}^{n/p})} \Vert\nabla\pi\Vert_{\widetilde{L}_{T}^{1}(%
\mathcal{N}_{p,q,r}^{s})},  \label{est:commutator_pi}
\end{equation}
for all $a\in\widetilde{L}_{T}^{\infty}(\mathcal{N}_{p,q,\infty}^{n/p})$ and
$\nabla\pi\in\widetilde{L}_{T}^{1}(\mathcal{N}_{p,q,r}^{s})$.
\end{lemma}

\textbf{Proof.} Again, thanks to Bony decomposition (\ref%
{decomp_para_product}), we can write
\begin{equation}
\lbrack\Delta_{j},a]\nabla\pi=\underbrace{[\Delta_{j},T_{a}]\nabla\pi }_{=:%
\mathcal{A}_{j}^{1}}+\underbrace{\Delta_{j}\left( R(a,\nabla\pi)\right) }_{=:%
\mathcal{A}_{j}^{2}}\;\underbrace{-\;R(a,\Delta_{j}\nabla\pi )}_{=:\mathcal{A%
}_{j}^{3}}.  \label{est:comt_decomp_pi}
\end{equation}
Since $\Delta_{j}(\Delta_{k}aS_{k+2}\nabla\pi)\equiv0$ if $|j-k|\geq8$, we
have that
\begin{equation*}
\mathcal{A}_{j}^{2}=\Delta_{j}\left( R(a,\nabla\pi)\right) =\sum
_{|j-k|\leq7}\Delta_{j}\left( \Delta_{k}aS_{k+2}\nabla\pi\right) .
\end{equation*}
By H\"{o}lder inequality (\ref{est:holder_linfty}), it follows that
\begin{equation}
\Vert\mathcal{A}_{j}^{2}\Vert_{L_{T}^{1}(\mathcal{M}_{q}^{p})}\lesssim
\sum_{|j-k|\leq7}\Vert\Delta_{k}a\Vert_{L_{T}^{\infty}(\mathcal{M}%
_{q}^{p})}\Vert S_{k+2}\nabla\pi\Vert_{L_{T}^{1}(L^{\infty})}.
\label{est:comt_A_22}
\end{equation}
For $s<n/p$ or $s=n/p$ with $r=1$, we get
\begin{align*}
\Vert S_{k+2}\nabla\pi\Vert_{L_{T}^{1}(L^{\infty})} & \lesssim\sum
_{l-k\leq1}2^{nl/p}\Vert\Delta_{l}\nabla\pi\Vert_{L_{T}^{1}(\mathcal{M}%
_{q}^{p})} \\
& \lesssim2^{k(n/p-s)}\sum_{l-k\leq1}2^{(l-k)(n/p-s)}\left(
2^{sl}\Vert\Delta_{l}\nabla\pi\Vert_{L_{T}^{1}(\mathcal{M}_{q}^{p})}\right)
\\
& \lesssim2^{k(n/p-s)}\Vert\nabla\pi\Vert_{\widetilde{L}_{T}^{1}(\mathcal{N}%
_{p,q,r}^{s})}.
\end{align*}
Consequently,
\begin{align*}
\Vert\mathcal{A}_{j}^{2}\Vert_{L_{T}^{1}(\mathcal{M}_{q}^{p})} &
\lesssim\left( \sum_{|j-k|\leq7}2^{k(n/p-s)}\Vert\Delta_{k}a\Vert
_{L_{T}^{\infty}(\mathcal{M}_{q}^{p})}\right) \Vert\nabla\pi\Vert _{%
\widetilde{L}_{T}^{1}(\mathcal{N}_{p,q,r}^{s})} \\
& \lesssim2^{-sj}\left( \sum_{|j-k|\leq7}2^{(j-k)s}\left(
2^{nk/p}\Vert\Delta_{k}a\Vert_{L_{T}^{\infty}(\mathcal{M}_{q}^{p})}\right)
\right) \Vert\nabla\pi\Vert_{\widetilde{L}_{T}^{1}(\mathcal{N}_{p,q,r}^{s})}.
\end{align*}
Multiplying both sides of the above inequality by $2^{sj}$ and applying the $%
\ell^{r}$-norm, we obtain
\begin{equation}
\Vert2^{sj}\Vert\mathcal{A}_{j}^{2}\Vert_{L_{T}^{1}(\mathcal{M}%
_{q}^{p})}\Vert_{\ell^{r}}\lesssim \Vert a\Vert_{\widetilde{L}_{T}^{\infty}(%
\mathcal{N}_{p,q,\infty}^{n/p})} \Vert\nabla\pi\Vert_{\widetilde{L}_{T}^{1}(%
\mathcal{N}_{p,q,r}^{s})},  \label{est:comt_A_2}
\end{equation}
for $s<n/p$ or $s=n/p$ with $r=1$. For $\mathcal{A}_{j}^{3}$, we can express
\begin{equation*}
\mathcal{A}_{j}^{3}=-R(a,\Delta_{j}\nabla\pi)=-\sum_{j-k\leq2}\Delta
_{k}aS_{k+2}(\Delta_{j}(\nabla\pi)),
\end{equation*}
and then
\begin{equation}
\Vert\mathcal{A}_{j}^{3}\Vert_{L_{T}^{1}(\mathcal{M}_{q}^{p})}\lesssim
\sum_{j-k\leq2}\Vert\Delta_{k}a\Vert_{L_{T}^{\infty}(\mathcal{M}%
_{q}^{p})}\Vert S_{k+2}\nabla\pi\Vert_{L_{T}^{1}(L^{\infty})}.
\label{est:comt_A_32}
\end{equation}
Note that the general term of (\ref{est:comt_A_32}) is identical to the
general term of (\ref{est:comt_A_22}), so a similar argument gives us
\begin{equation}
\Vert2^{sj}\Vert\mathcal{A}_{j}^{3}\Vert_{L_{T}^{1}(\mathcal{M}%
_{q}^{p})}\Vert_{\ell^{r}}\lesssim \Vert a\Vert_{\widetilde{L}_{T}^{\infty}(%
\mathcal{N}_{p,q,\infty}^{n/p})} \Vert\nabla\pi\Vert_{\widetilde{L}_{T}^{1}(%
\mathcal{N}_{p,q,r}^{s})},  \label{est:comt_A_3}
\end{equation}
for $0<s<n/p$ or $s=n/p$ with $r=1$. For the term $\mathcal{A}_{j}^{1}$ in (%
\ref{est:comt_decomp_pi}), we proceed as follows
\begin{align*}
\Vert\mathcal{A}_{j}^{1}\Vert_{L_{T}^{1}(\mathcal{M}_{q}^{p})} &
\lesssim\sum_{|j-k|\leq4}2^{-j}\Vert\nabla S_{k-1}a\Vert_{L_{T}^{\infty
}(L^{\infty})}\Vert y\cdot\varphi^{\vee}\Vert_{L^{1}}\Vert\Delta_{k}\nabla
\pi\Vert_{L_{T}^{1}(\mathcal{M}_{q}^{p})} \\
& \lesssim\sum_{|j-k|\leq4}2^{-j}\Vert\nabla S_{k-1}a\Vert_{L_{T}^{\infty
}(L^{\infty})}\Vert\Delta_{k}\nabla\pi\Vert_{L_{T}^{1}(\mathcal{M}_{q}^{p})}.
\end{align*}
Furthermore, for $s\in\mathbb{R}$,
\begin{align*}
\Vert\nabla S_{k-1}a\Vert_{L_{T}^{\infty}(L^{\infty})} & \leq\sum
_{l-k\leq-2}\Vert\nabla\Delta_{l}a\Vert_{L_{T}^{\infty}(L^{\infty})}\lesssim%
\sum_{l-k\leq-2}2^{l(n/p+1)}\Vert\Delta_{l}a\Vert_{L_{T}^{\infty }(\mathcal{M%
}_{q}^{p})} \\
& \lesssim2^{k}\sum_{l-k\leq-2}2^{l-k}\left(
2^{nl/p}\Vert\Delta_{l}a\Vert_{L_{T}^{\infty}(\mathcal{M}_{q}^{p})}\right) \\
& \lesssim 2^{k} \Vert a\Vert_{\widetilde{L}_{T}^{\infty}(\mathcal{N}%
_{p,q,\infty}^{n/p})}.
\end{align*}
So, it follows that
\begin{align*}
\Vert\mathcal{A}_{j}^{1}\Vert_{L_{T}^{1}(\mathcal{M}_{q}^{p})} &
\lesssim\left( \sum_{|j-k|\leq4}2^{-j}2^{k}\Vert\Delta_{k}\nabla\pi
\Vert_{L_{T}^{1}(\mathcal{M}_{q}^{p})}\right) \Vert a\Vert_{\widetilde{L}%
_{T}^{\infty}(\mathcal{N}_{p,q,\infty}^{n/p})} \\
& \lesssim2^{-sj}\left( \sum_{|j-k|\leq4}2^{(j-k)(s-1)}\left(
2^{sk}\Vert\Delta_{k}\nabla\pi\Vert_{L_{T}^{1}(\mathcal{M}_{q}^{p})}\right)
\right) \Vert a\Vert_{\widetilde{L}_{T}^{\infty}(\mathcal{N}%
_{p,q,\infty}^{n/p})}.
\end{align*}
Therefore,
\begin{equation}
\Vert2^{sj}\Vert\mathcal{A}_{j}^{1}\Vert_{L_{T}^{1}(\mathcal{M}%
_{q}^{p})}\Vert_{\ell^{r}}\lesssim \Vert a\Vert_{\widetilde{L}_{T}^{\infty}(%
\mathcal{N}_{p,q,\infty}^{n/p})} \Vert\nabla\pi\Vert_{\widetilde{L}_{T}^{1}(%
\mathcal{N}_{p,q,r}^{s})},\text{ for }s\in\mathbb{R}.\text{ }
\label{est:comt_A_1}
\end{equation}
From (\ref{est:comt_A_2}), (\ref{est:comt_A_3}) and (\ref{est:comt_A_1}), we
conclude (\ref{est:commutator_pi}).

\fin

Next, adapting some arguments from \cite{Chemin_2}, we obtain estimates for
the localizations of the heat semigroup $\{e^{t\Delta}\}_{t\geq0}$ in our
setting.

\begin{lemma}
\label{lem:exp_h_2} Let $1\leq q\leq p<\infty$, $j\in\mathbb{Z}$ and $t>0.$
There are two constants $c,C>0$ (independent of $j$ and $t$) such that%
\begin{equation*}
\Vert\Delta_{j}(e^{t\Delta}u)\Vert_{\mathcal{M}_{q}^{p}}\leq
Ce^{-ct2^{2j}}\Vert\Delta_{j}u\Vert_{\mathcal{M}_{q}^{p}},
\end{equation*}
for all $u\in\mathcal{S}^{\prime}(\mathbb{R}^{n})$ satisfying $\Delta_{j}u\in%
\mathcal{M}_{q}^{p}$.
\end{lemma}

\textbf{Proof:} Consider the function $\phi\in C_{0}^{\infty}(\mathbb{R}%
^{n}\setminus\{0\})$ with $\phi\equiv1$ in a neighborhood of the ring $%
\mathcal{C}$. Recalling that $\mathrm{supp}(\widehat{\Delta_{j}u}%
)\subset2^{j}\mathcal{C}$, we have that
\begin{equation*}
e^{t\Delta}\Delta_{j}u=\phi(2^{-j}D)e^{t\Delta}\Delta_{j}u=\left[ \phi
(2^{-j}\xi)e^{-t|\xi|^{2}}\widehat{\Delta_{j}u}(\xi)\right]
^{\vee}=G_{\phi,j}(\cdot,t)\ast\Delta_{j}u,
\end{equation*}
where
\begin{equation*}
G_{\phi,j}(x,t)=\left[ \phi(2^{-j}\xi)e^{-t|\xi|^{2}}\right] ^{\vee}.
\end{equation*}
Since $\Delta_{j}(e^{t\Delta}u)=e^{t\Delta}\Delta_{j}u$, it follows that $%
\Delta_{j}(e^{t\Delta}u)=G_{\phi,j}(\cdot,t)\ast\Delta_{j}u$. Then, Young
inequality in Morrey spaces (see \cite[Lemma 1.8]{Kozono}) and estimate $%
\Vert G_{\phi,j}(\cdot,t)\Vert_{L^{1}}\leq Ce^{-ct2^{2j}}$ yield
\begin{align*}
\Vert\Delta_{j}(e^{t\Delta}u)\Vert_{\mathcal{M}_{q}^{p}} & \leq\Vert
G_{\phi,j}(\cdot,t)\Vert_{L^{1}}\Vert\Delta_{j}u\Vert_{\mathcal{M}_{q}^{p}}
\\
& \leq Ce^{-ct2^{2j}}\Vert\Delta_{j}u\Vert_{\mathcal{M}_{q}^{p}},
\end{align*}
as desired. \fin

In the lemma below we provide some estimates for the heat semigroup in
Chemin-Lerner norms based on Besov-Morrey spaces.

\begin{lemma}
\label{cor:exp_heat_2} Let $1\leq q\leq p<\infty,$ $s\in\mathbb{R}$ and $%
1\leq\beta<\infty$. There is a universal constant $C>0$ such that%
\begin{align}
\Vert e^{t\Delta}u_{0}\Vert_{\widetilde{L}_{T}^{\infty}(\mathcal{N}%
_{p,q,1}^{s})} & \leq C\Vert u_{0}\Vert_{\mathcal{N}_{p,q,1}^{s}},
\label{H-aux-1} \\
\Vert e^{t\Delta}u_{0}\Vert_{\widetilde{L}_{T}^{\beta}(\mathcal{N}%
_{p,q,1}^{s+2/\beta})} & \leq C\sum_{j\in\mathbb{Z}}2^{sj}\Vert\Delta
_{j}u_{0}\Vert_{\mathcal{M}_{q}^{p}}\left( \frac{1-e^{-cT2^{2j}\beta}}{c\beta%
}\right) ^{1/\beta}.  \label{H-aux-2}
\end{align}
Moreover, for $u_{0}\in\mathcal{N}_{p,q,1}^{s}$, we have that
\begin{equation}
\lim_{t\rightarrow0^{+}}\sum_{j\in\mathbb{Z}}2^{sj}\Vert\Delta_{j}u_{0}%
\Vert_{\mathcal{M}_{q}^{p}}\left( \frac{1-e^{-ct2^{2j}\beta}}{c\beta}\right)
^{1/\beta}=0.  \label{H-aux-3}
\end{equation}
\end{lemma}

\textbf{Proof:} By Lemma \ref{lem:exp_h_2}, we obtain that
\begin{align*}
\Vert e^{t\Delta}u_{0}\Vert_{\widetilde{L}_{T}^{\beta}(\mathcal{N}%
_{p,q,1}^{s+2/\beta})} & =\sum_{j\in\mathbb{Z}}2^{(s+2/\beta)j}\left(
\int_{0}^{T}\left( \Vert\Delta_{j}(e^{t\Delta}u_{0})\Vert_{\mathcal{M}%
_{q}^{p}}\right) ^{\beta}dt\right) ^{1/\beta} \\
& \leq C\sum_{j\in\mathbb{Z}}2^{(s+2/\beta)j}\left( \int_{0}^{T}\left(
e^{-ct2^{2j}}\Vert\Delta_{j}u_{0}\Vert_{\mathcal{M}_{q}^{p}}\right) ^{\beta
}dt\right) ^{1/\beta} \\
& =C\sum_{j\in\mathbb{Z}}2^{sj}\left( \left(
\int_{0}^{T}2^{2j}e^{-ct2^{2j}\beta}\;dt\right) \Vert\Delta_{j}u_{0}\Vert_{%
\mathcal{M}_{q}^{p}}^{\beta}\right) ^{1/\beta} \\
& =C\sum_{j\in\mathbb{Z}}2^{sj}\Vert\Delta_{j}u_{0}\Vert_{\mathcal{M}%
_{q}^{p}}\left( \frac{1-e^{-cT2^{2j}\beta}}{c\beta}\right) ^{1/\beta},
\end{align*}
which gives (\ref{H-aux-1}) and (\ref{H-aux-2}), with the natural
modification in the case $\beta=\infty$. Using $\displaystyle\frac{%
1-e^{-ct2^{2j}\beta}}{c\beta}\leq C$ and $\sum_{j\in\mathbb{Z}%
}2^{sj}\Vert\Delta_{j}u_{0}\Vert_{\mathcal{M}_{q}^{p}}=$ $\Vert u_{0}\Vert_{%
\mathcal{N}_{p,q,1}^{s}}<\infty$, we can commute the limit with the sum and
then conclude (\ref{H-aux-3}). \fin

Since we need an estimate for the density, the following lemma about the
action of a volume-preserving map $X$ in Morrey spaces will be useful. A
similar result in modified Besov-weak-Morrey spaces was proved in \cite%
{Fer-Perez-1}.

\begin{lemma}
\label{lem:X_general} Let $1<q\leq p<\infty$ and assume that $X:\mathbb{R}%
^{n}\rightarrow\mathbb{R}^{n}$ is a volume-preserving diffeomorphism such
that
\begin{equation}
|X^{\pm1}(x_{0})-X^{\pm1}(y_{0})|\leq\gamma|x_{0}-y_{0}|, \hspace{0.5cm} %
\mbox{for all \; $x_{0},y_{0}\in\mathbb{R}^{n}$},  \label{eq:4.1}
\end{equation}
and some fixed $\gamma\geq1$. Then, there exists a constant $%
C:=C(n,p,q,\gamma )>0$ such that
\begin{equation}
C^{-1}\Vert u\Vert_{\mathcal{M}_{q}^{p}}\leq\Vert u\circ X\Vert_{\mathcal{M}%
_{q}^{p}}\leq C\Vert u\Vert_{\mathcal{M}_{q}^{p}},\hspace{0.5cm} \text{for
all \; $u\in\mathcal{M}_{q}^{p}$}.  \label{est:X_composition}
\end{equation}
\end{lemma}

\textbf{Proof.} Let $x_{0}\in\mathbb{R}^{n}$, $R>0$ and consider the closed
ball $B(x_{0},R)\subset\mathbb{R}^{n}$. In view of (\ref{eq:4.1}), it
follows that
\begin{equation}
X^{\pm1}\left( B(x_{0},R)\right) \subset B\left( X^{\pm1}(x_{0}),\gamma
R\right) .  \label{eq:xb}
\end{equation}
On the other hand, for $x_{0}\in\mathbb{R}^{n},$ there is a unique $y_{0}\in%
\mathbb{R}^{n}$ such that $x_{0}=X(y_{0})$ and
\begin{equation*}
\int_{B(x_{0},R)}|u(y)|^{q}\;dy=\int_{X^{-1}\left( B(x_{0},R)\right)
}|u(X(y))|^{q}\;dy.
\end{equation*}
Then, we multiply both sides of the above equality by $R^{n/p-n/q}$ and use (%
\ref{eq:xb}) to get
\begin{align*}
R^{n/p-n/q}\Vert u\Vert_{L^{q}(B(x_{0},R))} & =R^{n/p-n/q}\Vert u\circ
X\Vert_{L^{q}\left( X^{-1}\left( B(x_{0},R)\right) \right) } \\
& \leq R^{n/p-n/q}\Vert u\circ X\Vert_{L^{q}\left( B\left(
X^{-1}(x_{0}),\gamma R\right) \right) } \\
& =R^{n/p-n/q}\Vert u\circ X\Vert_{L^{q}\left( B(y_{0},\gamma R)\right) } \\
& \leq\gamma^{n/q-n/p}\Vert u\circ X\Vert_{\mathcal{M}_{q}^{p}}.
\end{align*}
Taking the supremum over $x_{0}\in\mathbb{R}^{n}$ and $R>0$, it follows that
\begin{equation}
\Vert u\Vert_{\mathcal{M}_{q}^{p}}\leq\gamma^{n/q-n/p}\Vert u\circ X\Vert_{%
\mathcal{M}_{q}^{p}}.  \label{est:fx2}
\end{equation}
Replacing $X$ by $X^{-1}$, we arrive at
\begin{equation}
\Vert u\Vert_{\mathcal{M}_{q}^{p}}\leq\gamma^{n/q-n/p}\Vert u\circ
X^{-1}\Vert_{\mathcal{M}_{q}^{p}}.  \label{eq:fx}
\end{equation}
Thus, taking $u=v\circ X$ in (\ref{eq:fx}), we conclude that
\begin{equation}
\Vert v\circ X\Vert_{\mathcal{M}_{q}^{p}}\leq\gamma^{n/q-n/p}\Vert v\Vert_{%
\mathcal{M}_{q}^{p}}.  \label{eq:hx}
\end{equation}
From (\ref{est:fx2}) and (\ref{eq:hx}), we obtain (\ref{est:X_composition})
with $C(n,p,q,\gamma)=\gamma^{n/q-n/p}$. \fin

\begin{remark}
\label{Rem:X_with_u}We recall that the flow $X$ generated by a field $u$
with bounded gradient verifies the Lipschitz-type estimate
\begin{equation}
\left\vert X(y,t)-X(z,t)\right\vert \leq\exp\left( \int_{0}^{t}\Vert\nabla
u(\tau)\Vert_{L^{\infty}}\;d\tau\right) |y-z|.  \label{est:Xu}
\end{equation}
Also, it is not difficult to verify that an estimate similar to (\ref{est:Xu}%
) holds for $X^{-1}$.
\end{remark}

\begin{remark}
\label{rem:composition_infty} Estimate \eqref{est:X_composition} in Lemma %
\ref{lem:X_general} also holds in the $L^{\infty }$-setting, that is,
\begin{equation}
C^{-1}\Vert u\Vert _{L^{\infty }}\leq \Vert u\circ X\Vert _{L^{\infty }}\leq
C\Vert u\Vert _{L^{\infty }},\hspace{0.5cm}\text{for all}\;u\in L^{\infty }(%
\mathbb{R}^{n}).  \label{est:X_composition_infty}
\end{equation}
\end{remark}

\section{Linear estimates}

This section is dedicated to estimates for linear problems associated to
system (\ref{sist:navier_stokes_1}), which will play a key role in the proof
of Theorem \ref{the:navier_stokes_2}.

\begin{proposition}
\label{prop:transport} (Transport equation) Let $0<s<n/p+1$ with $1\leq
r\leq \infty $, or $s=n/p+1$ with $r=1$, and $1\leq q\leq p<\infty $.
Consider $a_{0}\in \mathcal{N}_{p,q,r}^{s}\cap L^{\infty }$ and a
divergence-free vector field $u\in L_{T}^{1}(\mathcal{N}_{p,q,1}^{n/p+1})$
with $\nabla u\in L_{T}^{1}(L^{\infty })$ for $T>0$. Assume also that $a\in
\widetilde{L}_{T}^{\infty }(\mathcal{N}_{p,q,r}^{s})$ solves the transport
equation
\begin{equation}
\left\{
\begin{array}{l}
\partial _{t}a+u\cdot \nabla a=0 \\
a(\cdot ,0)=a_{0}%
\end{array}%
,\right. \hspace{0.5cm}(x,t)\in \mathbb{R}^{n}\times \mathbb{R}^{+}.
\label{sist:transport}
\end{equation}%
Then, there holds
\begin{equation}
\Vert a\Vert _{\widetilde{L}_{T}^{\infty }(\mathcal{N}_{p,q,r}^{s})}\leq
\gamma ^{2(n/q-n/p)}\left( \Vert a_{0}\Vert _{\mathcal{N}_{p,q,r}^{s}}+C%
\int_{0}^{T}\Vert a(\tau )\Vert _{\mathcal{N}_{p,q,r}^{s}}\left( \Vert
u(\tau )\Vert _{\mathcal{N}_{p,q,\infty }^{n/p+1}}+\Vert \nabla u(\tau
)\Vert _{L^{\infty }}\right) d\tau \right) ,  \label{est:transport}
\end{equation}%
where $\displaystyle\gamma :=\exp \left( \int_{0}^{t}\Vert \nabla u(\tau
)\Vert _{L^{\infty }}\;d\tau \right) $. As a consequence, we have that
\begin{equation}
\Vert a\Vert _{\widetilde{L}_{T}^{\infty }(\mathcal{N}_{p,q,r}^{s}\cap
L^{\infty })}\lesssim \gamma ^{2(n/q-n/p)}\Vert a_{0}\Vert _{\mathcal{N}%
_{p,q,r}^{s}\cap L^{\infty }}\exp \left( C\gamma
^{2(n/q-n/p)}\int_{0}^{T}\Vert u(\tau )\Vert _{\mathcal{N}%
_{p,q,1}^{n/p+1}}\;d\tau \right) .  \label{est:transport_1}
\end{equation}
\end{proposition}

\textbf{Proof.} Applying the frequency projection $\Delta _{j}$ to (\ref%
{sist:transport}) and using commutator properties, we have that
\begin{equation}
\left\{
\begin{array}{l}
\partial _{t}\Delta _{j}a+u\cdot \nabla \Delta _{j}a=[u\cdot \nabla ,\Delta
_{j}]a \\
\Delta _{j}a(\cdot ,0)=\Delta _{j}a_{0}%
\end{array}%
,\right. \hspace{0.5cm}(x,t)\in \mathbb{R}^{n}\times \mathbb{R}^{+}.
\label{sist:transport_aj}
\end{equation}%
Now, consider $X:\mathbb{R}^{n}\times (0,\infty )\rightarrow \mathbb{R}^{n}$
the flow associated with the problem
\begin{equation}
\left\{
\begin{array}{l}
\partial _{t}X(y,t)=u(X(y,t),t) \\
X(y,0)=y%
\end{array}%
,\right. \;\;\;\;(y,t)\in \mathbb{R}^{n}\times \mathbb{R}^{+}.
\label{sist:field_X}
\end{equation}%
Since $\mathrm{div}\;u=0$, it follows that $X$ is a volume-preserving
diffeomorphism for every $t\geq 0$ and satisfies Lemma \ref{lem:X_general}
and Remark \ref{Rem:X_with_u}. Thus, it is not difficult to see that
\begin{equation}
\partial _{t}\left( \Delta _{j}a\left( X(y,t),t\right) \right) =\partial
_{t}\Delta _{j}a\left( X(y,t),t\right) +\left( u\cdot \nabla \Delta
_{j}a\right) \left( X(y,t),t\right) .  \label{eq:transport_X_aj}
\end{equation}%
Then, we obtain from (\ref{sist:transport_aj}) and (\ref{eq:transport_X_aj})
that
\begin{equation*}
\partial _{t}\Delta _{j}a\left( X(y,t),t\right) =\left[ u\cdot \nabla
,\Delta _{j}\right] a\left( X(y,t),t\right) .
\end{equation*}%
Integrating from $0$ to $t$ and using the initial condition in (\ref%
{sist:field_X}), we arrive at
\begin{equation*}
\Delta _{j}a\left( X(y,t),t\right) =\Delta _{j}a_{0}(y)+\int_{0}^{t}[u\cdot
\nabla ,\Delta _{j}]a\left( X(y,\tau ),\tau \right) \;d\tau ,
\end{equation*}%
which along with Lemma \ref{lem:X_general} (see (\ref{est:fx2})-(\ref{eq:hx}%
)) yields that
\begin{align*}
\Vert \Delta _{j}a(t)\Vert _{\mathcal{M}_{q}^{p}}& \leq \gamma
^{n/q-n/p}\Vert \Delta _{j}a\left( X(\cdot ,t),t\right) \Vert _{\mathcal{M}%
_{q}^{p}} \\
& \leq \gamma ^{n/q-n/p}\left( \Vert \Delta _{j}a_{0}\Vert _{\mathcal{M}%
_{q}^{p}}+\int_{0}^{t}\Vert \lbrack u\cdot \nabla ,\Delta _{j}]a\left(
X(\cdot ,\tau ),\tau \right) \Vert _{\mathcal{M}_{q}^{p}}\;d\tau \right) \\
& \leq \gamma ^{2(n/q-n/p)}\left( \Vert \Delta _{j}a_{0}\Vert _{\mathcal{M}%
_{q}^{p}}+\int_{0}^{t}\Vert \lbrack u\cdot \nabla ,\Delta _{j}]a(\tau )\Vert
_{\mathcal{M}_{q}^{p}}\;d\tau \right) .
\end{align*}%
Now, multiplying both sides of the above inequality by $2^{sj}$, taking the $%
\ell ^{r}$-norm and using the Minkowski inequality, we can estimate
\begin{equation}
\Vert a\Vert _{\widetilde{L}_{T}^{\infty }(\mathcal{N}_{p,q,r}^{s})}\leq
\gamma ^{2(n/q-n/p)}\left( \Vert a_{0}\Vert _{\mathcal{N}_{p,q,r}^{s}}+%
\int_{0}^{T}\Vert 2^{sj}\Vert \lbrack u\cdot \nabla ,\Delta _{j}]a(\tau
)\Vert _{\mathcal{M}_{q}^{p}}\Vert _{\ell ^{r}}\;d\tau \right) .
\label{est:transport_a1}
\end{equation}%
Estimate (\ref{est:transport_a1}) along with Lemma \ref{lem:commutator_uv}
gives (\ref{est:transport}). For (\ref{est:transport_1}), it follows from (%
\ref{sist:transport}) and (\ref{sist:field_X}) that $\partial
_{t}a(X(y,t),t)=0$. Integrating from $0$ at $t$, we get $%
a(X(y,t),t)=a_{0}(y) $. Finally, using estimate (\ref%
{est:X_composition_infty}) in Remark \ref{rem:composition_infty}, we obtain
\begin{equation}
\Vert a(t)\Vert _{L^{\infty }}\leq \gamma ^{n/q-n/p}\Vert a(X(\cdot
,t),t)\Vert _{L^{\infty }}\leq \gamma ^{n/q-n/p}\Vert a_{0}\Vert _{L^{\infty
}}.  \label{est:a_infty}
\end{equation}%
Combining (\ref{est:transport}) with (\ref{est:a_infty}), recalling that $%
\mathcal{N}_{p,q,1}^{n/p+1}\subset \mathcal{N}_{p,q,\infty }^{n/p+1}$ (see ({%
\ref{inclusion-1}})) and $\mathcal{N}_{p,q,1}^{n/p}\hookrightarrow L^{\infty
}$ (see ({\ref{linfty_imersion}})), and using Gr\"{o}nwall inequality, we
get (\ref{est:transport_1}).

\fin

The following proposition provides estimates for the Stokes problem.

\begin{proposition}
\label{prop:stokes}(Stokes problem) Let $T>0$, $1<q\leq p<\infty $, $1\leq
\beta _{2}\leq \beta _{1}\leq \infty $, and $1/\beta _{3}=1+1/\beta
_{1}-1/\beta _{2}.$ For $s\in \mathbb{R}$, consider $u_{0}\in \mathcal{N}%
_{p,q,1}^{s}$ and $f\in \widetilde{L}_{T}^{\beta _{2}}(\mathcal{N}%
_{p,q,1}^{s-2+2/\beta _{2}})$. Assume that $u$ is a divergence-free vector
field satisfying
\begin{equation}
\left\{
\begin{array}{l}
\partial _{t}u-\Delta u+\nabla \pi =f, \\
u(\cdot ,0)=u_{0}.%
\end{array}%
\right.   \label{sist:stokes_1_2}
\end{equation}%
Then, there exist two constants $c,C>0$ such that
\begin{align}
\Vert u\Vert _{\widetilde{L}_{T}^{\beta _{1}}(\mathcal{N}_{p,q,1}^{s+2/\beta
_{1}})}& \leq C\sum_{j\in \mathbb{Z}}2^{sj}\Vert \Delta _{j}u_{0}\Vert _{%
\mathcal{M}_{q}^{p}}\left( \frac{1-e^{-cT2^{2j}\beta _{1}}}{c\beta _{1}}%
\right) ^{1/\beta _{1}}  \notag \\
& \hspace{4cm}+\;C\sum_{j\in \mathbb{Z}}2^{(s-2+2/\beta _{2})j}\Vert \Delta
_{j}f\Vert _{L_{T}^{\beta _{2}}(\mathcal{M}_{q}^{p})}\left( \frac{%
1-e^{-cT2^{2j}\beta _{3}}}{c\beta _{3}}\right) ^{1/\beta _{3}}.
\label{est:upi_2}
\end{align}%
Moreover, for $1\leq \beta _{2}<\infty $, we have the estimates
\begin{equation}
\Vert \nabla \pi \Vert _{\widetilde{L}_{T}^{\beta _{2}}(\mathcal{N}%
_{p,q,1}^{s-2+2/\beta _{2}})}\leq C\Vert f\Vert _{\widetilde{L}_{T}^{\beta
_{2}}(\mathcal{N}_{p,q,1}^{s-2+2/\beta _{2}})},  \label{est:pressure}
\end{equation}%
and
\begin{equation}
\Vert u\Vert _{\widetilde{L}_{T}^{\infty }(\mathcal{N}_{p,q,1}^{s})}+\Vert
u\Vert _{L_{T}^{1}(\mathcal{N}_{p,q,1}^{s+2})}+\Vert \nabla \pi \Vert
_{L_{T}^{1}(\mathcal{N}_{p,q,1}^{s})}\leq C\left( \Vert u_{0}\Vert _{%
\mathcal{N}_{p,q,1}^{s}}+\Vert f\Vert _{L_{T}^{1}(\mathcal{N}%
_{p,q,1}^{s})}\right) .  \label{est:u_pressure}
\end{equation}
\end{proposition}

\textbf{Proof.} Applying the Leray projector $\mathbb{P}$ in the first
equation of system (\ref{sist:stokes_1_2}) and using $\mathrm{div}\;u=0$, we
arrive at
\begin{equation}
\partial_{t}u-\Delta u=\mathbb{P}f\text{ \ with }u(\cdot,0)=u_{0},
\label{sist:stokes_2}
\end{equation}
which can be rewritten in the integral form
\begin{equation}
u(t)=e^{t\Delta}u_{0}+\int_{0}^{t}e^{(t-\tau)\Delta}\mathbb{P}f(\tau)\;d\tau.
\label{equality:1}
\end{equation}
Applying the frequency operator $\Delta_{j}$ and Morrey norm in (\ref%
{equality:1}), we get
\begin{equation*}
\Vert\Delta_{j}u(t)\Vert_{\mathcal{M}_{q}^{p}}\leq\Vert\Delta_{j}(e^{t\Delta
}u_{0})\Vert_{\mathcal{M}_{q}^{p}}+\int_{0}^{t}\Vert\Delta_{j}(e^{(t-\tau
)\Delta}\mathbb{P}f)(\tau)\Vert_{\mathcal{M}_{q}^{p}}\;d\tau.
\end{equation*}
Also, we have that
\begin{equation}
\Vert u\Vert_{\widetilde{L}_{T}^{\beta_{1}}(\mathcal{N}_{p,q,1}^{s+2/\beta
_{1}})}\leq\Vert e^{t\Delta}u_{0}\Vert_{\widetilde{L}_{T}^{\beta_{1}}(%
\mathcal{N}_{p,q,1}^{s+2/\beta_{1}})}+\sum_{j\in\mathbb{Z}}2^{(s+2/\beta
_{1})j}\left[ \int_{0}^{T}\left( \int_{0}^{t}\Vert\Delta_{j}(e^{(t-\tau
)\Delta}\mathbb{P}f)(\tau)\Vert_{\mathcal{M}_{q}^{p}}\;d\tau\right)
^{\beta_{1}}dt\right] ^{1/\beta_{1}}.  \label{est:upi_0000}
\end{equation}
By Lemma \ref{cor:exp_heat_2}, it follows that
\begin{equation}
\Vert e^{t\Delta}u_{0}\Vert_{\widetilde{L}_{T}^{\beta_{1}}(\mathcal{N}%
_{p,q,1}^{s+2/\beta_{1}})}\leq C\sum_{j\in\mathbb{Z}}2^{sj}\Vert\Delta
_{j}u_{0}\Vert_{\mathcal{M}_{q}^{p}}\left( \frac{1-e^{-cT2^{2j}\beta_{1}}}{%
c\beta_{1}}\right) ^{1/\beta_{1}}.  \label{est:eu_0}
\end{equation}
Using Lemma \ref{lem:exp_h_2}, we can estimate
\begin{equation*}
\Vert\Delta_{j}(e^{(t-\tau)\Delta}\mathbb{P}f)(\tau)\Vert_{\mathcal{M}%
_{q}^{p}}\leq Ce^{-c(t-\tau)2^{2j}}\Vert\Delta_{j}(\mathbb{P}f)(\tau )\Vert_{%
\mathcal{M}_{q}^{p}},
\end{equation*}
and then
\begin{align}
\sum_{j\in\mathbb{Z}}2^{(s+2/\beta_{1})j}\left[ \int_{0}^{T}\left( \int
_{0}^{t}\Vert\Delta_{j}(e^{(t-\tau)\Delta}\mathbb{P}f)(\tau)\Vert _{\mathcal{%
M}_{q}^{p}}\;d\tau\right) ^{\beta_{1}}dt\right] ^{1/\beta_{1}} &
\label{est:pf} \\
& \hspace{-6cm}\lesssim\sum_{j\in\mathbb{Z}}2^{(s-2+2/\beta_{2})j}\left[
\int_{0}^{T}\left(
\int_{0}^{t}2^{2j/\beta_{3}}e^{-c(t-\tau)2^{2j}}\Vert\Delta_{j}(\mathbb{P}%
f)(\tau)\Vert_{\mathcal{M}_{q}^{p}}\;d\tau\right) ^{\beta_{1}}dt\right]
^{1/\beta_{1}},  \notag
\end{align}
since $1/\beta_{1}=1/\beta_{3}+1/\beta_{2}-1$. Now, setting $%
G_{j}(t)=2^{2j/\beta_{3}}e^{-ct2^{2j}}$ and $F_{j}(t)=\Vert\Delta_{j}(%
\mathbb{P}f)(t)\Vert_{\mathcal{M}_{q}^{p}}$, we arrive at
\begin{align}
\sum_{j\in\mathbb{Z}}2^{(s+2/\beta_{1})j}\left[ \int_{0}^{T}\left( \int
_{0}^{t}\Vert\Delta_{j}(e^{(t-\tau)\Delta}\mathbb{P}f)(\tau)\Vert _{\mathcal{%
M}_{q}^{p}}\;d\tau\right) ^{\beta_{1}}dt\right] ^{1/\beta_{1}} &  \notag \\
& \hspace{-3cm} \lesssim\sum_{j\in\mathbb{Z}}2^{(s-2+2/\beta_{2})j}\left[
\int_{0}^{T}\left( G_{j}\ast F_{j}(t)\right) ^{\beta_{1}}\;dt\right]
^{1/\beta_{1}}  \notag \\
& \hspace{-3cm} \lesssim\sum_{j\in\mathbb{Z}}2^{(s-2+2/\beta_{2})j}\Vert
G_{j}\Vert_{L^{\beta_{3}}([0,T])}\Vert F_{j}\Vert_{L^{\beta_{2}}([0,T])}
\notag \\
& \hspace{-3cm} \lesssim\sum_{j\in\mathbb{Z}}2^{(s-2+2/\beta_{2})j}\left(
\int_{0}^{T}2^{2j}e^{-ct2^{2j}\beta_{3}}\;dt\right) ^{1/\beta_{3}}\Vert
\Delta_{j}(\mathbb{P}f)\Vert_{L_{T}^{\beta_{2}}(\mathcal{M}_{q}^{p})}  \notag
\\
& \hspace{-3cm} \lesssim\sum_{j\in\mathbb{Z}}2^{(s-2+2/\beta_{2})j}\Vert%
\Delta_{j}f\Vert_{L_{T}^{\beta_{2}}(\mathcal{M}_{q}^{p})}\left( \frac{%
1-e^{-cT2^{2j}\beta_{3}}}{c\beta_{3}}\right) ^{1/\beta_{3}},  \label{est:Pfe}
\end{align}
where above we have used that $1/\beta_{1}=1/\beta_{3}+1/\beta_{2}-1$,
Young-type inequality in the time-variable, and the continuity of $\mathbb{P}
$ in $\mathcal{M}_{q}^{p}$. Inserting (\ref{est:eu_0}) and (\ref{est:Pfe})
in (\ref{est:upi_0000}) yields (\ref{est:upi_2}).

For the second part, applying the divergence in (\ref{sist:stokes_1_2}), we
have $\mathrm{div}(\nabla\pi)=\mathrm{div}f$. Then, using Bernstein
inequality (\ref{bernstein_bm_homog}) and the boundedness of Riesz
transforms in the Besov-Morrey setting (see \cite{Kozono,Mazzucato}), we
obtain (\ref{est:pressure}) for $1\leq\beta_{2}<\infty.$ Finally, using (\ref%
{est:upi_2}) and (\ref{est:pressure}) with suitable values for $\beta
_{1},\beta_{2},\beta_{3}$ and that $0\leq1-e^{-cT2^{2j}\beta}\leq1$ for $%
\beta\geq0$, we conclude (\ref{est:u_pressure}).\fin

\begin{proposition}
\label{prop:navier_stokes}(Linearized inhomogeneous Navier-Stokes system)
Let $T>0$, $1\leq q\leq p<\infty $, and $n/p-1<s\leq n/p$ with $n/p\geq 1$
or $s=n/p-1$ with $n/p>1$. Let $u_{0}\in \mathcal{N}_{p,q,1}^{s},$ $a\in
\widetilde{L}_{T}^{\infty }(\mathcal{N}_{p,q,\infty }^{n/p})\cap
L_{T}^{\infty }(L^{\infty })$, and let $v$ be a divergence-free vector field
such that $\nabla v\in L_{T}^{1}(\mathcal{N}_{p,q,1}^{n/p})$. Assume also
that $u\in \widetilde{L}_{T}^{\infty }(\mathcal{N}_{p,q,1}^{s})\cap
L_{T}^{1}(\mathcal{N}_{p,q,1}^{s+2})$ and $\nabla \pi \in L_{T}^{1}(\mathcal{%
N}_{p,q,1}^{s})$ solve the system%
\begin{equation}
\left\{
\begin{array}{l}
\partial _{t}u+v\cdot \nabla u-(1+a)(\Delta u-\nabla \pi )=0 \\
\mathrm{div}\;u=0 \\
u(\cdot ,0)=u_{0}%
\end{array}%
\right. ,\hspace{0.5cm}(x,t)\in \mathbb{R}^{n}\times \mathbb{R}^{+}.
\label{sist:navier_stokes}
\end{equation}%
Then, we have the estimate
\begin{align}
\Vert u\Vert _{\widetilde{L}_{T}^{\infty }(\mathcal{N}_{p,q,1}^{s})}+\Vert
u\Vert _{L_{T}^{1}(\mathcal{N}_{p,q,1}^{s+2})}+\Vert \nabla \pi \Vert
_{L_{T}^{1}(\mathcal{N}_{p,q,1}^{s})}&   \notag  \label{est:navier_stokes} \\
& \hspace{-3cm}\leq C\left[ \Vert u_{0}\Vert _{\mathcal{N}%
_{p,q,1}^{s}}+\Vert a\Vert _{\widetilde{L}_{T}^{\infty }(\mathcal{N}%
_{p,q,\infty }^{n/p}\cap L^{\infty })}\left( \Vert u\Vert _{L_{T}^{1}(%
\mathcal{N}_{p,q,1}^{s+2})}+\Vert \nabla \pi \Vert _{L_{T}^{1}(\mathcal{N}%
_{p,q,1}^{s})}\right) \right.  \\
& \hspace{3.5cm}+\left. \int_{0}^{T}\Vert u(\tau )\Vert _{\mathcal{N}%
_{p,q,1}^{s}}\Vert v(\tau )\Vert _{\mathcal{N}_{p,q,1}^{n/p+1}}\;d\tau %
\right] .  \notag
\end{align}%
Moreover,
\begin{align}
\Vert u\Vert _{\widetilde{L}_{T}^{\infty }(\mathcal{N}_{p,q,1}^{s})}+\Vert
u\Vert _{L_{T}^{1}(\mathcal{N}_{p,q,1}^{s+2})}+\Vert \nabla \pi \Vert
_{L_{T}^{1}(\mathcal{N}_{p,q,1}^{s})}&   \notag  \label{est:gronwall} \\
& \hspace{-3.5cm}\leq C\exp \left( C\int_{0}^{T}\Vert v(\tau )\Vert _{%
\mathcal{N}_{p,q,1}^{n/p+1}}\;d\tau \right)  \\
& \hspace{-1.55cm}\times \left[ \Vert u_{0}\Vert _{\mathcal{N}%
_{p,q,1}^{s}}+\Vert a\Vert _{\widetilde{L}_{T}^{\infty }(\mathcal{N}%
_{p,q,\infty }^{n/p}\cap L^{\infty })}\left( \Vert u\Vert _{L_{T}^{1}(%
\mathcal{N}_{p,q,1}^{s+2})}+\Vert \nabla \pi \Vert _{L_{T}^{1}(\mathcal{N}%
_{p,q,1}^{s})}\right) \right] .  \notag
\end{align}
\end{proposition}

\textbf{Proof:} Setting $F:=a(\Delta u-\nabla\pi)-v\cdot\nabla u$ and using (%
\ref{est:u_pressure}), we obtain that
\begin{equation}
\Vert u\Vert_{\widetilde{L}_{T}^{\infty}(\mathcal{N}_{p,q,1}^{s})}+\Vert
u\Vert_{L_{T}^{1}(\mathcal{N}_{p,q,1}^{s+2})}+\Vert\nabla\pi\Vert_{L_{T}^{1}(%
\mathcal{N}_{p,q,1}^{s})}\leq C\left( \Vert u_{0}\Vert_{\mathcal{N}%
_{p,q,1}^{s}}+\Vert F\Vert_{L_{T}^{1}(\mathcal{N}_{p,q,1}^{s})}\right) .
\label{est:u_pressure_with_F}
\end{equation}
Note that $F$ can be estimated as
\begin{equation}
\Vert F\Vert_{\widetilde{L}_{T}^{1}(\mathcal{N}_{p,q,r}^{s})}\leq\Vert
a\Delta u\Vert_{\widetilde{L}_{T}^{1}(\mathcal{N}_{p,q,r}^{s})}+\Vert
a\nabla\pi \Vert_{\widetilde{L}_{T}^{1}\left( \mathcal{N}_{p,q,r}^{s}\right)
}+\int _{0}^{T}\Vert v\cdot\nabla u(\tau)\Vert_{\mathcal{N}%
_{p,q,r}^{s}}\;d\tau,  \label{est:FLT}
\end{equation}
where in the last term we have used Minkowski inequality. In what follows,
we treat separately the parcels in the R.H.S. of (\ref{est:FLT}).

\textbf{Statement 1.} For $0<s<n/p$ or $s=n/p$ with $r=1$, we have that
\begin{equation}
\Vert a\Delta u\Vert_{\widetilde{L}_{T}^{1}(\mathcal{N}_{p,q,r}^{s})}%
\lesssim\left( \Vert a\Vert_{L_{T}^{\infty}(L^{\infty})}+ \Vert a\Vert _{%
\widetilde{L}_{T}^{\infty}(\mathcal{N}_{p,q,\infty}^{n/p})} \right) \Vert
u\Vert_{\widetilde{L}_{T}^{1}(\mathcal{N}_{p,q,r}^{s+2})}.
\label{est:a.lapla}
\end{equation}
In fact, first we can write $\Delta_{j}\left( a\Delta u\right) =a\Delta
_{j}\Delta u+\left[ \Delta_{j},a\right] \Delta u$. Also, by Bony
decomposition,
\begin{equation*}
a\Delta_{j}\Delta u=T_{a}\Delta_{j}\Delta u+R(a,\Delta_{j}\Delta u).
\end{equation*}
As before, by H\"{o}lder inequality (\ref{est:holder_linfty}),
\begin{align*}
\Vert R(a,\Delta_{j}\Delta u)\Vert_{L_{T}^{1}(\mathcal{M}_{q}^{p})} &
\leq\sum_{j-k\leq2}\Vert\Delta_{k}aS_{k+2}(\Delta_{j}\Delta
u)\Vert_{L_{T}^{1}(\mathcal{M}_{q}^{p})} \\
& \lesssim\sum_{j-k\leq2}\Vert\Delta_{k}a\Vert_{L_{T}^{\infty}(\mathcal{M}%
_{q}^{p})}\Vert S_{k+2}(\Delta_{j}\Delta u)\Vert_{L_{T}^{1}(L^{\infty})}.
\end{align*}
However, the estimate
\begin{align*}
\Vert S_{k+2}(\Delta_{j}\Delta u)\Vert_{L_{T}^{1}(L^{\infty})} &
\lesssim\sum_{l-k\leq1}2^{nl/p}\Vert\Delta_{l}\Delta u\Vert_{L_{T}^{1}(%
\mathcal{M}_{q}^{p})} \\
& \lesssim2^{k(n/p-s)}\sum_{l-k\leq1}2^{(l-k)(n/p-s)}\left(
2^{sl}\Vert\Delta_{l}\Delta u\Vert_{L_{T}^{1}(\mathcal{M}_{q}^{p})}\right) \\
& \lesssim2^{k(n/p-s)}\Vert\Delta u\Vert_{\widetilde{L}_{T}^{1}(\mathcal{N}%
_{p,q,r}^{s})}
\end{align*}
holds for $s<n/p$ or $s=n/p$ with $r=1$. Consequently,
\begin{align*}
\Vert R(a,\Delta_{j}\Delta u)\Vert_{L_{T}^{1}(\mathcal{M}_{q}^{p})} &
\lesssim\left( \sum_{j-k\leq2}2^{k(n/p-s)}\Vert\Delta_{k}a\Vert
_{L_{T}^{\infty}(\mathcal{M}_{q}^{p})}\right) \Vert\Delta u\Vert _{%
\widetilde{L}_{T}^{1}(\mathcal{N}_{p,q,r}^{s})} \\
& \lesssim2^{-sj}\left( \sum_{j-k\leq2}2^{s(j-k)}\left( 2^{nk/p}\Vert
\Delta_{k}a\Vert_{L_{T}^{\infty}(\mathcal{M}_{q}^{p})}\right) \right)
\Vert\Delta u\Vert_{\widetilde{L}_{T}^{1}(\mathcal{N}_{p,q,r}^{s})}.
\end{align*}
Therefore, multiplying both sides of the above inequality by $2^{sj}$ and
applying the $\ell^{r}$-norm, it follows that
\begin{equation}
\Vert2^{sj}\Vert R(a,\Delta_{j}\Delta u)\Vert_{L_{T}^{1}(\mathcal{M}%
_{q}^{p})}\Vert_{\ell^{r}}\lesssim \Vert a\Vert_{\widetilde{L}_{T}^{\infty}(%
\mathcal{N}_{p,q,\infty}^{n/p})} \Vert\Delta u\Vert_{\widetilde{L}_{T}^{1}(%
\mathcal{N}_{p,q,r}^{s})},  \label{est:28}
\end{equation}
for $0<s<n/p$ or $s=n/p$ with $r=1$. Similarly, observe that
\begin{align*}
\Vert T_{a}\Delta_{j}\Delta u\Vert_{L_{T}^{1}(\mathcal{M}_{q}^{p})} &
\leq\sum_{|j-k|\leq1}\Vert S_{k-1}a\Delta_{k}\Delta_{j}\Delta u\Vert
_{L_{T}^{1}(\mathcal{M}_{q}^{p})} \\
& \lesssim\sum_{|j-k|\leq1}\Vert
S_{k-1}a\Vert_{L_{T}^{\infty}(L^{\infty})}\Vert\Delta_{k}\Delta_{j}\Delta
u\Vert_{L_{T}^{1}(\mathcal{M}_{q}^{p})} \\
& \lesssim\Vert
a\Vert_{L_{T}^{\infty}(L^{\infty})}\sum_{|j-k|\leq1}\Vert\Delta_{k}\Delta
u\Vert_{L_{T}^{1}(\mathcal{M}_{q}^{p})}.
\end{align*}
Then, we have that
\begin{equation*}
\Vert T_{a}\Delta_{j}\Delta u\Vert_{L_{T}^{1}(\mathcal{M}_{q}^{p})}%
\lesssim2^{-sj}\Vert a\Vert_{L_{T}^{\infty}(L^{\infty})}\sum_{|j-k|\leq
1}2^{s(j-k)}\left( 2^{sk}\Vert\Delta_{k}\Delta u\Vert_{L_{T}^{1}(\mathcal{M}%
_{q}^{p})}\right) ,
\end{equation*}
which implies that
\begin{equation}
\Vert2^{sj}\Vert T_{a}\Delta_{j}\Delta u\Vert_{L_{T}^{1}(\mathcal{M}%
_{q}^{p})}\Vert_{\ell^{r}}\lesssim\Vert
a\Vert_{L_{T}^{\infty}(L^{\infty})}\Vert\Delta u\Vert_{\widetilde{L}_{T}^{1}(%
\mathcal{N}_{p,q,r}^{s})},\text{ for all }s\in\mathbb{R}.  \label{est:29}
\end{equation}
Finally, analogously to the proof of Lemma \ref{lem:commutator_pi}, we have
that
\begin{equation}
\Vert2^{sj}\Vert\lbrack\Delta_{j},a]\Delta u\Vert_{L_{T}^{1}(\mathcal{M}%
_{q}^{p})}\Vert_{\ell^{r}}\lesssim \Vert a\Vert_{\widetilde{L}_{T}^{\infty }(%
\mathcal{N}_{p,q,\infty}^{n/p})} \Vert\Delta u\Vert_{\widetilde{L}_{T}^{1}(%
\mathcal{N}_{p,q,r}^{s})},  \label{est:30}
\end{equation}
for $0<s<n/p$ or $s=n/p$ with $r=1$. From (\ref{est:28}), (\ref{est:29}), (%
\ref{est:30}), and Bernstein inequality in homogeneous Besov-Morrey spaces (%
\ref{bernstein_bm_homog}), it follows (\ref{est:a.lapla}).

\textbf{Statement 2.} Similarly to the above, using Lemma \ref%
{lem:commutator_pi} and the equality $\Delta_{j}\left( a\nabla\pi\right) =
a\Delta_{j}\nabla\pi+ \left[ \Delta_{j},a\right] \nabla\pi$, for $0<s<n/p$
or $s=n/p$ with $r=1$, it follows that
\begin{equation}  \label{est:a.pres}
\Vert a\nabla\pi\Vert_{\widetilde{L}_{T}^{1}(\mathcal{N}_{p,q,r}^{s})}
\lesssim\left( \Vert a\Vert_{L_{T}^{\infty}\left( L^{\infty }\right) } +
\Vert a\Vert_{\widetilde{L}_{T}^{\infty}(\mathcal{N}_{p,q,\infty}^{n/p})}
\right) \Vert\nabla\pi\Vert_{\widetilde{L}_{T}^{1}(\mathcal{N}_{p,q,r}^{s})}.
\end{equation}

\textbf{Statement 3.} For $s<n/p$ or $s=n/p$ with $r=1$, it follows that
\begin{equation}
\Vert v\cdot\nabla u\Vert_{\mathcal{N}_{p,q,r}^{s}}\lesssim\Vert u\Vert_{%
\mathcal{N}_{p,q,r}^{s}}\left( \Vert v\Vert_{\mathcal{N}_{p,q,r}^{n/p+1}}+%
\Vert\nabla v\Vert_{L^{\infty}}\right) .  \label{est:uv}
\end{equation}
Since $\text{div}\;u=\text{div}\;v=0,$ we have that $\text{div}(v\cdot\nabla
u)=\text{div}(u\cdot\nabla v)$, and then
\begin{equation*}
\Vert\text{div}\left( \Delta_{j}(v\cdot\nabla u)\right) \Vert_{\mathcal{M}%
_{q}^{p}}=\Vert\text{div}\left( \Delta_{j}(u\cdot\nabla v)\right) \Vert_{%
\mathcal{M}_{q}^{p}}.
\end{equation*}
Bernstein inequality (\ref{bernstein_anel}) yields $\Vert\Delta_{j}(v\cdot%
\nabla u)\Vert_{\mathcal{M}_{q}^{p}}\lesssim\Vert\Delta_{j}(u\cdot\nabla
v)\Vert_{\mathcal{M}_{q}^{p}}$. Then, using (\ref{est:holder_linfty}), we
can estimate
\begin{align*}
\Vert\Delta_{j}\left( u\cdot\nabla v\right) \Vert_{\mathcal{M}_{q}^{p}} &
\leq\sum_{|j-k|\leq4}\Vert\Delta_{j}\left( S_{k-1}u\Delta_{k}\nabla v\right)
\Vert_{\mathcal{M}_{q}^{p}}+\sum_{|j-k|\leq7}\Vert\Delta_{j}\left( \Delta
_{k}uS_{k+2}\nabla v\right) \Vert_{\mathcal{M}_{q}^{p}} \\
& \lesssim\sum_{|j-k|\leq4}\Vert S_{k-1}u\Vert_{L^{\infty}}\Vert\Delta
_{k}\nabla v\Vert_{\mathcal{M}_{q}^{p}}+\sum_{|j-k|\leq7}\Vert\Delta_{k}u%
\Vert_{\mathcal{M}_{q}^{p}}\Vert S_{k+2}\nabla v\Vert_{L^{\infty}},
\end{align*}
and also $\Vert S_{k-1}u\Vert_{L^{\infty}}\lesssim2^{k(n/p-s)}\Vert u\Vert_{%
\mathcal{N}_{p,q,r}^{s}},$ for $s<n/p$ or $s=n/p$ with $r=1$, and $\Vert
S_{k+2}\nabla v\Vert_{L^{\infty}}\lesssim\Vert\nabla v\Vert_{L^{\infty }}$.
Thus, we have that
\begin{align*}
\Vert\Delta_{j}\left( u\cdot\nabla v\right) \Vert_{\mathcal{M}_{q}^{p}} &
\lesssim2^{-sj}\Vert u\Vert_{\mathcal{N}_{p,q,r}^{s}}\left( \sum_{|j-k|\leq
4}2^{(j-k)s}\left( 2^{nk/p}\Vert\Delta_{k}\nabla v\Vert_{\mathcal{M}%
_{q}^{p}}\right) \right) \\
& \hspace{6cm}+2^{-sj}\left( \sum_{|j-k|\leq7}2^{(j-k)s}\left(
2^{sk}\Vert\Delta_{k}u\Vert_{\mathcal{M}_{q}^{p}}\right) \right) \Vert\nabla
v\Vert_{L^{\infty}}.
\end{align*}
Therefore, we obtain that
\begin{equation*}
\Vert u\cdot\nabla v\Vert_{\mathcal{N}_{p,q,r}^{s}}\lesssim\Vert u\Vert_{%
\mathcal{N}_{p,q,r}^{s}}\left( \Vert\nabla v\Vert_{\mathcal{N}%
_{p,q,r}^{n/p}}+\Vert\nabla v\Vert_{L^{\infty}}\right) .
\end{equation*}
Using Bernstein inequality (\ref{bernstein_bm_homog}), we get (\ref{est:uv}%
). Next, substituting (\ref{est:a.lapla}), (\ref{est:a.pres}) and (\ref%
{est:uv}) into (\ref{est:FLT}), we arrive at
\begin{align}
\Vert F\Vert_{\widetilde{L}_{T}^{1}(\mathcal{N}_{p,q,r}^{s})} &
\lesssim\left( \Vert a\Vert_{L_{T}^{\infty}\left( L^{\infty}\right) }+ \Vert
a\Vert_{\widetilde{L}_{T}^{\infty}(\mathcal{N}_{p,q,\infty}^{n/p})} \right)
\left( \Vert u\Vert_{\widetilde{L}_{T}^{1}(\mathcal{N}_{p,q,r}^{s+2})}+\Vert%
\nabla \pi\Vert_{\widetilde{L}_{T}^{1}(\mathcal{N}_{p,q,r}^{s})}\right)
\notag \\
& \hspace{5cm}+\int_{0}^{T}\Vert u(\tau)\Vert_{\mathcal{N}%
_{p,q,r}^{s}}\left( \Vert v(\tau)\Vert_{\mathcal{N}_{p,q,r}^{n/p+1}}+\Vert%
\nabla v(\tau)\Vert_{L^{\infty}}\right) d\tau.  \label{est:FLT2}
\end{align}
Thus, if $r=1$, from the inclusion $\mathcal{N}_{p,q,1}^{n/p}\hookrightarrow
L^{\infty}$ given in \eqref{linfty_imersion} and by estimate %
\eqref{est:l_tilde}, from \eqref{est:FLT2} we have
\begin{equation}
\Vert F\Vert_{L_{T}^{1}(\mathcal{N}_{p,q,1}^{s})}\lesssim \Vert a\Vert_{%
\widetilde{L}_{T}^{\infty}(\mathcal{N}_{p,q,\infty}^{n/p}\cap L^\infty)}
\left( \Vert u\Vert_{L_{T}^{1}(\mathcal{N}_{p,q,1}^{s+2})}+\Vert\nabla\pi%
\Vert_{L_{T}^{1}(\mathcal{N}_{p,q,1}^{s})}\right) +\int_{0}^{T}\Vert u(\tau
)\Vert_{\mathcal{N}_{p,q,r}^{s}}\Vert v(\tau)\Vert_{\mathcal{N}%
_{p,q,1}^{n/p+1}}\;d\tau.  \label{est:F_finaly}
\end{equation}
Therefore, from \eqref{est:u_pressure_with_F}, \eqref{est:FLT} and %
\eqref{est:F_finaly}, we conclude \eqref{est:navier_stokes}. From %
\eqref{est:F_finaly}, we assume $n/p+1\leq s+2$, so that the regularity of $%
u $ on the right side in \eqref{est:F_finaly} does not exceed its regularity
on the left side in \eqref{est:F_finaly}. Finally, \eqref{est:gronwall} is a
direct consequence of estimate \eqref{est:navier_stokes} along with the Gr%
\"{o}nwall inequality.

\fin

\section{Existence and uniqueness of solution}

This section is dedicated to the proof of Theorem \ref{the:navier_stokes_2}.

\subsection{Construction of approximate solutions}

The proof of the our main theorem is based on an interaction argument. For $%
m\in\mathbb{N}_{0}$, we define $\{a^{m+1}\}_{m\in\mathbb{N}_{0}}$ as the
solution of the linear transport equation
\begin{equation}
\left\{
\begin{array}{l}
\partial_{t}a^{m+1}+u^{m}\cdot\nabla a^{m+1}=0, \\
a^{m+1}(\cdot,0)=a_{0}^{m+1}:=S_{m+1}a_{0},%
\end{array}
\right.  \label{sist:transporte_m}
\end{equation}
and $\{(u^{m+1},\nabla\pi^{m+1})\}_{m\in\mathbb{N}_{0}}$ as the solution of
the Navier-Stokes equation%
\begin{equation}
\left\{
\begin{array}{l}
\partial_{t}u^{m+1}+u^{m}\cdot\nabla u^{m+1}+(1+a^{m+1})(\nabla\pi
^{m+1}-\Delta u^{m+1})=0, \\
\mathrm{div}\;u^{m+1}=0, \\
u^{m+1}(\cdot,0)=u_{0}^{m+1}:=S_{m+1}u_{0}.%
\end{array}
\right.  \label{sist:navier_stokes_m}
\end{equation}

\subsection{Uniform bounds for the approximate solutions}

\subsubsection{Local case for $n/p-1<s\leq n/p$ with $n/p\geq1$}

\label{Subsection-unif-bound-1}Based on the estimates in Propositions \ref%
{prop:transport} and \ref{prop:navier_stokes}, for $T>0$, it follows that
\begin{equation}
\Vert a^{m+1}\Vert_{\widetilde{L}_{T}^{\infty}(\mathcal{N}_{p,q,r}^{s}\cap
L^\infty)}\leq\gamma^{2(n/q-n/p)} \Vert a_{0}\Vert_{\mathcal{N}%
_{p,q,r}^{s}\cap L^\infty} \exp\left( C\gamma^{2(n/q-n/p)}\int_{0}^{T}\Vert
u^{m}(\tau)\Vert_{\mathcal{N}_{p,q,1}^{n/p+1}}\;d\tau\right) ,
\label{est:a.mmais1}
\end{equation}
for $1\leq r\leq\infty$, since $\Vert S_{m+1}a_{0}\Vert_{\mathcal{N}%
_{p,q,r}^{s}\cap L^\infty}\lesssim\Vert a_{0}\Vert_{\mathcal{N}%
_{p,q,r}^{s}\cap L^\infty}$. Also, we have that
\begin{align}
\Vert u^{m+1}\Vert_{\widetilde{L}_{T}^{\infty}(\mathcal{N}_{p,q,1}^{s})} &
+\Vert u^{m+1}\Vert_{L_{T}^{1}(\mathcal{N}_{p,q,1}^{s+2})}+\Vert\nabla
\pi^{m+1}\Vert_{L_{T}^{1}\left( \mathcal{N}_{p,q,1}^{s}\right) }  \notag \\
& \hspace{-0.5cm}\leq C\exp\left( C\int_{0}^{T}\Vert u^{m}(\tau )\Vert_{%
\mathcal{N}_{p,q,1}^{n/p+1}}\;d\tau\right)  \notag \\
& \hspace{1.15cm}\times\left[ \Vert u_{0}\Vert_{\mathcal{N}_{p,q,1}^{s}}+
\Vert a^{m+1}\Vert_{\widetilde{L}_{T}^{\infty}(\mathcal{N}%
_{p,q,\infty}^{n/p} \cap L^\infty)} \left( \Vert u^{m+1}\Vert_{\widetilde{L}%
_{T}^{1}(\mathcal{N}_{p,q,1}^{s+2})}+\Vert\nabla\pi^{m+1}\Vert_{\widetilde{L}%
_{T}^{1}(\mathcal{N}_{p,q,1}^{s})}\right) \right] ,  \label{est:u.mmais1}
\end{align}
for all $n/p-1<s\leq n/p$ and $n/p\geq1$, since $\Vert S_{m+1}u_{0}\Vert_{%
\mathcal{N}_{p,q,1}^{s}}\lesssim\Vert u_{0}\Vert_{\mathcal{N}_{p,q,1}^{s}}$.

Furthermore, note that (\ref{est:a.mmais1}) holds true when the maps $%
(X^{m})^{\pm1}$ satisfies
\begin{equation}
\left\vert (X^{m})^{\pm1}(y,t)-(X^{m})^{\pm1}(z,t)\right\vert
\leq\gamma|y-z|,  \label{est:Xm}
\end{equation}
for all $m\geq0$ and $t\in\lbrack0,T]$, where
\begin{equation}
\gamma:=\exp\left( \int_{0}^{t}\Vert\nabla
u^{m}(\tau)\Vert_{L^{\infty}}\;d\tau\right) \leq\exp\left(
C\int_{0}^{t}\Vert u^{m}(\tau)\Vert _{\mathcal{N}_{p,q,1}^{n/p+1}}\;d\tau%
\right) .  \label{est:new_field}
\end{equation}
Since $n/p-1<s\leq n/p$, by H\"{o}lder inequality and interpolation estimate
(\ref{est:interpolation}) with $\theta=(s+1-n/p)/2\in(0,1)$, we arrive at
\begin{equation*}
\int_{0}^{t}\Vert u^{m}(\tau)\Vert_{\mathcal{N}_{p,q,1}^{n/p+1}}\;d\tau\leq
t^{\theta}\Vert u^{m}\Vert_{L_{t}^{1/(1-\theta)}(\mathcal{N}%
_{p,q,1}^{n/p+1})}\leq t^{1/2}\left( \Vert u^{m}\Vert_{\widetilde{L}%
_{t}^{\infty}(\mathcal{N}_{p,q,1}^{s})}+\Vert u^{m}\Vert_{L_{t}^{1}(\mathcal{%
N}_{p,q,1}^{s+2})}\right) ,
\end{equation*}
for all $m\geq0$ and $t\in\lbrack0,T]$. Considering
\begin{equation}
\gamma=\exp\left( Ct^{\theta}\left( \Vert u^{m}\Vert_{\widetilde{L}%
_{t}^{\infty}(\mathcal{N}_{p,q,1}^{s})}+\Vert u^{m}\Vert_{L_{t}^{1}(\mathcal{%
N}_{p,q,1}^{s+2})}\right) \right) ,  \label{gamma_new}
\end{equation}
estimate (\ref{est:Xm}) holds for $\gamma$ as in (\ref{gamma_new}).

Now, setting
\begin{equation*}
\Vert(u^{m+1},\nabla\pi^{m+1})\Vert_{F_{T}^{s}}:=\Vert u^{m+1}\Vert _{%
\widetilde{L}_{T}^{\infty}(\mathcal{N}_{p,q,1}^{s})}+\Vert
u^{m+1}\Vert_{L_{T}^{1}(\mathcal{N}_{p,q,1}^{s+2})}+\Vert\nabla\pi^{m+1}%
\Vert _{L_{T}^{1}(\mathcal{N}_{p,q,1}^{s})},
\end{equation*}
and using (\ref{est:u.mmais1}) and interpolation inequality, we obtain that
\begin{align}
\Vert(u^{m+1},\nabla\pi^{m+1})\Vert_{F_{T}^{s}} & \leq C{\exp\left(
CT^{\theta}\left( \Vert u^{m}\Vert_{\widetilde{L}_{T}^{\infty}(\mathcal{N}%
_{p,q,1}^{s})}+\Vert u^{m}\Vert_{L_{T}^{1}(\mathcal{N}_{p,q,1}^{s+2})}%
\right) \right) }  \notag \\
& \hspace{3.8cm}\times\left[ \Vert u_{0}\Vert_{\mathcal{N}%
_{p,q,1}^{s}}+\Vert a^{m+1}\Vert_{\widetilde{L}_{T}^{\infty}(\mathcal{N}%
_{p,q,\infty}^{n/p}\cap
L^\infty)}\Vert(u^{m+1},\nabla\pi^{m+1})\Vert_{F_{T}^{s}}\right] .
\label{est:upmmais1}
\end{align}
Substituting (\ref{est:a.mmais1}) with $s=n/p$ into (\ref{est:upmmais1})
yields
\begin{align}
\Vert(u^{m+1},\nabla\pi^{m+1})\Vert_{F_{T}^{s}} & \leq C{\exp\left(
CT^{\theta}\left( \Vert u^{m}\Vert_{\widetilde{L}_{T}^{\infty}(\mathcal{N}%
_{p,q,1}^{s})}+\Vert u^{m}\Vert_{L_{T}^{1}(\mathcal{N}_{p,q,1}^{s+2})}%
\right) \right) }  \notag \\
& \hspace{1cm}\times\left[ \Vert u_{0}\Vert_{\mathcal{N}_{p,q,1}^{s}}+%
\gamma^{2(n/q-n/p)} \Vert a_{0}\Vert_{\mathcal{N}_{p,q,\infty}^{n/p}\cap
L^\infty}\right.  \notag \\
& \left. \hspace{2cm}\cdot\exp\left( C_{\gamma}{T^{\theta}\left( \Vert
u^{m}\Vert_{\widetilde{L}_{T}^{\infty}(\mathcal{N}_{p,q,1}^{s})}+\Vert
u^{m}\Vert_{L_{T}^{1}(\mathcal{N}_{p,q,1}^{s+2})}\right) }\right)
\Vert(u^{m+1},\nabla\pi^{m+1})\Vert_{F_{T}^{s}}\right] ,  \label{est:up.Fst}
\end{align}
for all $m\geq0$, where $C_{\gamma}:=C\gamma^{2(n/q-n/p)}$. From the above
inequalities and proceeding by induction, it is not difficult to show the
uniform boundedness. {In fact, choosing $T_{1}>0$ such that
\begin{equation}
\exp\left( CT_{1}^{\theta}\left( \Vert u^{m}\Vert_{\widetilde{L}%
_{T_{1}}^{\infty}(\mathcal{N}_{p,q,1}^{s})}+\Vert u^{m}\Vert_{L_{T_{1}}^{1}(%
\mathcal{N}_{p,q,1}^{s+2})}\right) \right) \leq2,  \label{est:exp2}
\end{equation}
estimate (\ref{est:Xm}) holds with $\gamma=2$. Define $\lambda:=2(n/q-n/p)$.
Then, for $m=0$, it follows from (\ref{est:up.Fst}) that
\begin{align}
\Vert(u^{1},\nabla\pi^{1})\Vert_{F_{T_{1}}^{s}} & \leq2C\left[ \Vert
u_{0}\Vert_{\mathcal{N}_{p,q,1}^{s}}+ 2^{\lambda}\Vert a_{0}\Vert_{\mathcal{N%
}_{p,q,\infty}^{n/p}\cap L^\infty}\cdot2^{(2^{\lambda})}
\Vert(u^{1},\nabla\pi ^{1})\Vert_{F_{T_{1}}^{s}}\right]  \notag \\
& =2C\Vert u_{0}\Vert_{\mathcal{N}_{p,q,1}^{s}}+ 2^{\lambda+2^{\lambda}}
\Vert a_{0}\Vert_{\mathcal{N}_{p,q,\infty}^{n/p}\cap
L^\infty}\Vert(u^{1},\nabla\pi^{1})\Vert_{F_{T_{1}}^{s}}.  \label{est:5.10}
\end{align}
Thus, taking $2^{\lambda+2^{\lambda}} \Vert a_{0}\Vert_{\mathcal{N}%
_{p,q,\infty}^{n/p}\cap L^\infty}\leq1/2$, estimate (\ref{est:5.10}) leads
us to
\begin{equation}
\Vert(u^{1},\nabla\pi^{1})\Vert_{F_{T_{1}}^{s}}\leq4C\Vert u_{0}\Vert_{%
\mathcal{N}_{p,q,1}^{s}}.  \label{est:induction_1-2}
\end{equation}
Now consider $\widetilde{C}>0$ and }$0<T_{2}\leq T_{1}${\ such that
\begin{equation*}
4C\Vert u_{0}\Vert_{\mathcal{N}_{p,q,1}^{s}}\leq(\widetilde{C}/2)\Vert
u_{0}\Vert_{\mathcal{N}_{p,q,1}^{s}}\;\text{ and }\;\exp(2C\widetilde{C}%
T_{2}^{\theta}\Vert u_{0}\Vert_{\mathcal{N}_{p,q,1}^{s}})\leq2.
\end{equation*}
We are going to show that
\begin{equation}
\Vert(u^{m},\nabla\pi^{m})\Vert_{F_{T_{2}}^{s}}\leq\widetilde{C}\Vert
u_{0}\Vert_{\mathcal{N}_{p,q,1}^{s}},\text{ for all }m\geq0.
\label{est:induction_new}
\end{equation}
In fact, since }$4C<\widetilde{C}${, by (\ref{est:induction_1-2}) it follows
that $(u^{1},\nabla\pi^{1})$ satisfies (\ref{est:induction_new}). Next,
suppose that $(u^{m},\nabla\pi^{m})$ satisfies (\ref{est:induction_new}). By
(\ref{gamma_new}) and induction hypotheses (\ref{est:induction_new}),
estimate (\ref{est:Xm}) holds for
\begin{equation*}
\gamma=\exp(2C\widetilde{C}T_{2}^{\theta}\Vert u_{0}\Vert_{\mathcal{N}%
_{p,q,1}^{s}})\leq2.
\end{equation*}
Then, we deduce from (\ref{est:up.Fst}) that
\begin{equation*}
\Vert(u^{m+1},\nabla\pi^{m+1})\Vert_{F_{T_{2}}^{s}}\leq2C\left( \Vert
u_{0}\Vert_{\mathcal{N}_{p,q,1}^{s}}+2^{\lambda} \Vert a_{0}\Vert_{\mathcal{N%
}_{p,q,\infty}^{n/p}\cap L^\infty}\cdot2^{(2^{\lambda})}\Vert(u^{m+1},\nabla
\pi^{m+1})\Vert_{F_{T_{2}}^{s}}\right) .
\end{equation*}
Thus, considering }$2^{1+\lambda+2^{\lambda}}C\Vert a_{0}\Vert_{\mathcal{N}%
_{p,q,\infty}^{n/p}\cap L^\infty}\leq1/2,${\ we arrive at
\begin{equation*}
\Vert(u^{m+1},\nabla\pi^{m+1})\Vert_{F_{T_{2}}^{s}}\leq2C\Vert u_{0}\Vert_{%
\mathcal{N}_{p,q,1}^{s}}+\frac{1}{2}\Vert(u^{m+1},\nabla\pi^{m+1})%
\Vert_{F_{T_{2}}^{s}},
\end{equation*}
which implies
\begin{equation}
\Vert(u^{m+1},\nabla\pi^{m+1})\Vert_{F_{T_{2}}^{s}}\leq2\cdot2C\Vert
u_{0}\Vert_{\mathcal{N}_{p,q,1}^{s}}\leq\widetilde{C}\Vert u_{0}\Vert_{%
\mathcal{N}_{p,q,1}^{s}}.  \label{est:induction_new_2}
\end{equation}
} Therefore, from (\ref{est:induction_1-2}), (\ref{est:induction_new}) and (%
\ref{est:induction_new_2}), we conclude the uniform boundedness of $%
\{(u^{m},\nabla\pi^{m})\}_{m\in\mathbb{N}_{0}}$ in the corresponding space.
The uniform boundedness of $\{a^{m}\}_{m\in\mathbb{N}_{0}}$ follows
similarly by using the boundedness of $\{u^{m}\}_{m\in\mathbb{N}_{0}}$.

\subsubsection{Local case for $s=n/p-1$ with $n/p>1$}

\label{Subsection-unif-bound-2}Since {$u^{m+1}$ satisfies (\ref%
{sist:navier_stokes_m}), we can decompose $u^{m+1}=u_{H}^{m+1}+w^{m+1}$,
where $u_{H}^{m}=e^{t\Delta}u_{0}^{m}$ is the solution of the heat equation }%
$\partial_{t}u-\Delta u=0$ with $u(x,0)=u_{0}^{m}.${\ In turn, by (\ref%
{sist:navier_stokes_m}), we have that $(w^{m+1},\nabla\pi^{m+1})$ verifies
\begin{equation}
\left\{
\begin{array}{l}
\partial_{t}w^{m+1}-\Delta w^{m+1}+\nabla\pi^{m+1}=a^{m+1}(\Delta
u^{m+1}-\nabla\pi^{m+1})-u^{m}\cdot\nabla u^{m+1}, \\
\mathrm{div}\;w^{m+1}=0, \\
w^{m+1}(\cdot,0)=0.%
\end{array}
\right.  \label{sist:navier_stokes_m_2}
\end{equation}
Applying $\Delta_{j}$ to (\ref{sist:navier_stokes_m_2}), we arrive at the
system
\begin{equation}
\left\{
\begin{array}{l}
\partial_{t}\Delta_{j}w^{m+1}-\Delta\Delta_{j}w^{m+1}+\nabla\Delta_{j}%
\pi^{m+1}=\Delta_{j}F_{(a,\pi,u)}^{m+1}, \\
\mathrm{div}\;\Delta_{j}w^{m+1}=0, \\
\Delta_{j}w^{m+1}(\cdot,0)=0,%
\end{array}
\right.  \label{sist:navier_stokes_m_3}
\end{equation}
where $\Delta_{j}F_{(a,\pi,u)}^{m+1}:=\Delta_{j}(a^{m+1}(\Delta
u^{m+1}-\nabla\pi^{m+1}))-\Delta_{j}(u^{m}\cdot\nabla u^{m+1})$. Using the
linear estimates given in Proposition \ref{prop:navier_stokes} for $w_{0}=0$%
, it follows that
\begin{equation*}
\Vert w^{m+1}\Vert_{\widetilde{L}_{T}^{\infty}(\mathcal{N}%
_{p,q,1}^{s})}+\Vert w^{m+1}\Vert_{L_{T}^{1}(\mathcal{N}_{p,q,1}^{s+2})}+%
\Vert\nabla\pi^{m+1}\Vert_{L_{T}^{1}(\mathcal{N}_{p,q,1}^{s})}\lesssim\Vert
F_{(a,\pi,u)}^{m+1}\Vert_{L_{T}^{1}(\mathcal{N}_{p,q,1}^{s})}.
\end{equation*}
Moreover, since $s=n/p-1$, by estimate (\ref{est:gronwall}) in Proposition %
\ref{prop:navier_stokes}, we have that
\begin{align}
\Vert w^{m+1}\Vert_{\widetilde{L}_{T}^{\infty}(\mathcal{N}%
_{p,q,1}^{n/p-1})}+\Vert w^{m+1}\Vert_{L_{T}^{1}(\mathcal{N}%
_{p,q,1}^{n/p+1})} & +\Vert\nabla\pi^{m+1}\Vert_{L_{T}^{1}(\mathcal{N}%
_{p,q,1}^{n/p-1})}  \label{est:wm+1_12} \\
& \hspace{-5.75cm}\leq C\exp\left( C\left( \Vert
e^{t\Delta}u_{0}^{m}\Vert_{L_{T}^{1}(\mathcal{N}_{p,q,1}^{n/p+1})}+\Vert
w^{m}\Vert_{L_{T}^{1}(\mathcal{N}_{p,q,1}^{n/p+1})}\right) \right)  \notag \\
& \hspace{-2.5cm}\times\left[ \Vert e^{t\Delta}u_{0}^{m+1}\Vert_{L_{T}^{1}(%
\mathcal{N}_{p,q,1}^{n/p+1})}\left( \Vert e^{t\Delta}u_{0}^{m}\Vert_{%
\widetilde{L}_{T}^{\infty}(\mathcal{N}_{p,q,1}^{n/p-1})}+\Vert w^{m}\Vert_{%
\widetilde{L}_{T}^{\infty}(\mathcal{N}_{p,q,1}^{n/p-1})}\right) \right.
\notag \\
& \hspace{-5.25cm} + \left. \Vert a^{m+1}\Vert_{\widetilde{L}_{T}^{\infty }(%
\mathcal{N}_{p,q,\infty}^{n/p}\cap L^\infty)}\left( \Vert
e^{t\Delta}u_{0}^{m+1}\Vert _{L_{T}^{1}(\mathcal{N}_{p,q,1}^{n/p+1})}+\Vert%
\nabla\pi^{m+1}\Vert_{L_{T}^{1}(\mathcal{N}_{p,q,1}^{n/p-1})}+\Vert
w^{m+1}\Vert_{L_{T}^{1}(\mathcal{N}_{p,q,1}^{n/p+1})}\right) \right].  \notag
\end{align}
Substituting the first estimate of Lemma \ref{cor:exp_heat_2} in (\ref%
{est:wm+1_12}), we can choose }$C_{0}\geq1/4$ to estimate{\
\begin{align}
\Vert w^{m+1}\Vert_{\widetilde{L}_{T}^{\infty}(\mathcal{N}%
_{p,q,1}^{n/p-1})}+\Vert w^{m+1}\Vert_{L_{T}^{1}(\mathcal{N}%
_{p,q,1}^{n/p+1})} & +\Vert\nabla\pi^{m+1}\Vert_{L_{T}^{1}(\mathcal{N}%
_{p,q,1}^{n/p-1})}  \label{est:wm+1_13} \\
& \hspace{-5.5cm}\leq C_{0}\exp\left( C_{0}\left( \Vert
e^{t\Delta}u_{0}\Vert_{L_{T}^{1}(\mathcal{N}_{p,q,1}^{n/p+1})}+\Vert
w^{m}\Vert_{L_{T}^{1}(\mathcal{N}_{p,q,1}^{n/p+1})}\right) \right)  \notag \\
& \hspace{-2.5cm}\times\left[ \Vert e^{t\Delta}u_{0}\Vert_{L_{T}^{1}(%
\mathcal{N}_{p,q,1}^{n/p+1})}\left( \Vert u_{0}\Vert_{\mathcal{N}%
_{p,q,1}^{n/p-1}}+\Vert w^{m}\Vert_{\widetilde{L}_{T}^{\infty}(\mathcal{N}%
_{p,q,1}^{n/p-1})}\right) \right.  \notag \\
& \hspace{-4.75cm}+ \left. \Vert a^{m+1}\Vert_{\widetilde{L}_{T}^{\infty }(%
\mathcal{N}_{p,q,\infty}^{n/p}\cap L^\infty)}\left( \Vert
e^{t\Delta}u_{0}\Vert_{L_{T}^{1}(\mathcal{N}_{p,q,1}^{n/p+1})}+\Vert\nabla%
\pi^{m+1}\Vert_{L_{T}^{1}(\mathcal{N}_{p,q,1}^{n/p-1})}+\Vert
w^{m+1}\Vert_{L_{T}^{1}(\mathcal{N}_{p,q,1}^{n/p+1})}\right) \right] ,
\notag
\end{align}
where we have used that $\Vert S_{m}u_{0}\Vert_{\mathcal{M}%
_{q}^{p}}\lesssim\Vert u_{0}\Vert_{\mathcal{M}_{q}^{p}}$. }Moreover, by (\ref%
{est:a.mmais1}), it follows that
\begin{equation}
\Vert a^{m+1}\Vert_{\widetilde{L}_{T}^{\infty}(\mathcal{N}_{p,q,r}^{n/p}\cap
L^\infty)}\leq\gamma^{\lambda} \Vert a_{0}\Vert_{\mathcal{N}%
_{p,q,r}^{n/p}\cap L^\infty} \cdot \exp\left( C_{0}\gamma^{\lambda}\left(
\Vert e^{t\Delta}u_{0}\Vert _{L_{T}^{1}(\mathcal{N}_{p,q,1}^{n/p+1})}+\Vert
w^{m}\Vert_{L_{T}^{1}(\mathcal{N}_{p,q,1}^{n/p+1})}\right) \right) ,
\label{est:transport-2}
\end{equation}
where $\lambda=2(n/q-n/p)\geq0$. Furthermore, note that (\ref%
{est:transport-2}) holds true when (\ref{est:Xm}) also holds, for all $%
m\geq0 $ and $t\in \lbrack0,T]$, where we can consider
\begin{equation}
\gamma:=\exp\left( \int_{0}^{t}\Vert\nabla
u^{m}(\tau)\Vert_{L^{\infty}}\;d\tau\right) \lesssim\exp\left( C_{0}\left(
\Vert e^{t\Delta}u_{0}\Vert_{L_{T}^{1}(\mathcal{N}_{p,q,1}^{n/p+1})}+\Vert
w^{m}\Vert_{L_{t}^{1}(\mathcal{N}_{p,q,1}^{n/p+1})}\right) \right) .
\label{est:new_field_2}
\end{equation}
So, setting
\begin{equation*}
\Vert(w^{m+1},\nabla\pi^{m+1})\Vert_{F_{T}^{n/p-1}}:=\Vert w^{m+1}\Vert_{%
\widetilde{L}_{T}^{\infty}(\mathcal{N}_{p,q,1}^{n/p-1})}+\Vert
w^{m+1}\Vert_{L_{T}^{1}(\mathcal{N}_{p,q,1}^{n/p+1})}+\Vert\nabla\pi
^{m+1}\Vert_{L_{T}^{1}(\mathcal{N}_{p,q,1}^{n/p-1})},
\end{equation*}
we can rewrite (\ref{est:wm+1_13}) as
\begin{align}
\Vert(w^{m+1},\nabla\pi^{m+1})\Vert_{F_{T}^{n/p-1}} & \leq C_{0}\exp\left(
C_{0}\left( \Vert e^{t\Delta}u_{0}\Vert_{L_{T}^{1}(\mathcal{N}%
_{p,q,1}^{n/p+1})}+\Vert w^{m}\Vert_{L_{T}^{1}(\mathcal{N}%
_{p,q,1}^{n/p+1})}\right) \right)  \label{est:wm+1_14} \\
& \hspace{2cm}\times\left[ \Vert e^{t\Delta}u_{0}\Vert_{L_{T}^{1}(\mathcal{N}%
_{p,q,1}^{n/p+1})}\left( \Vert u_{0}\Vert_{\mathcal{N}_{p,q,1}^{n/p-1}}+%
\Vert w^{m}\Vert_{\widetilde{L}_{T}^{\infty}(\mathcal{N}_{p,q,1}^{n/p-1})}%
\right) \right.  \notag \\
& \hspace{0.5cm} + \left. \Vert a^{m+1}\Vert_{\widetilde{L}_{T}^{\infty }(%
\mathcal{N}_{p,q,\infty}^{n/p}\cap L^\infty)}\left( \Vert
e^{t\Delta}u_{0}\Vert_{L_{T}^{1}(\mathcal{N}_{p,q,1}^{n/p+1})}+%
\Vert(w^{m+1},\nabla\pi^{m+1})\Vert _{F_{T}^{n/p-1}}\right) \right] .  \notag
\end{align}

In the sequel we prove the uniform boundedness of $\{(w^{m},\nabla\pi
^{m})\}_{m\in\mathbb{N}_{0}}$ by induction. Recall that $w^{0}=0$ in (\ref%
{sist:navier_stokes_m_2}). For $m=0$, by (\ref{est:new_field_2}), note that (%
\ref{est:Xm}) holds for
\begin{equation*}
\gamma=\exp\left( C_{0}\Vert e^{t\Delta}u_{0}\Vert_{L_{T}^{1}(\mathcal{N}%
_{p,q,1}^{n/p+1})}\right) .
\end{equation*}
Using Lemma \ref{cor:exp_heat_2}, we can choose $T_{1}>0$ such that
\begin{equation*}
C_{0}\Vert e^{t\Delta}u_{0}\Vert_{L_{T_{1}}^{1}(\mathcal{N}%
_{p,q,1}^{n/p+1})}\leq\frac{\eta}{4}\;\text{ and }\;\exp\left( C_{0}\Vert
e^{t\Delta}u_{0}\Vert_{L_{T_{1}}^{1}(\mathcal{N}_{p,q,1}^{n/p+1})}\right)
\leq2,
\end{equation*}
for some small constant $\eta>0$. Using transport estimate (\ref%
{est:transport-2}) and $\gamma\leq2$, we have that
\begin{equation}
\Vert a^{1}\Vert_{\widetilde{L}_{T_{1}}^{\infty}(\mathcal{N}%
_{p,q,r}^{n/p}\cap L^\infty)}\leq2^{\lambda} \Vert a_{0}\Vert_{\mathcal{N}%
_{p,q,r}^{n/p}\cap L^\infty} \cdot\exp\left( C_{0}\Vert
e^{t\Delta}u_{0}\Vert_{L_{T_{1}}^{1}(\mathcal{N}_{p,q,1}^{n/p+1})}\right)
^{2^{\lambda}}\leq2^{\lambda+2^{\lambda}} \Vert a_{0}\Vert _{\mathcal{N}%
_{p,q,r}^{n/p}\cap L^\infty}.  \label{est:a^1}
\end{equation}
Then, we can use (\ref{est:a^1}) in (\ref{est:wm+1_14}) with $m=0$ to obtain
\begin{align}
\Vert(w^{1},\nabla\pi^{1})\Vert_{F_{T_{1}}^{n/p-1}} & \leq C_{0}\exp\left(
C_{0}\Vert e^{t\Delta}u_{0}\Vert_{L_{T_{1}}^{1}(\mathcal{N}%
_{p,q,1}^{n/p+1})}\right) \times\left[ \Vert
e^{t\Delta}u_{0}\Vert_{L_{T_{1}}^{1}(\mathcal{N}_{p,q,1}^{n/p+1})}\Vert
u_{0}\Vert_{\mathcal{N}_{p,q,1}^{n/p-1}}\right.  \notag  \label{est:wm+1_15}
\\
& \hspace{3cm}+ \left. \Vert a^{1}\Vert_{\widetilde{L}_{T_{1}}^{\infty }(%
\mathcal{N}_{p,q,\infty}^{n/p}\cap L^\infty)} \left( \Vert
e^{t\Delta}u_{0}\Vert_{L_{T_{1}}^{1}(\mathcal{N}_{p,q,1}^{n/p+1})}+%
\Vert(w^{1},\nabla\pi^{1})\Vert_{F_{T_{1}}^{n/p-1}}\right) \right]  \notag \\
& \leq 2C_{0}\left[ \frac{\eta}{4C_{0}}\Vert u_{0}\Vert_{\mathcal{N}%
_{p,q,1}^{n/p-1}}+2^{\lambda+2^{\lambda}} \Vert a_{0}\Vert_{\mathcal{N}%
_{p,q,\infty}^{n/p}\cap L^\infty} \left( \frac{\eta}{4C_{0}}%
+\Vert(w^{1},\nabla\pi^{1})\Vert_{F_{T_{1}}^{n/p-1}}\right) \right]  \notag
\\
& =\frac{\eta}{2}\Vert u_{0}\Vert_{\mathcal{N}_{p,q,1}^{n/p-1}}+\frac{\eta }{%
2}2^{\lambda+2^{\lambda}} \Vert a_{0}\Vert_{\mathcal{N}_{p,q,\infty}^{n/p}%
\cap L^\infty}+C_{0}2^{1+\lambda+2^{\lambda}} \Vert a_{0}\Vert_{\mathcal{N}%
_{p,q,\infty}^{n/p}\cap
L^\infty}\Vert(w^{1},\nabla\pi^{1})\Vert_{F_{T_{1}}^{n/p-1}}.
\end{align}
As $C_{0}\geq1/4$, note that $2^{\lambda+2^{\lambda}} \Vert a_{0}\Vert_{%
\mathcal{N}_{p,q,\infty}^{n/p}\cap L^\infty}\leq1$ provided that $%
C_{0}2^{1+\lambda +2^{\lambda}} \Vert a_{0}\Vert_{\mathcal{N}%
_{p,q,\infty}^{n/p}\cap L^\infty}\leq1/2.$ Then, it follows from (\ref%
{est:wm+1_15}) that
\begin{equation*}
\Vert(w^{1},\nabla\pi^{1})\Vert_{F_{T_{1}}^{n/p-1}}\leq\frac{\eta}{2}\Vert
u_{0}\Vert_{\mathcal{N}_{p,q,1}^{n/p-1}}+\frac{\eta}{2}+\frac{1}{2}\Vert
(w^{1},\nabla\pi^{1})\Vert_{F_{T_{1}}^{n/p-1}}
\end{equation*}
and thus
\begin{equation}
\Vert(w^{1},\nabla\pi^{1})\Vert_{F_{T_{1}}^{n/p-1}}\leq\eta\left( \Vert
u_{0}\Vert_{\mathcal{N}_{p,q,1}^{n/p-1}}+1\right) .  \label{est:induction_1}
\end{equation}
Now, let $T_{2}>0$ be such that $\eta\leq1$,
\begin{equation}
C_{0}\Vert e^{t\Delta}u_{0}\Vert_{L_{T_{2}}^{1}(\mathcal{N}%
_{p,q,1}^{n/p+1})}\leq\frac{\eta}{8}\;\text{ and }\;\exp\left( \frac{\eta}{8}%
+C_{0}\eta\left( \Vert u_{0}\Vert_{\mathcal{N}_{p,q,1}^{n/p-1}}+1\right)
\right) \leq2.  \label{hip:induction}
\end{equation}
Note that $\exp(\eta/8)\leq2$. Moreover, $\eta/8\leq\eta/4$, and then (\ref%
{est:induction_1}) holds for $T_{1}=T_{2}$. By induction, suppose that
\begin{equation}
\Vert(w^{m},\nabla\pi^{m})\Vert_{F_{T_{2}}^{n/p-1}}\leq\eta\left( \Vert
u_{0}\Vert_{\mathcal{N}_{p,q,1}^{n/p-1}}+1\right) .  \label{est:induction_m}
\end{equation}
Let us prove that $(w^{m+1},\nabla\pi^{m+1})$ satisfies (\ref%
{est:induction_m}). In fact, using (\ref{hip:induction}) and (\ref%
{est:induction_m}), we arrive at
\begin{equation*}
\exp\left( C_{0}\left( \Vert e^{t\Delta}u_{0}\Vert_{L_{T_{2}}^{1}(\mathcal{N}%
_{p,q,1}^{n/p+1})}+\Vert w^{m}\Vert_{L_{T_{2}}^{1}(\mathcal{N}%
_{p,q,1}^{n/p+1})}\right) \right) \leq\exp\left( C_{0}\left( \frac{\eta }{%
8C_{0}}+\eta\left( \Vert u_{0}\Vert_{\mathcal{N}_{p,q,1}^{n/p+1}}+1\right)
\right) \right) \leq2,
\end{equation*}
and (\ref{est:Xm}) holds with $\gamma\leq2$. In view of (\ref%
{est:transport-2}), we can estimate
\begin{align}
\Vert a^{m+1}\Vert_{\widetilde{L}_{T_{2}}^{\infty}(\mathcal{N}%
_{p,q,r}^{n/p}\cap L^\infty)} & \leq2^{\lambda}\Vert a_{0}\Vert_{\mathcal{N}%
_{p,q,r}^{n/p}\cap L^\infty}\exp\left( C_{0}\left( \Vert
e^{t\Delta}u_{0}\Vert_{L_{T_{2}}^{1}(\mathcal{N}_{p,q,1}^{n/p+1})}+\Vert
w^{m}\Vert_{L_{T_{2}}^{1}(\mathcal{N}_{p,q,1}^{n/p+1})}\right) \right)
^{2^{\lambda}}  \label{est:a^m+1} \\
& \leq2^{\lambda}\Vert a_{0}\Vert_{\mathcal{N}_{p,q,r}^{n/p}\cap
L^\infty}\exp\left( C_{0}\left( \frac{\eta}{8C_{0}}+\eta\left( \Vert
u_{0}\Vert_{\mathcal{N}_{p,q,1}^{n/p+1}}+1\right) \right) \right)
^{2^{\lambda}}  \notag \\
& \leq2^{\lambda+2^{\lambda}}\Vert a_{0}\Vert_{\mathcal{N}_{p,q,r}^{n/p}\cap
L^\infty}.  \notag
\end{align}
Now applying the induction hypothesis in (\ref{est:wm+1_14}) leads us to
\begin{align*}
\Vert(w^{m+1},\nabla\pi^{m+1})\Vert_{F_{T_{2}}^{n/p-1}} & \leq2C_{0}\left[
\frac{\eta}{8C_{0}}\left( \Vert u_{0}\Vert_{\mathcal{N}_{p,q,1}^{n/p-1}}+%
\eta\left( \Vert u_{0}\Vert_{\mathcal{N}_{p,q,1}^{n/p-1}}+1\right) \right)
\right. \\
& \hspace{3cm}+ \left. \Vert a^{m+1}\Vert_{\widetilde{L}_{T_{2}}^{\infty }(%
\mathcal{N}_{p,q,\infty}^{n/p}\cap L^\infty)}\left( \frac{\eta}{8C_{0}}%
+\Vert(w^{m+1},\nabla\pi^{m+1})\Vert_{F_{T_{2}}^{n/p-1}}\right) \right] .
\end{align*}
Substituting (\ref{est:a^m+1}) into the above inequality, we obtain that
\begin{align*}
\Vert(w^{m+1},\nabla\pi^{m+1})\Vert_{F_{T_{2}}^{n/p-1}} & \leq2C_{0}\left[
\frac{\eta}{8C_{0}}\left( \Vert u_{0}\Vert_{\mathcal{N}_{p,q,1}^{n/p-1}}+%
\eta\left( \Vert u_{0}\Vert_{\mathcal{N}_{p,q,1}^{n/p-1}}+1\right) \right)
\right. \\
& \hspace{3.5cm} + \left. 2^{\lambda+2^{\lambda}} \Vert a_{0}\Vert _{%
\mathcal{N}_{p,q,\infty}^{n/p}\cap L^\infty}\left( \frac{\eta}{8C_{0}}%
+\Vert(w^{m+1},\nabla\pi^{m+1})\Vert_{F_{T_{2}}^{n/p-1}}\right) \right] \\
& =\frac{\eta}{4}\Vert u_{0}\Vert_{\mathcal{N}_{p,q,1}^{n/p-1}}+\frac{\eta }{%
4}\cdot\eta\left( \Vert u_{0}\Vert_{\mathcal{N}_{p,q,1}^{n/p-1}}+1\right) +%
\frac{\eta}{4}2^{\lambda+2^{\lambda}}\Vert a_{0}\Vert _{\mathcal{N}%
_{p,q,\infty}^{n/p}\cap L^\infty} \\
& \hspace{4.9cm}+2C_{0}2^{\lambda+2^{\lambda}}\Vert a_{0}\Vert _{\mathcal{N}%
_{p,q,\infty}^{n/p}\cap
L^\infty}\Vert(w^{m+1},\nabla\pi^{m+1})\Vert_{F_{T_{2}}^{n/p-1}}.
\end{align*}
Since $2^{\lambda+2^{\lambda}} \Vert a_{0}\Vert _{\mathcal{N}%
_{p,q,\infty}^{n/p}\cap L^\infty} \leq1,$ $2C_{0}2^{\lambda+2^{\lambda}}
\Vert a_{0}\Vert _{\mathcal{N}_{p,q,\infty}^{n/p}\cap L^\infty}\leq1/2$ and $%
\eta\leq1,$ we deduce that
\begin{align*}
\Vert(w^{m+1},\nabla\pi^{m+1})\Vert_{F_{T_{2}}^{n/p-1}} & \leq2\left( \frac{%
\eta}{4}\Vert u_{0}\Vert_{\mathcal{N}_{p,q,1}^{n/p-1}}+\frac{\eta}{4}%
\cdot\eta\left( \Vert u_{0}\Vert_{\mathcal{N}_{p,q,1}^{n/p-1}}+1\right) +%
\frac{\eta}{4}\cdot1\right) \\
& \leq\frac{\eta}{2}\Vert u_{0}\Vert_{\mathcal{N}_{p,q,1}^{n/p-1}}+\frac {%
\eta}{2}\cdot1\left( \Vert u_{0}\Vert_{\mathcal{N}_{p,q,1}^{n/p-1}}+1\right)
+\frac{\eta}{2} \\
& =\eta\left( \Vert u_{0}\Vert_{\mathcal{N}_{p,q,1}^{n/p-1}}+1\right) .
\end{align*}
Thus, we have that $\{(w^{m+1},\nabla\pi^{m+1})\}_{m\in\mathbb{N}_{0}}$
satisfies (\ref{est:induction_m}). Furthermore, according to the satisfied
hypotheses for $T_{2}>0$, it follows from (\ref{est:a^1}) and (\ref%
{est:a^m+1}) the uniform boundedness of $\{a^{m}\}$ for $T_{2}>0$. Recalling
that $u^{m}=e^{t\Delta}u_{0}^{m}+w^{m}$, using the uniform boundedness of $%
\{w^{m}\}_{m\in\mathbb{N}_{0}}$, and the estimates for $e^{t\Delta}u_{0}$ in
$\widetilde{L}_{T_{2}}^{\infty}$ and $L_{T_{2}}^{1}$, we can conclude the
uniform boundedness of $\{u^{m}\}_{m\in\mathbb{N}_{0}}$ and then the one of $%
\{(a^{m},u^{m},\nabla\pi^{m})\}_{m\in\mathbb{N}_{0}}$ in their respective
spaces. 

\begin{remark}
\label{Summary-Boundedness-1}In summary, setting
\begin{equation}
F_{T_{2}}^{s,n/p}:= \widetilde{L}_{T_{2}}^{\infty}(\mathcal{N}%
_{p,q,r}^{n/p}\cap L^\infty)\times\widetilde{L}_{T_{2}}^{\infty}(\mathcal{N}%
_{p,q,1}^{s})\cap L_{T_{2}}^{1}(\mathcal{N}_{p,q,1}^{s+2})\times
L_{T_{2}}^{1}(\mathcal{N}_{p,q,1}^{s}),  \label{space-aux-1}
\end{equation}
and considering the estimates obtained in the last two subsections, we have
proved that $\{(a^{m},u^{m},\nabla\pi^{m})\}_{m\in\mathbb{N}_{0}}$ belongs
to $F_{T_{2}}^{s,n/p}$ and $\Vert(a^{m},u^{m},\nabla\pi^{m})%
\Vert_{F_{T_{2}}^{s,n/p}}$ is bounded uniformly with respect to $m\in\mathbb{%
N}_{0}$, where
\begin{equation*}
\Vert(a^{m},u^{m},\nabla\pi^{m})\Vert_{F_{T_{2}}^{s,n/p}}:=\Vert a^{m}\Vert_{%
\widetilde{L}_{T_{2}}^{\infty}(\mathcal{N}_{p,q,r}^{n/p}\cap
L^\infty)}+\Vert u^{m}\Vert_{\widetilde{L}_{T_{2}}^{\infty}(\mathcal{N}%
_{p,q,1}^{s})}+\Vert u^{m}\Vert_{L_{T_{2}}^{1}(\mathcal{N}%
_{p,q,1}^{s+2})}+\Vert\nabla\pi^{m}\Vert_{L_{T_{2}}^{1}(\mathcal{N}%
_{p,q,1}^{s})}.
\end{equation*}
\end{remark}

\subsubsection{Global case for $s=n/p-1$ with $n/p>1$}

\label{Subsection-unif-bound-3}{For the global case, i.e., $T=\infty$, it is
sufficient to consider $u^{0}=0,$ to use (\ref{est:a.mmais1}), (\ref%
{est:u.mmais1}) and (\ref{est:Xm}) with (\ref{est:new_field}) and proceed by
induction to obtain the uniform boundedness. In fact, since $u^{0}=0$ it
follows that (\ref{est:a.mmais1}) with $m=0$ holds for $\gamma=1$. Then, we
obtain from (\ref{est:a.mmais1})-(\ref{est:u.mmais1}) that
\begin{equation*}
\Vert(u^{1},\nabla\pi^{1})\Vert_{F_{T}^{s}}\leq C\left( \Vert u_{0}\Vert_{%
\mathcal{N}_{p,q,1}^{s}}+\Vert a_{0}\Vert _{\mathcal{N}_{p,q,\infty}^{n/p}%
\cap L^\infty}\Vert(u^{1},\nabla\pi^{1})\Vert_{F_{T}^{s}}\right) ,
\end{equation*}
which leads us to
\begin{equation}
\Vert(u^{1},\nabla\pi^{1})\Vert_{F_{T}^{s}}\leq2C\Vert u_{0}\Vert _{\mathcal{%
N}_{p,q,1}^{s}},  \label{est:u1}
\end{equation}
provided that $C\Vert a_{0}\Vert _{\mathcal{N}_{p,q,\infty}^{n/p}\cap
L^\infty}\leq1/2$. Now, consider $\widetilde{C}>0$ and $c^{\prime}>0$ such
that $2C\Vert u_{0}\Vert_{\mathcal{N}_{p,q,1}^{s}}\leq(\widetilde{C}/2)\Vert
u_{0}\Vert _{\mathcal{N}_{p,q,1}^{s}} $ and also $\exp(C\widetilde{C}\Vert
u_{0}\Vert_{\mathcal{N}_{p,q,1}^{s}})\leq2$ if $\Vert u_{0}\Vert_{\mathcal{N}%
_{p,q,1}^{s}}\leq c^{\prime}$. By induction, we prove that
\begin{equation}
\Vert(u^{m},\nabla\pi^{m})\Vert_{F_{T}^{s}}\leq\widetilde{C}\Vert
u_{0}\Vert_{\mathcal{N}_{p,q,1}^{s}},\hspace{0.25cm}\text{for all $m\in
\mathbb{N}_{0}$}.  \label{est:um}
\end{equation}
Note that (\ref{est:u1}) satisfies (\ref{est:um}) for $m=1$. Suppose that (%
\ref{est:um}) holds for $(u^{m},\nabla\pi^{m})$. Then, we obtain that (\ref%
{est:Xm}) holds for $\gamma=\exp(C\widetilde{C}\Vert u_{0}\Vert _{\mathcal{N}%
_{p,q,1}^{s}})\leq2$. From (\ref{est:a.mmais1})-(\ref{est:u.mmais1}), we
deduce that
\begin{equation*}
\Vert(u^{m+1},\nabla\pi^{m+1})\Vert_{F_{T}^{s}}\leq2C\left( \Vert
u_{0}\Vert_{\mathcal{N}_{p,q,1}^{s}}+2^{\lambda} \Vert a_{0}\Vert _{\mathcal{%
N}_{p,q,\infty}^{n/p}\cap L^\infty}
\cdot2^{2\lambda}\Vert(u^{m+1},\nabla\pi^{m+1})\Vert_{F_{T}^{s}}\right) .
\end{equation*}
It follows that
\begin{equation}
\Vert(u^{m+1},\nabla\pi^{m+1})\Vert_{F_{T}^{s}}\leq4C\Vert u_{0}\Vert_{%
\mathcal{N}_{p,q,1}^{s}}\leq\widetilde{C}\Vert u_{0}\Vert _{\mathcal{N}%
_{p,q,1}^{s}},  \label{est:um+1}
\end{equation}
provided that $C2^{3\lambda+1} \Vert a_{0}\Vert _{\mathcal{N}%
_{p,q,\infty}^{n/p}\cap L^\infty} \leq1/2$. Thus, using (\ref{est:u1}), (\ref%
{est:um}) and (\ref{est:um+1}), we conclude the uniform boundedness in the
global case. As in subsection \ref{Subsection-unif-bound-1}, the boundedness
of $\{a^{m}\}_{m\in \mathbb{N}_{0}}$ follows from that of $\{u^{m}\}_{m\in%
\mathbb{N}_{0}}$. Moreover, Remark \ref{Summary-Boundedness-1} holds for $%
T_{2}=\infty$ and $s=n/p-1$. }

\subsection{Convergence of the approximate solutions}

\subsubsection{Local case for $n/p-1<s\leq n/p$ with $n/p\geq1$}

\label{Subsec-conv1}We are going to prove that $\{(a^{m},u^{m},\nabla\pi
^{m})\}_{m\in\mathbb{N}_{0}}$ is a Cauchy sequence in the space
\begin{align}  \label{space-aux-2}
G_{T}^{s,n/p-\epsilon} := \widetilde{L}_{T}^{\infty}(\mathcal{N}%
_{p,q,\infty}^{n/p-\epsilon})\times\widetilde{L}_{T}^{\infty}(\mathcal{N}%
_{p,q,1}^{s-\epsilon})\cap L_{T}^{1}(\mathcal{N}_{p,q,1}^{s+2-\epsilon})%
\times L_{T}^{1}(\mathcal{N}_{p,q,1}^{s-\epsilon}),
\end{align}
for some $0<T\leq T_{2}$ and $\epsilon>0$. For that, we estimate the
difference of the interactions
\begin{equation*}
\delta a^{m+1}:=a^{m+1}-a^{m},\hspace{0.25cm}\delta u^{m+1}:=u^{m+1}-u^{m}%
\hspace{0.25cm}\mbox{and}\hspace{0.25cm}\delta\pi^{m+1}:=\pi^{m+1}-\pi^{m},
\end{equation*}
for all $m\geq0$. From (\ref{sist:transporte_m}), it follows that
\begin{equation*}
\left\{
\begin{array}{l}
\partial_{t}\delta a^{m+1}+u^{m}\cdot\nabla\delta a^{m+1}+\delta
u^{m}\cdot\nabla a^{m}=0, \\
\delta a^{m+1}(\cdot,0)=\Delta_{m+1}a_{0}.%
\end{array}
\right.
\end{equation*}
Consider a fixed $\epsilon\in(0,1)$ and $s_{\epsilon} := s-\epsilon$ such
that $0\leq n/p-1<s_{\epsilon}<s$. Proceeding as in the proof of Proposition %
\ref{prop:transport} and using Bernstein inequality, we can estimate
\begin{align}
\Vert\delta a^{m+1}\Vert_{\widetilde{L}_{T}^{\infty}(\mathcal{N}%
_{p,q,\infty}^{s_{\epsilon}})} & \lesssim \Vert\Delta_{m+1}a_{0}\Vert _{%
\mathcal{N}_{p,q,\infty}^{s_{\epsilon}}} + \Vert\delta a^{m+1}\Vert_{%
\widetilde {L}_{T}^{\infty}(\mathcal{N}_{p,q,\infty}^{s_{\epsilon}})} \Vert
u^{m}\Vert _{L_{T}^{1}(\mathcal{N}_{p,q,1}^{n/p+1})} + \Vert\delta
u^{m}\cdot\nabla a^{m}\Vert_{L_{T}^{1}(\mathcal{N}_{p,q,\infty}^{s_{%
\epsilon}})}  \notag \\
& \lesssim2^{-\epsilon(m+1)} \Vert a_{0}\Vert_{\mathcal{N}%
_{p,q,\infty}^{s}}+\Vert\delta a^{m+1}\Vert_{\widetilde{L}_{T}^{\infty}(%
\mathcal{N}_{p,q,\infty}^{s_{\epsilon}})} \Vert u^{m}\Vert_{L_{T}^{1}(%
\mathcal{N}_{p,q,1}^{n/p+1})}  \label{est:a.cauchy} \\
& \hspace{8cm}+\Vert\delta u^{m}\Vert_{L_{T}^{1}(\mathcal{N}%
_{p,q,1}^{s_{\epsilon}+1})} \Vert a^{m}\Vert_{\widetilde{L}_{T}^{\infty}(%
\mathcal{N}_{p,q,\infty}^{n/p})}  \notag \\
& {\leq C\left[ 2^{-\epsilon(m+1)} \Vert a_{0}\Vert_{\mathcal{N}%
_{p,q,\infty}^{s}} +T^{\theta} \Vert\delta a^{m+1}\Vert_{\widetilde{L}%
_{T}^{\infty}(\mathcal{N}_{p,q,\infty}^{s_{\epsilon}})} \left( \Vert
u^{m}\Vert_{\widetilde {L}_{T}^{\infty}(\mathcal{N}_{p,q,1}^{s^{\prime}})}+%
\Vert u^{m}\Vert _{L_{T}^{1}(\mathcal{N}_{p,q,1}^{s^{\prime}+2})}\right)
\right. }  \notag \\
& {\hspace{8.5cm}\left. +\Vert\delta u^{m}\Vert_{L_{T}^{1}(\mathcal{N}%
_{p,q,1}^{s_{\epsilon}+1})} \Vert a^{m}\Vert_{\widetilde{L}_{T}^{\infty }(%
\mathcal{N}_{p,q,\infty}^{n/p})} \right] ,\label{est:a.cauchy_2}}
\end{align}
where the last inequality follows from the interpolation inequality (\ref%
{est:interpolation}) with $\theta=(s^{\prime}+1-n/p)/2\in(0,1)$ and estimate
(\ref{est:l_tilde}). From (\ref{sist:navier_stokes_m}), it is not difficult
to see that the difference $(\delta u^{m+1},\delta\nabla\pi^{m+1})$
satisfies
\begin{equation*}
\left\{
\begin{array}{l}
\partial_{t}\Delta_{j}\delta u^{m+1}-\Delta_{j}\left( \Delta\delta
u^{m+1}\right) +\Delta_{j}\left( \nabla\delta\pi^{m+1}\right) =\Delta
_{j}\left( a^{m+1}(\Delta\delta u^{m+1}-\nabla\delta\pi^{m+1})\right)
-\;\Delta_{j}(u^{m}\cdot\nabla\delta u^{m+1}) \\
\hspace{7.8cm}+\;\Delta_{j}(\delta a^{m+1}(\Delta
u^{m}-\nabla\pi^{m}))-\Delta_{j}(\delta u^{m}\cdot\nabla u^{m}), \\
\mathrm{div}\;\Delta_{j}\delta u^{m+1}=0, \\
\delta u^{m+1}(\cdot,0)=\Delta_{m+1}u_{0}.%
\end{array}
\right.
\end{equation*}
Considering $n/p-1<s_{\epsilon}<s$ as above, we can proceed as in
Proposition \ref{prop:navier_stokes} to obtain
\begin{align}
\Vert\delta u^{m+1}\Vert_{\widetilde{L}_{T}^{\infty}(\mathcal{N}%
_{p,q,1}^{s_{\epsilon}})} & +\Vert\delta u^{m+1}\Vert_{L_{T}^{1}(\mathcal{N}%
_{p,q,1}^{s_{\epsilon}+2})}+\Vert\nabla\delta\pi^{m+1}\Vert _{L_{T}^{1}(%
\mathcal{N}_{p,q,1}^{s_{\epsilon}})}  \notag \\
&\hspace{-0.4cm} \lesssim \Vert\Delta_{m+1}u_{0}\Vert_{\mathcal{N}%
_{p,q,1}^{s_{\epsilon}}}+ \Vert a^{m+1}\Vert_{\widetilde{L}_{T}^{\infty}(%
\mathcal{N}_{p,q,\infty}^{n/p} \cap L^\infty)} \left( \Vert\delta
u^{m+1}\Vert_{L_{T}^{1}(\mathcal{N}_{p,q,1}^{s_{\epsilon}+2})}+\Vert\nabla%
\delta\pi^{m+1}\Vert_{L_{T}^{1}(\mathcal{N}_{p,q,1}^{s_{\epsilon}})}\right)
\notag \\
& \hspace{1.5cm}+ \Vert\delta a^{m+1}\Vert_{\widetilde{L}_{T}^{\infty }(%
\mathcal{N}_{p,q,\infty}^{n/p-\epsilon})} \left( \Vert u^{m}\Vert_{L_{T}^{1}(%
\mathcal{N}_{p,q,1}^{s+2})}+\Vert\nabla\pi^{m}\Vert_{L_{T}^{1}(\mathcal{N}%
_{p,q,1}^{s})}\right)  \notag \\
& \hspace{3cm}+\left( \Vert\delta u^{m}\Vert_{\widetilde{L}_{T}^{\infty }(%
\mathcal{N}_{p,q,1}^{s_{\epsilon}})}+\Vert\delta u^{m+1}\Vert_{\widetilde {L}%
_{T}^{\infty}(\mathcal{N}_{p,q,1}^{s_{\epsilon}})}\right) \Vert
u^{m}\Vert_{L_{T}^{1}(\mathcal{N}_{p,q,1}^{n/p+1})}  \label{est:up.cauchy} \\
&\hspace{-0.4cm} \lesssim 2^{-\epsilon(m+1)}\Vert u_{0}\Vert_{\mathcal{N}%
_{p,q,1}^{s}}+ \Vert a^{m+1}\Vert_{\widetilde{L}_{T}^{\infty}(\mathcal{N}%
_{p,q,\infty}^{n/p}\cap L^\infty)} \left( \ \Vert\delta
u^{m+1}\Vert_{L_{T}^{1}(\mathcal{N}_{p,q,1}^{s_{\epsilon}+2})}+
\Vert\nabla\delta\pi^{m+1}\Vert_{L_{T}^{1}(\mathcal{N}_{p,q,1}^{s_{%
\epsilon}})}\right) {\ }  \notag \\
& \hspace{0.4cm}+ \Vert\delta a^{m+1}\Vert_{\widetilde{L}_{T}^{\infty }(%
\mathcal{N}_{p,q,\infty}^{n/p-\epsilon})} \left( \Vert u^{m}\Vert_{L_{T}^{1}(%
\mathcal{N}_{p,q,1}^{s+2})}+\Vert\nabla\pi^{m}\Vert_{L_{T}^{1}(\mathcal{N}%
_{p,q,1}^{s})}\right)  \notag \\
& {\hspace{0.9cm}+}\left( \Vert\delta u^{m}\Vert_{\widetilde{L}_{T}^{\infty
}(\mathcal{N}_{p,q,1}^{s_{\epsilon}})}+\Vert\delta u^{m+1}\Vert_{\widetilde {%
L}_{T}^{\infty}(\mathcal{N}_{p,q,1}^{s_{\epsilon}})}\right) \cdot T^{\theta
}\left( \Vert u^{m}\Vert_{\widetilde{L}_{T}^{\infty}(\mathcal{N}%
_{p,q,1}^{s})}+\Vert u^{m}\Vert_{L_{T}^{1}(\mathcal{N}_{p,q,1}^{s+2})}%
\right) ,  \label{est:up.cauchy_2}
\end{align}
where the last inequality follows from Bernstein inequality and (\ref%
{est:interpolation}) with $\theta=(s+1-n/p)/2\in(0,1)$. Now let
\begin{align*}
\Vert(\delta a^{m+1},\delta
u^{m+1},\nabla\delta\pi^{m+1})\Vert_{G_{T}^{s,n/p-\epsilon}} & \\
& \hspace{-2.9cm}:= \Vert\delta a^{m+1}\Vert_{\widetilde{L}_{T}^{\infty }(%
\mathcal{N}_{p,q,\infty}^{n/p-\epsilon})} +\Vert\delta u^{m+1}\Vert_{%
\widetilde {L}_{T}^{\infty}(\mathcal{N}_{p,q,1}^{s_{\epsilon}})}+\Vert\delta
u^{m+1}\Vert_{L_{T}^{1}(\mathcal{N}_{p,q,1}^{s_{\epsilon}+2})}+\Vert\nabla%
\delta \pi^{m+1}\Vert_{L_{T}^{1}(\mathcal{N}_{p,q,1}^{s_{\epsilon}})}.
\end{align*}
Using (\ref{est:up.cauchy_2}), we can estimate
\begin{align}
\Vert(\delta a^{m+1},\delta
u^{m+1},\nabla\delta\pi^{m+1})\Vert_{G_{T}^{s,n/p-\epsilon}} &  \notag \\
& \hspace{-2.5cm}\lesssim2^{-\epsilon(m+1)}\Vert u_{0}\Vert_{\mathcal{N}%
_{p,q,1}^{s}}+ \Vert a^{m+1}\Vert_{\widetilde{L}_{T}^{\infty}(\mathcal{N}%
_{p,q,\infty}^{n/p}\cap L^\infty)}\Vert(\delta a^{m+1},\delta
u^{m+1},\nabla\delta\pi ^{m+1})\Vert_{G_{T}^{s,n/p-\epsilon}}  \notag \\
& \hspace{-1cm}+ \Vert\delta a^{m+1}\Vert_{\widetilde{L}_{T}^{\infty }(%
\mathcal{N}_{p,q,\infty}^{n/p-\epsilon})} \left( \Vert u^{m}\Vert_{L_{T}^{1}(%
\mathcal{N}_{p,q,1}^{s+2})}+\Vert\nabla\pi^{m}\Vert_{L_{T}^{1}(\mathcal{N}%
_{p,q,1}^{s})}+1\right)  \notag \\
& \hspace{-2cm}+T^{\theta}\left( \Vert\delta u^{m}\Vert_{\widetilde{L}%
_{T}^{\infty}(\mathcal{N}_{p,q,1}^{s_{\epsilon}})}+\Vert\delta u^{m+1}\Vert_{%
\widetilde{L}_{T}^{\infty}(\mathcal{N}_{p,q,1}^{s_{\epsilon}})}\right)
\left( \Vert u^{m}\Vert_{\widetilde{L}_{T}^{\infty}(\mathcal{N}%
_{p,q,1}^{s})}+\Vert u^{m}\Vert_{L_{T}^{1}(\mathcal{N}_{p,q,1}^{s+2})}%
\right) .  \label{est:up.cauchy_3}
\end{align}
By the uniform boundedness of the sequence, it follows that
\begin{equation*}
\Vert u^{m}\Vert_{\widetilde{L}_{T}^{\infty}(\mathcal{N}_{p,q,1}^{s})}+\Vert
u^{m}\Vert_{L_{T}^{1}(\mathcal{N}_{p,q,1}^{s+2})}+\Vert\nabla\pi^{m}%
\Vert_{L_{T}^{1}(\mathcal{N}_{p,q,1}^{s})}\leq C_{1}\text{ \ and \ } \Vert
a^{m+1}\Vert_{\widetilde{L}_{T}^{\infty}(\mathcal{N}_{p,q,\infty}^{n/p}\cap
L^\infty)}\leq C_{2}\Vert a_{0}\Vert_{\mathcal{N}_{p,q,\infty}^{n/p}\cap
L^\infty},
\end{equation*}
and estimate (\ref{est:up.cauchy_3}) gives
\begin{align}
\Vert(\delta a^{m+1},\delta
u^{m+1},\nabla\delta\pi^{m+1})\Vert_{G_{T}^{s,n/p-\epsilon}} &  \notag \\
& \hspace{-2.5cm}\leq C\left[ 2^{-\epsilon(m+1)}\Vert u_{0}\Vert _{\mathcal{N%
}_{p,q,1}^{s}}+C_{2} \Vert a_{0}\Vert_{\mathcal{N}_{p,q,\infty}^{n/p}\cap
L^\infty}\Vert(\delta a^{m+1},\delta u^{m+1},\nabla\delta\pi^{m+1})\Vert
_{G_{T}^{s,n/p-\epsilon}}\right.  \notag \\
& \hspace{-0.8cm}\left. +\;C_{1} \Vert\delta a^{m+1}\Vert_{\widetilde{L}%
_{T}^{\infty}(\mathcal{N}_{p,q,\infty}^{n/p-\epsilon})}+C_{1}T^{\theta}%
\left( \Vert\delta u^{m}\Vert_{\widetilde{L}_{T}^{\infty}(\mathcal{N}%
_{p,q,1}^{s_{\epsilon}})}+\Vert\delta u^{m+1}\Vert_{\widetilde{L}%
_{T}^{\infty }(\mathcal{N}_{p,q,1}^{s_{\epsilon}})}\right) \right] .
\label{est:up.cauchy_4}
\end{align}
On the other hand, the estimate of $\Vert\delta a^{m+1}\Vert_{\widetilde {L}%
_{T}^{\infty}(\mathcal{N}_{p,q,r}^{n/p-\epsilon})}$ in (\ref{est:a.cauchy_2}%
), with $s=n/p$, leads us to
\begin{align}
\Vert\delta a^{m+1}\Vert_{\widetilde{L}_{T}^{\infty}(\mathcal{N}%
_{p,q,\infty}^{n/p-\epsilon})} & \leq C\left[ 2^{-\epsilon(m+1)} \Vert
a_{0}\Vert_{\mathcal{N}_{p,q,\infty}^{n/p}} + T^{\theta} \Vert\delta
a^{m+1}\Vert _{\widetilde{L}_{T}^{\infty}(\mathcal{N}_{p,q,\infty}^{n/p-%
\epsilon})} \left( \Vert u^{m}\Vert_{\widetilde{L}_{T}^{\infty}(\mathcal{N}%
_{p,q,1}^{s^{\prime}})}+\Vert u^{m}\Vert_{L_{T}^{1}(\mathcal{N}%
_{p,q,1}^{s^{\prime}+2})}\right) \right.  \notag \\
& \hspace{4.35cm} \left. + T^{\theta}\left( \Vert\delta u^{m}\Vert_{%
\widetilde{L}_{T}^{\infty}(\mathcal{N}_{p,q,1}^{s_{\epsilon}})}+\Vert\delta
u^{m}\Vert_{L_{T}^{1}(\mathcal{N}_{p,q,1}^{s_{\epsilon}+2})}\right) \Vert
a^{m}\Vert_{\widetilde{L}_{T}^{\infty}(\mathcal{N}_{p,q,\infty}^{n/p})} %
\right] \hspace{-0.075cm} ,  \label{est:a_n/p}
\end{align}
where, in the last inequality, we use interpolation inequality (\ref%
{est:interpolation}). Substituting (\ref{est:a_n/p}) into the R.H.S. of (\ref%
{est:up.cauchy_4}) and using the uniform boundedness as above, we get
\begin{align*}
\Vert(\delta a^{m+1},\delta
u^{m+1},\nabla\delta\pi^{m+1})\Vert_{G_{T}^{s,n/p-\epsilon}} & \\
&\hspace{-3.5cm} \leq C2^{-\epsilon(m+1)}\left( CC_{1} \Vert a_{0}\Vert_{%
\mathcal{N}_{p,q,\infty}^{n/p}} + \Vert u_{0}\Vert_{\mathcal{N}%
_{p,q,1}^{s}}\right) \\
& \hspace{-1.75cm} + C\left( (C_{1}+C_{1}^{2})T^{\theta}+C_{2} \Vert
a_{0}\Vert_{\mathcal{N}_{p,q,\infty}^{n/p}\cap L^\infty}\right) \Vert(\delta
a^{m+1},\delta u^{m+1},\nabla\delta\pi^{m+1})\Vert_{G_{T}^{s,n/p-\epsilon}}
\notag \\
& \hspace{-0.5cm}+CC_{1}T^{\theta}\left( 1+C_{2}\Vert a_{0}\Vert_{\mathcal{N}%
_{p,q,\infty}^{n/p}\cap L^\infty} \right) \left( \Vert\delta u^{m}\Vert _{%
\widetilde{L}_{T}^{\infty}(\mathcal{N}_{p,q,1}^{s_{\epsilon}})}+\Vert\delta
u^{m}\Vert_{L_{T}^{1}(\mathcal{N}_{p,q,1}^{s_{\epsilon}+2})}\right).
\end{align*}
Therefore, there exist $0<T_{3}\leq T_{2}$ and a small constant $c>0$ such
that if $\Vert a_{0}\Vert_{\mathcal{N}_{p,q,\infty}^{n/p}\cap L^\infty}\leq
c $, then
\begin{equation*}
C\left( (C_{1}+C_{1}^{2})T_{3}^{\theta}+C_{2}\Vert a_{0}\Vert _{\mathcal{N}%
_{p,q,\infty}^{n/p}\cap L^\infty}\right) \leq1-\frac{1}{2^{\epsilon}}\hspace{%
0.5cm} \text{and} \hspace{0.5cm} CC_{1}T_{3}^{\theta}\left( 1+C_{2} \Vert
a_{0}\Vert _{\mathcal{N}_{p,q,\infty}^{n/p}\cap L^\infty} \right) \leq\frac{1%
}{4^{\epsilon}},
\end{equation*}
which implies that
\begin{equation*}
\Vert(\delta a^{m+1},\delta
u^{m+1},\nabla\delta\pi^{m+1})\Vert_{G_{T_{3}}^{s,n/p-\epsilon}}\leq
C_{3}2^{-\epsilon m} \Vert(a_{0},u_{0})\Vert _{\mathcal{N}%
_{p,q,1,\infty}^{s,n/p}} +\frac{1}{2^{\epsilon}}\left( \Vert\delta
u^{m}\Vert_{\widetilde{L}_{T_{3}}^{\infty}(\mathcal{N}_{p,q,1}^{s_{%
\epsilon}})}+\Vert\delta u^{m}\Vert_{L_{T_{3}}^{1}(\mathcal{N}%
_{p,q,1}^{s_{\epsilon}+2})}\right) \hspace{-0.075cm} ,
\end{equation*}
for all $m\geq1,$ where $\Vert(a_{0},u_{0})\Vert_{\mathcal{N}%
_{p,q,1,\infty}^{s,n/p}}:=\Vert a_{0}\Vert_{\mathcal{N}_{p,q,\infty}^{n/p}%
\cap L^\infty} + \Vert u_{0}\Vert _{\mathcal{N}_{p,q,1}^{s}}$. In summary,
we have that
\begin{equation*}
\Vert(\delta a^{m+1},\delta
u^{m+1},\nabla\delta\pi^{m+1})\Vert_{G_{T_{3}}^{s,n/p-\epsilon}}\leq
C_{3}\cdot \frac{m}{2^{\epsilon m}} \Vert(a_{0},u_{0})\Vert_{\mathcal{N}%
_{p,q,1, \infty}^{s,n/p}} + \frac{1}{2^{\epsilon m}}\left( \Vert\delta
u^{1}\Vert_{\widetilde{L}_{T}^{\infty}(\mathcal{N}_{p,q,1}^{s_{\epsilon}})}+%
\Vert\delta u^{1}\Vert_{L_{T}^{1}(\mathcal{N}_{p,q,1}^{s_{\epsilon}+2})}%
\right) \hspace{-0.075cm} ,
\end{equation*}
and then $\{(a^{m},u^{m},\nabla\pi^{m})\}_{m\in\mathbb{N}_{0}}$ is a Cauchy
sequence in $F_{T_{3}}^{s,n/p-\epsilon}$ whose limit is denoted by $%
(a,u,\nabla\pi)$.

\subsubsection{Local case for $s=n/p-1$ with $n/p>1$}

\label{Subsec-conv2} For the critical case, we shall show that, up to a
subsequence, $\{(a^{m},u^{m},\nabla\pi^{m})\}_{m\in\mathbb{N}_{0}}$ is
convergent. The proof is based on compactness arguments. For that, we show
that the first derivative in time of $\{(a^{m},u^{m})\}_{m\in\mathbb{N}_{0}}$
is uniformly bounded in suitable spaces, which allows to employ Arzel\`{a}%
-Ascoli theorem and then obtain the existence of a subsequence of $%
\{(a^{m},u^{m},\nabla\pi^{m})\}_{m\in\mathbb{N}_{0}}$ converging to a limit $%
(a,u,\nabla\pi)$.

\hfill

\textbf{Step 1:} We are going to show that $\{\partial_{t}a^{m}\}_{m\in
\mathbb{N}_{0}}$ is uniformly bounded in $L_{T_{3}}^{2}(\mathcal{N}%
_{p,q,r}^{n/p-1})$. Initially, since $\partial_{t}a^{m+1}=-u^{m}\cdot\nabla
a^{m+1}$, we have
\begin{equation*}
\Vert\partial_{t}a^{m+1}\Vert_{L_{T_{3}}^{2}(\mathcal{N}_{p,q,r}^{n/p-1})}=%
\Vert u^{m}\cdot\nabla a^{m+1}\Vert_{L_{T_{3}}^{2}(\mathcal{N}%
_{p,q,r}^{n/p-1})}\lesssim\Vert u^{m}\Vert_{L_{T_{3}}^{2}(\mathcal{N}%
_{p,q,1}^{n/p})}\Vert a^{m+1}\Vert_{\widetilde{L}_{T_{3}}^{\infty}(\mathcal{N%
}_{p,q,r}^{n/p})}.
\end{equation*}
By interpolation inequality (\ref{est:interpolation}) with $\theta=1/2$, it
follows that
\begin{equation*}
\Vert u^{m}\Vert_{L_{T_{3}}^{2}(\mathcal{N}_{p,q,1}^{n/p})}\lesssim\Vert
u^{m}\Vert_{L_{T_{3}}^{\infty}(\mathcal{N}_{p,q,1}^{n/p-1})}^{1/2}\Vert
u^{m}\Vert_{L_{T_{3}}^{1}(\mathcal{N}_{p,q,1}^{n/p+1})}^{1/2},
\end{equation*}
which implies that
\begin{equation*}
\Vert\partial_{t}a^{m+1}\Vert_{L_{T_{3}}^{2}(\mathcal{N}_{p,q,r}^{n/p-1})}
\lesssim \Vert u^{m}\Vert_{L_{T_{3}}^{\infty}(\mathcal{N}%
_{p,q,1}^{n/p-1})}^{1/2}\Vert u^{m}\Vert_{L_{T_{3}}^{1}(\mathcal{N}%
_{p,q,1}^{n/p+1})}^{1/2} \Vert a^{m+1}\Vert_{\widetilde{L}_{T_{3}}^{\infty}(%
\mathcal{N}_{p,q,r}^{n/p})}.
\end{equation*}
By the uniform boundedness of $\{a^{m}\}_{m\in\mathbb{N}_{0}}$ and $%
\{u^{m}\}_{m\in\mathbb{N}_{0}}$, we can conclude that $\{\partial_{t}a^{m}%
\}_{m\in\mathbb{N}_{0}}$ is uniformly bounded in $L_{T_{3}}^{2}(\mathcal{N}%
_{p,q,r}^{n/p-1})$.

\hfill

\textbf{Step 2:} Let $\epsilon>0$ be such that $n/p-1-\epsilon>0$. Now we
are going to prove that $\{\nabla\pi^{m}\}_{m\in\mathbb{N}_{0}}$ is
uniformly bounded in $L_{T_{3}}^{2/(2-\epsilon)}(\mathcal{N}%
_{p,q,1}^{n/p-1-\epsilon})$. Observe that
\begin{align*}
\partial_tu^{m+1} - \Delta u^{m+1} + \nabla\pi^{m+1} = a^{m+1}( \Delta
u^{m+1} - \nabla\pi^{m+1}) - u^{m}\cdot\nabla u^{m+1}.
\end{align*}
{Then, using estimate (\ref{est:pressure}) in Proposition \ref{prop:stokes}
with $\beta_{2}=2/(2-\epsilon)$ and $s=n/p-1$ and proceeding as in
Proposition \ref{prop:navier_stokes},} we can estimate
\begin{align*}
\Vert\nabla\pi^{m+1}\Vert_{L_{T_{3}}^{2/(2-\epsilon)}(\mathcal{N}%
_{p,q,1}^{n/p-1-\epsilon})} & \lesssim \Vert a^{m+1}\Vert_{\widetilde {L}%
_{T_{3}}^{\infty}(\mathcal{N}_{p,q,\infty}^{n/p}\cap
L^\infty)}\Vert\nabla\pi^{m+1}\Vert_{L_{T_{3}}^{2/(2-\epsilon)}(\mathcal{N}%
_{p,q,1}^{n/p-1-\epsilon})} \\
& \hspace{1.15cm} + \left(\Vert a^{m+1}\Vert_{\widetilde{L}_{T_{3}}^{\infty
}(\mathcal{N}_{p,q,\infty}^{n/p}\cap L^\infty)} + \Vert u^{m}\Vert_{%
\widetilde{L}_{T_{3}}^{\infty }(\mathcal{N}_{p,q,1}^{n/p-1})}\right) \Vert
u^{m+1}\Vert_{L_{T_{3}}^{2/(2-\epsilon)}(\mathcal{N}_{p,q,1}^{n/p+1-%
\epsilon})}.
\end{align*}
By interpolation inequality (\ref{est:interpolation}) with $%
\theta=\epsilon/2 $ and $n/p+1-\epsilon=\theta(n/p-1)+(1-\theta)(n/p+1)$, we
have that
\begin{align*}
\Vert\nabla\pi^{m+1}\Vert_{L_{T_{3}}^{2/(2-\epsilon)}(\mathcal{N}%
_{p,q,1}^{n/p-1-\epsilon})} & \\
&\hspace{-2.5cm} \lesssim \Vert a^{m+1}\Vert_{\widetilde {L}%
_{T_{3}}^{\infty}(\mathcal{N}_{p,q,\infty}^{n/p}\cap
L^\infty)}\Vert\nabla\pi^{m+1}\Vert_{L_{T_{3}}^{2/(2-\epsilon)}(\mathcal{N}%
_{p,q,1}^{n/p-1-\epsilon})} \\
& \hspace{-1.25cm}+\left(\Vert a^{m+1}\Vert_{\widetilde{L}_{T_{3}}^{\infty }(%
\mathcal{N}_{p,q,\infty}^{n/p}\cap L^\infty)} + \Vert u^{m}\Vert_{\widetilde{%
L}_{T_{3}}^{\infty }(\mathcal{N}_{p,q,1}^{n/p-1})}\right) \Vert
u^{m+1}\Vert_{\widetilde {L}_{T_{3}}^{\infty}(\mathcal{N}%
_{p,q,1}^{n/p-1})}^{\epsilon/2}\Vert u^{m+1}\Vert_{L_{T_{3}}^{1}(\mathcal{N}%
_{p,q,1}^{n/p+1})}^{1-\epsilon/2}.
\end{align*}
Using the smallness condition on $\Vert a_{0}\Vert_{\mathcal{N}%
_{p,q,\infty}^{n/p} \cap L^\infty}$ and the uniform estimates obtained in
previous sections, we can conclude that $\{\nabla\pi^{m}\}_{m\in\mathbb{N}%
_{0}}$ is uniformly bounded in $L_{T_{3}}^{2/(2-\epsilon)}(\mathcal{N}%
_{p,q,1}^{n/p-1-\epsilon})$.

\hfill

\textbf{Step 3:} Finally, we prove that $\{\partial_{t}u^{m}\}_{m\in \mathbb{%
N}_{0}}$ is uniformly bounded in $L_{T_{3}}^{2/(2-\epsilon )}(\mathcal{N}%
_{p,q,1}^{n/p-1-\epsilon})$, where $\epsilon>0$ is as before. Analogously,
since $\partial_{t}u^{m}=(1+a^{m+1})(\Delta u^{m}-\nabla\pi
^{m+1})-u^{m}\cdot\nabla u^{m+1},$ it follows that
\begin{align*}
\Vert\partial_{t}u^{m+1}\Vert_{L_{T_{3}}^{2/(2-\epsilon)}(\mathcal{N}%
_{p,q,1}^{n/p-1-\epsilon})} & \\
& \hspace{-3cm} \lesssim \left( 1+ \Vert a^{m+1}\Vert_{\widetilde{L}%
_{T_{3}}^{\infty}(\mathcal{N}_{p,q,\infty}^{n/p}\cap L^\infty)} \right)
\left( \Vert u^{m+1}\Vert_{L_{T_{3}}^{2/(2-\epsilon)}(\mathcal{N}%
_{p,q,1}^{n/p+1-\epsilon})}+\Vert\nabla\pi^{m+1}\Vert_{L_{T_{3}}^{2/(2-%
\epsilon)}(\mathcal{N}_{p,q,1}^{n/p-1-\epsilon})}\right) \\
& \hspace{5.5cm}+\Vert u^{m}\Vert_{\widetilde{L}_{T_{3}}^{\infty}(\mathcal{N}%
_{p,q,1}^{n/p-1})}\Vert u^{m+1}\Vert_{L_{T_{3}}^{2/(2-\epsilon )}(\mathcal{N}%
_{p,q,1}^{n/p+1-\epsilon})}.
\end{align*}
Using the interpolation inequality with $\theta=\epsilon/2$ and $%
n/p+1-\epsilon=\theta(n/p-1)+(1-\theta)(n/p+1)$, the uniform boundedness of $%
\{\nabla\pi^{m}\}_{m\in\mathbb{N}_{0}}$ in $L_{T_{3}}^{2/(2-\epsilon )}(%
\mathcal{N}_{p,q,1}^{n/p-1-\epsilon})$, and the previous uniform estimates,
we obtain that $\{\partial_{t}u^{m}\}_{m\in\mathbb{N}_{0}}$ is uniformly
bounded in $L_{T_{3}}^{2/(2-\epsilon)}(\mathcal{N}_{p,q,1}^{n/p-1-\epsilon})$%
.

\hfill

Before proceeding, given {$\alpha\in(0,1)$ and a Banach space }$X$, recall
that {$C^{\alpha}([0,T];X)$ is the space of all $u\in C([0,T];X)$ such that }%
$\Vert u(s)-u(t)\Vert_{X}\leq C|s-t|^{\alpha},$ {for all $s,t\in\lbrack0,T]$%
, for some constant $C>0$. }

For $t<s$, using H\"{o}lder inequality in time, we can estimate{\
\begin{align}
\frac{\Vert u(s)-u(t)\Vert_{\mathcal{N}_{p,q,1}^{n/p-1-\epsilon}}}{%
|s-t|^{\epsilon/2}}\leq\frac{\displaystyle\int_{t}^{s}\Vert\partial
_{z}u\Vert_{\mathcal{N}_{p,q,1}^{n/p-1-\epsilon}}\;dz}{|s-t|^{\epsilon/2}} &
\leq\frac{(s-t)^{\epsilon/2}\Vert\partial_{t}u\Vert_{L_{T}^{2/(2-\epsilon )}(%
\mathcal{N}_{p,q,1}^{n/p-1-\epsilon})}}{|s-t|^{\epsilon/2}}  \notag \\
& =\Vert\partial_{t}u\Vert_{L_{T}^{2/(2-\epsilon)}(\mathcal{N}%
_{p,q,1}^{n/p-1-\epsilon})}.  \label{est:u-u-aux}
\end{align}
}In view of the uniform boundedness of $\{\partial_{t}u^{m}\}_{m\in \mathbb{N%
}_{0}}$ in $L_{T_{3}}^{2/(2-\epsilon)}(\mathcal{N}_{p,q,1}^{n/p-1-\epsilon})$
and estimate (\ref{est:u-u-aux}), it follows that $\{u^{m}\}_{m\in\mathbb{N}%
_{0}}$ is uniformly bounded in $C^{\epsilon /2}([0,T_{3}];\mathcal{N}%
_{p,q,1}^{n/p-1-\epsilon})$. Analogously, we also can conclude that $%
\{a^{m}\}_{m\in\mathbb{N}_{0}}$ is uniformly bounded in $C^{1/2}([0,T_{3}];%
\mathcal{N}_{p,q,r}^{n/p-1} )$.

Moreover, we have the continuous inclusion $\mathcal{N}_{p,q,r,\mathrm{loc}%
}^{s}\hookrightarrow N_{p,q,r,\mathrm{loc}}^{s},$ for $s\geq0,$ and the
compact embeddings (see, e.g., \cite{Runst-Sickel})%
\begin{equation*}
N_{p,q,r,\mathrm{loc}}^{n/p}\hookrightarrow N_{p,q,r,\mathrm{loc}}^{n/p-1}
\hspace{0.5cm} \text{and} \hspace{0.5cm}N_{p,q,1,\mathrm{loc}%
}^{n/p-1}\hookrightarrow N_{p,q,1,\mathrm{loc}}^{n/p-1-\epsilon}.
\end{equation*}
Then, we can apply Arzel\`{a}-Ascoli Theorem (similarly to \cite%
{Danchin_2000}) and obtain the convergence (up to a subsequence) of $%
\{(a^{m},u^{m},\nabla\pi^{m})\}_{m\in\mathbb{N}_{0}}$ to a triple $%
(a,u,\nabla\pi)$ {in $\mathcal{D}^{\prime}([0,T_{3}]\times\mathbb{R}^{n})$}.

\subsubsection{Global case for $s=n/p-1$ with $n/p>1$}

\label{Subsec-conv3}For the global case ($T=\infty$), the proof of the
Cauchy property follows similarly to subsection \ref{Subsec-conv1} by
assuming a smallness condition on $\Vert a_{0}\Vert_{\mathcal{N}%
_{p,q,r}^{n/p}\cap L^\infty}$ and $\Vert u_{0}\Vert_{\mathcal{N}%
_{p,q,1}^{s}} $. In fact, it is sufficient to consider (\ref{est:a.cauchy})
and (\ref{est:up.cauchy}), without applying interpolation, instead of (\ref%
{est:a.cauchy_2}) and (\ref{est:up.cauchy_2}), respectively. Furthermore,
the proof via a compactness argument developed in subsection \ref%
{Subsec-conv2} also can be extended to the global case. The details are left
to the reader. Thus, we obtain that $\{(a^{m},u^{m},\nabla \pi^{m})\}_{m\in%
\mathbb{N}_{0}}$ converges to a triple $\{(a,u,\nabla \pi)\}_{m\in\mathbb{N}%
_{0}}$ {in $\mathcal{D}^{\prime}(\mathbb{R}^{+}\times\mathbb{R}^{n})$}.

\hfill

\subsection{Existence of solution}

{Let $T\in (0,\infty ]$ be the existence time of the limit $(a,u,\nabla \pi
) $ of the approximate solutions }$\{(a^{m},u^{m},\nabla \pi ^{m})\}_{m\in
\mathbb{N}_{0}}$ (up to a subsequence) {obtained in subsections (\ref%
{Subsec-conv1}), (\ref{Subsec-conv2}), and (\ref{Subsec-conv3}), according
to the respective case. }By the uniform boundedness of $\{(a^{m},u^{m},%
\nabla \pi ^{m})\}_{m\in \mathbb{N}_{0}}$, there exists a subsequence $%
\{(a^{m_{k}},u^{m_{k}},\nabla \pi ^{m_{k}})\}_{m_{k}\in \mathbb{N}_{0}}$ in
such a way that $\{(a^{m_{k}},u^{m_{k}},\nabla \pi ^{m_{k}})\}_{m_{k}\in
\mathbb{N}_{0}}$ converges to $(\widetilde{a},\widetilde{u},\nabla
\widetilde{\pi })\in F_{T}^{s,n/p}$ in a weak sense. By the uniqueness of
the limit in the sense of distributions, it follows that $(a,u,\nabla \pi )=(%
\widetilde{a},\widetilde{u},\nabla \widetilde{\pi })$ and, consequently, $%
(a,u,\nabla \pi )\in F_{T}^{s,n/p}$ (see (\ref{space-aux-1})). Moreover, for
the case $n/p-1<s\leq n/p$, using that
\begin{equation*}
\{(a^{m},u^{m})\}_{m\in \mathbb{N}_{0}}\subset C([0,T);\mathcal{N}%
_{p,q,\infty }^{n/p-\epsilon })\times C([0,T);\mathcal{N}_{p,q,1}^{s-%
\epsilon }),
\end{equation*}%
that $\{(a^{m},u^{m})\}_{m\in \mathbb{N}_{0}}$ converges to $(a,u)$ in $%
\widetilde{L}_{T}^{\infty }(\mathcal{N}_{p,q,\infty }^{n/p-\epsilon })\times
\widetilde{L}_{T}^{\infty }(\mathcal{N}_{p,q,1}^{s-\epsilon })$, and
recalling (\ref{est:l_tilde}), it follows that
\begin{equation}
(a,u)\in C([0,T);\mathcal{N}_{p,q,\infty }^{n/p-\epsilon })\times C([0,T);%
\mathcal{N}_{p,q,1}^{s-\epsilon }),  \label{eq:limit_2}
\end{equation}%
where we recall that the time-continuity at $t=0^{+}$ for the density $a$ is
taken in the $\mathcal{S}^{\prime }$-sense. Similarly, for the case $s=n/p-1$%
, we have
\begin{equation}
(a,u)\in C([0,T);\mathcal{N}_{p,q,r}^{n/p-1})\times C([0,T);\mathcal{N}%
_{p,q,1}^{n/p-1-\epsilon }).  \label{eq:limit_3}
\end{equation}%
Using $(a,u,\nabla \pi )\in F_{T}^{s,n/p}$ and (\ref{eq:limit_2})-%
\eqref{eq:limit_3}, we can pass the limit in (\ref{sist:transporte_m}) and (%
\ref{sist:navier_stokes_m}) to obtain that $(a,u,\nabla \pi )$ is a
solution. Furthermore, by standard regularity argument, it follows from $%
(a,u,\nabla \pi )\in F_{T}^{s,n/p}$ and (\ref{eq:limit_2})-\eqref{eq:limit_3}
that
\begin{equation*}
(a,u)\in C([0,T);\mathcal{N}_{p,q,r}^{n/p})\times C([0,T);\mathcal{N}%
_{p,q,1}^{s}).
\end{equation*}

\subsection{Uniqueness of solution}

In this subsection we show the uniqueness part in Theorem \ref%
{the:navier_stokes_2}. Without loss of generality, we consider only the
critical case $s=n/p-1$ with $n/p>1$. The other cases can be proved in a
similar way. Suppose that
\begin{equation}
(a^{i},u^{i},\nabla\pi^{i})\in \widetilde{C}([0,T);\mathcal{N}%
_{p,q,\infty}^{n/p})\cap L_{T}^{\infty}(L^\infty)\times\widetilde{C}([0,T);%
\mathcal{N}_{p,q,1}^{s})\cap L_{T}^{1}(\mathcal{N}_{p,q,1}^{s+2})\times
L_{T}^{1}(\mathcal{N}_{p,q,1}^{s}),  \label{unique_space}
\end{equation}
for $i=1,2$, are two solutions of (\ref{sist:navier_stokes_1}) with the same
initial data $(a_{0},u_{0})$.

Let $(\delta a,\delta u,\nabla\delta\pi):=(a^{2}-a^{1},u^{2}-u^{1},\nabla
\pi^{2}-\nabla\pi^{1})$ and assume that $(a^{1},u^{1},\nabla\pi^{1})$ is the
solution obtained in the previous subsection. First note that
\begin{equation*}
\left\{
\begin{array}{l}
\partial_{t}\delta a+u^{2}\cdot\nabla\delta a+\delta u\cdot\nabla a^{1}=0,
\\
\delta a(\cdot,0)=0.%
\end{array}
\right.
\end{equation*}
Let $\epsilon\in(0,1)$ such that $n/p-1-\epsilon>0$. Proceeding similarly to
proof of (\ref{est:a.cauchy}) and using $\delta a(\cdot,0)=0,$ for all $%
0<t<T,$ we obtain that
\begin{align*}
\Vert\delta a\Vert_{\widetilde{L}_{t}^{\infty}(\mathcal{N}%
_{p,q,\infty}^{n/p-\epsilon})} & \lesssim\Vert\delta u\Vert_{L_{t}^{1}(%
\mathcal{N}_{p,q,1}^{n/p+1-\epsilon})} \Vert a^{1}\Vert_{\widetilde{L}%
_{t}^{\infty }(\mathcal{N}_{p,q,\infty}^{n/p})} + \int_{0}^{t} \Vert\delta
a(\tau)\Vert _{N_{p,q,\infty}^{n/p-\epsilon}} \Vert u^{2}(\tau)\Vert_{%
\mathcal{N}_{p,q,1}^{n/p+1}}\;d\tau.
\end{align*}
Now, Gr\"{o}nwall inequality yields
\begin{align}
\Vert\delta a\Vert_{\widetilde{L}_{t}^{\infty}(\mathcal{N}%
_{p,q,\infty}^{n/p-\epsilon})} & \leq C\exp\left( C\int_{0}^{t}\Vert
u^{2}(\tau)\Vert_{\mathcal{N}_{p,q,1}^{n/p+1}}\;d\tau\right) \Vert\delta
u\Vert_{L_{t}^{1}(\mathcal{N}_{p,q,1}^{n/p+1-\epsilon})} \Vert a^{1}\Vert_{%
\widetilde{L}_{t}^{\infty }(\mathcal{N}_{p,q,\infty}^{n/p})},
\label{est:a.unic}
\end{align}
In turn, we have that the pair $(\delta u,\nabla\delta\pi)$ satisfies
\begin{equation*}
\left\{
\begin{array}{l}
\partial_{t}\delta u-(1+a^{1})(\Delta\delta u-\nabla\delta\pi)=\delta
a(\Delta u^{2}-\nabla\pi^{2})-u^{2}\cdot\nabla\delta u-\delta u\cdot\nabla
u^{1}, \\
\mathrm{div}\;\delta u=0, \\
\delta u(\cdot,0)=0.%
\end{array}
\right.
\end{equation*}
Using the same arguments for obtaining (\ref{est:up.cauchy}), we arrive at
\begin{align}
\Vert\delta u\Vert_{\widetilde{L}_{t}^{\infty}(\mathcal{N}%
_{p,q,1}^{s-\epsilon})}+\Vert\delta u\Vert_{L_{t}^{1}(\mathcal{N}%
_{p,q,1}^{s+2-\epsilon})}+\Vert\nabla\delta\pi\Vert_{L_{t}^{1}(N_{p,q,1}^{s-%
\epsilon })} &  \notag \\
& \hspace{-3cm}\leq C_{1}\left[\Vert a^{1}\Vert_{\widetilde{L}_{t}^{\infty }(%
\mathcal{N}_{p,q,\infty}^{n/p}\cap L^\infty)} \left( \Vert\delta
u\Vert_{L_{t}^{1}(\mathcal{N}_{p,q,1}^{s+2-\epsilon})}+\Vert\nabla\delta\pi%
\Vert_{L_{t}^{1}(\mathcal{N}_{p,q,1}^{s-\epsilon})}\right) \right.  \notag \\
& \hspace{-1.75cm}+\int_{0}^{t} \Vert\delta a(\tau)\Vert_{\mathcal{N}%
_{p,q,\infty}^{n/p-\epsilon}} \left( \Vert u^{2}(\tau)\Vert_{\mathcal{N}%
_{p,q,1}^{s+2}}+\Vert\nabla\pi^{2}(\tau)\Vert_{\mathcal{N}%
_{p,q,1}^{s}}\right) d\tau  \label{est:up.unic} \\
& \hspace{-0.75cm}\left. +\int_{0}^{t}\Vert\delta u(\tau)\Vert _{\mathcal{N}%
_{p,q,1}^{s-\epsilon}}\left( \Vert u^{1}(\tau)\Vert _{\mathcal{N}%
_{p,q,1}^{n/p+1}}+\Vert u^{2}(\tau)\Vert_{\mathcal{N}_{p,q,1}^{n/p+1}}%
\right) d\tau\right] ,  \notag
\end{align}
since $\delta u(\cdot,0)=0$. Then, using $C_{1} \Vert a^{1}\Vert_{%
\widetilde {L}_{t}^{\infty}(\mathcal{N}_{p,q,\infty}^{n/p}\cap L^\infty)}
\leq1/2$ in (\ref{est:up.unic}), it follows that
\begin{align*}
\Vert\delta u\Vert_{\widetilde{L}_{t}^{\infty}(\mathcal{N}%
_{p,q,1}^{s-\epsilon})}+\Vert\delta u\Vert_{L_{t}^{1}(\mathcal{N}%
_{p,q,1}^{s+2-\epsilon})} & \leq C_{3}\left[ \int_{0}^{t} \Vert\delta a(\tau
)\Vert_{\mathcal{N}_{p,q,\infty}^{n/p-\epsilon})} \left( \Vert u^{2}(\tau
)\Vert_{\mathcal{N}_{p,q,1}^{s+2}}+\Vert\nabla\pi^{2}(\tau)\Vert _{\mathcal{N%
}_{p,q,1}^{s}}\right) d\tau\right. \\
& \left. \hspace{2.35cm}+\int_{0}^{t}\Vert\delta u(\tau)\Vert_{\mathcal{N}%
_{p,q,1}^{s-\epsilon}}\left( \Vert u^{1}(\tau)\Vert_{\mathcal{N}%
_{p,q,1}^{n/p+1}}+\Vert u^{2}(\tau)\Vert_{\mathcal{N}_{p,q,1}^{n/p+1}}%
\right) d\tau\right] ,
\end{align*}
for $C_{3}=2C_{1}$. From (\ref{unique_space}) and (\ref{est:a.unic}), we
have that $\Vert\delta a\Vert_{\widetilde{L}_{t}^{\infty}(\mathcal{N}%
_{p,q,\infty}^{n/p-\epsilon})} \leq C_{4}\Vert\delta u\Vert_{L_{t}^{1}(%
\mathcal{N}_{p,q,1}^{n/p+1-\epsilon})}$. Thus, since $n/p-1-\epsilon>0$,
using the above inequality with $s=n/p-1$, we arrive at
\begin{align*}
\Vert\delta u\Vert_{\widetilde{L}_{t}^{\infty}(\mathcal{N}%
_{p,q,1}^{s-\epsilon})}+\Vert\delta u\Vert_{L_{t}^{1}(\mathcal{N}%
_{p,q,1}^{s+2-\epsilon})} & \leq C_{3}\left[ \int_{0}^{t}\Vert\delta u\Vert
_{L_{\tau}^{1}(\mathcal{N}_{p,q,1}^{s+2-\epsilon})}\left( \Vert
u^{2}(\tau)\Vert_{\mathcal{N}_{p,q,1}^{n/p+1}}+\Vert\nabla\pi^{2}(\tau
)\Vert_{\mathcal{N}_{p,q,1}^{n/p-1}}\right) d\tau\right. \\
& \left. \hspace{2.3cm}+\int_{0}^{t}\Vert\delta u(\tau)\Vert_{\mathcal{N}%
_{p,q,1}^{s-\epsilon}}\left( \Vert u^{1}(\tau)\Vert_{\mathcal{N}%
_{p,q,1}^{n/p+1}}+\Vert u^{2}(\tau)\Vert_{\mathcal{N}_{p,q,1}^{n/p+1}}%
\right) d\tau\right] \\
& \hspace{-4.45cm}\leq C_{3}\int_{0}^{t}\left( \Vert\delta u(\tau )\Vert_{%
\mathcal{N}_{p,q,1}^{s-\epsilon}}+\Vert\delta u\Vert_{L_{\tau}^{1}(\mathcal{N%
}_{p,q,1}^{s+2-\epsilon})}\right) \left( \Vert u^{1}(\tau)\Vert_{\mathcal{N}%
_{p,q,1}^{n/p+1}}+\Vert u^{2}(\tau)\Vert _{\mathcal{N}_{p,q,1}^{n/p+1}}+%
\Vert\nabla\pi^{2}(\tau)\Vert_{\mathcal{N}_{p,q,1}^{n/p-1}}\right) d\tau.
\end{align*}

Thus, by (\ref{unique_space}) and Gr\"{o}nwall inequality, we conclude that $%
\delta u=0$. Using (\ref{est:a.unic}) and (\ref{est:up.unic}) for the
pressure, we also obtain that $\delta a=0$ and $\nabla\delta\pi=0,$ for all $%
0<t<T$. Therefore $(a^{1},u^{1},\nabla\pi^{1})=(a^{2},u^{2},\nabla\pi^{2})$
and we are done.

\noindent\textbf{Conflict of interest statement:} The authors declare that
they have no conflict of interest. All co-authors have seen and agree with
the contents of the manuscript and there is no financial interest to report.

\

\noindent\textbf{Data availability statement:} This manuscript has no
associated data.

\end{document}